\input amstex

%Sekretariats-Stildatei
%
%AMSPPT.STY modifiziert und erweitert von C. Krattenthaler

\def\next{AMS-SEKR}\ifx\styname\next \endinput\fi
\catcode`\@=11
\def\styname{AMS-SEKR}
\def\styversion{2.0}
{\W@{}\W@{\styname.STY - Version \styversion}\W@{}}
\hyphenation{acad-e-my acad-e-mies af-ter-thought anom-aly anom-alies
an-ti-deriv-a-tive an-tin-o-my an-tin-o-mies apoth-e-o-ses apoth-e-o-sis
ap-pen-dix ar-che-typ-al as-sign-a-ble as-sist-ant-ship as-ymp-tot-ic
asyn-chro-nous at-trib-uted at-trib-ut-able bank-rupt bank-rupt-cy
bi-dif-fer-en-tial blue-print busier busiest cat-a-stroph-ic
cat-a-stroph-i-cally con-gress cross-hatched data-base de-fin-i-tive
de-riv-a-tive dis-trib-ute dri-ver dri-vers eco-nom-ics econ-o-mist
elit-ist equi-vari-ant ex-quis-ite ex-tra-or-di-nary flow-chart
for-mi-da-ble forth-right friv-o-lous ge-o-des-ic ge-o-det-ic geo-met-ric
griev-ance griev-ous griev-ous-ly hexa-dec-i-mal ho-lo-no-my ho-mo-thetic
ideals idio-syn-crasy in-fin-ite-ly in-fin-i-tes-i-mal ir-rev-o-ca-ble
key-stroke lam-en-ta-ble light-weight mal-a-prop-ism man-u-script
mar-gin-al meta-bol-ic me-tab-o-lism meta-lan-guage me-trop-o-lis
met-ro-pol-i-tan mi-nut-est mol-e-cule mono-chrome mono-pole mo-nop-oly
mono-spline mo-not-o-nous mul-ti-fac-eted mul-ti-plic-able non-euclid-ean
non-iso-mor-phic non-smooth par-a-digm par-a-bol-ic pa-rab-o-loid
pa-ram-e-trize para-mount pen-ta-gon phe-nom-e-non post-script pre-am-ble
pro-ce-dur-al pro-hib-i-tive pro-hib-i-tive-ly pseu-do-dif-fer-en-tial
pseu-do-fi-nite pseu-do-nym qua-drat-ics quad-ra-ture qua-si-smooth
qua-si-sta-tion-ary qua-si-tri-an-gu-lar quin-tes-sence quin-tes-sen-tial
re-arrange-ment rec-tan-gle ret-ri-bu-tion retro-fit retro-fit-ted
right-eous right-eous-ness ro-bot ro-bot-ics sched-ul-ing se-mes-ter
semi-def-i-nite semi-ho-mo-thet-ic set-up se-vere-ly side-step sov-er-eign
spe-cious spher-oid spher-oid-al star-tling star-tling-ly
sta-tis-tics sto-chas-tic straight-est strange-ness strat-a-gem strong-hold
sum-ma-ble symp-to-matic syn-chro-nous topo-graph-i-cal tra-vers-a-ble
tra-ver-sal tra-ver-sals treach-ery turn-around un-at-tached un-err-ing-ly
white-space wide-spread wing-spread wretch-ed wretch-ed-ly Brown-ian
Eng-lish Euler-ian Feb-ru-ary Gauss-ian Grothen-dieck Hamil-ton-ian
Her-mit-ian Jan-u-ary Japan-ese Kor-te-weg Le-gendre Lip-schitz
Lip-schitz-ian Mar-kov-ian Noe-ther-ian No-vem-ber Rie-mann-ian
Schwarz-schild Sep-tem-ber
%Zus"tze (Sind auch in besp.exc einzutragen!):
form per-iods Uni-ver-si-ty cri-ti-sism for-ma-lism}
\Invalid@\nofrills
\Invalid@\usualspace
\newif\ifnofrills@
\def\nofrills@#1#2{\relaxnext@
  \DN@{\ifx\next\nofrills
    \nofrills@true\let#2\relax\DN@\nofrills{\nextii@}%
  \else
    \nofrills@false\def#2{#1}\let\next@\nextii@\fi
\next@}}
\def\usualspace@#1{\ifnofrills@\def\usualspace{#1}\fi}
\def\addto#1#2{\csname \expandafter\eat@\string#1@\endcsname
  \expandafter{\the\csname \expandafter\eat@\string#1@\endcsname#2}}
\newdimen\bigsize@
\def\big@#1#2{{\hbox{$\left#2\vcenter to#1\bigsize@{}%
  \right.\nulldelimiterspace\z@\m@th$}}}
\def\big{\big@\@ne}
\def\Big{\big@{1.5}}
\def\bigg{\big@\tw@}
\def\Bigg{\big@{2.5}}
\def\raggedcenter@{\leftskip\z@ plus.4\hsize \rightskip\leftskip
 \parfillskip\z@ \parindent\z@ \spaceskip.3333em \xspaceskip.5em
 \pretolerance9999\tolerance9999 \exhyphenpenalty\@M
 \hyphenpenalty\@M \let\\\linebreak}
\def\upperspecialchars{\def\ss{SS}\let\i=I\let\j=J\let\ae\AE\let\oe\OE
  \let\o\O\let\aa\AA\let\l\L}
\def\uppercasetext@#1{%
  {\spaceskip1.2\fontdimen2\the\font plus1.2\fontdimen3\the\font
   \upperspecialchars\uctext@#1$\m@th\aftergroup\eat@$}}
\def\uctext@#1$#2${\endash@#1-\endash@$#2$\uctext@}
\def\endash@#1-#2\endash@{\uppercase{#1}\if\notempty{#2}--\endash@#2\endash@\fi}
\def\runaway@#1{\DN@{#1}\ifx\envir@\next@
  \Err@{You seem to have a missing or misspelled \string\end#1 ...}%
  \let\envir@\empty\fi}
\newif\iftemp@
\def\notempty#1{TT\fi\def\test@{#1}\ifx\test@\empty\temp@false
  \else\temp@true\fi \iftemp@}
\font@\tensmc=cmcsc10
\font@\sevenex=cmex7
\font@\sevenit=cmti7
\font@\eightrm=cmr8 % preloaded in plain.tex
\font@\sixrm=cmr6 % preloaded in plain.tex
\font@\eighti=cmmi8     \skewchar\eighti='177 % preloaded
\font@\sixi=cmmi6       \skewchar\sixi='177   % preloaded
\font@\eightsy=cmsy8    \skewchar\eightsy='60 % preloaded
\font@\sixsy=cmsy6      \skewchar\sixsy='60   % preloaded
\font@\eightex=cmex8
\font@\eightbf=cmbx8 % preloaded in plain.tex
\font@\sixbf=cmbx6   % preloaded in plain.tex
\font@\eightit=cmti8 % preloaded in plain.tex
\font@\eightsl=cmsl8 % preloaded in plain.tex
\font@\eightsmc=cmcsc8
\font@\eighttt=cmtt8 % preloaded in plain.tex
%\font@\ninerm=cmr9
%\font@\ninei=cmmi9    \skewchar\ninei='177
%\font@\ninesy=cmsy9   \skewchar\ninesy='60
%\font@\nineex=cmex9
%\font@\ninebf=cmbx9
%\font@\nineit=cmti9
%\font@\ninesl=cmsl9
%\font@\ninesmc=cmcsc9
%\font@\ninemsa=msam9
%\font@\ninemsb=msbm9
%\font@\nineeufm=eufm9

%Erg"nzung des fetten Small-Capitals-Fonts:
%\font@\eightbsmc=cmbcsc10 scaled 833
%\font@\tenbsmc=cmbcsc10
%\font@\twelvebsmc=cmbcsc10 scaled \magstep1
%\font@\fourteenbsmc=cmbcsc10 scaled \magstep2
%\font@\seventeenbsmc=cmbcsc10 scaled \magstep3
%\font@\twentybsmc=cmbcsc10 scaled \magstep4

\loadmsam
\loadmsbm
\loadeufm
\UseAMSsymbols
\newtoks\tenpoint@
\def\tenpoint{\normalbaselineskip12\p@
 \abovedisplayskip12\p@ plus3\p@ minus9\p@
 \belowdisplayskip\abovedisplayskip
 \abovedisplayshortskip\z@ plus3\p@
 \belowdisplayshortskip7\p@ plus3\p@ minus4\p@
 \textonlyfont@\rm\tenrm \textonlyfont@\it\tenit
 \textonlyfont@\sl\tensl \textonlyfont@\bf\tenbf
 \textonlyfont@\smc\tensmc \textonlyfont@\tt\tentt
%Erg"nzung des fetten Small-Capitals-Fonts:
 \textonlyfont@\bsmc\tenbsmc
 \ifsyntax@ \def\big##1{{\hbox{$\left##1\right.$}}}%
  \let\Big\big \let\bigg\big \let\Bigg\big
 \else
  \textfont\z@=\tenrm  \scriptfont\z@=\sevenrm  \scriptscriptfont\z@=\fiverm
  \textfont\@ne=\teni  \scriptfont\@ne=\seveni  \scriptscriptfont\@ne=\fivei
  \textfont\tw@=\tensy \scriptfont\tw@=\sevensy \scriptscriptfont\tw@=\fivesy
  \textfont\thr@@=\tenex \scriptfont\thr@@=\sevenex
        \scriptscriptfont\thr@@=\sevenex
  \textfont\itfam=\tenit \scriptfont\itfam=\sevenit
        \scriptscriptfont\itfam=\sevenit
  \textfont\bffam=\tenbf \scriptfont\bffam=\sevenbf
        \scriptscriptfont\bffam=\fivebf
  \setbox\strutbox\hbox{\vrule height8.5\p@ depth3.5\p@ width\z@}%
  \setbox\strutbox@\hbox{\lower.5\normallineskiplimit\vbox{%
        \kern-\normallineskiplimit\copy\strutbox}}%
 \setbox\z@\vbox{\hbox{$($}\kern\z@}\bigsize@=1.2\ht\z@
 \fi
 \normalbaselines\rm\ex@.2326ex\jot3\ex@\the\tenpoint@}
\newtoks\eightpoint@
\def\eightpoint{\normalbaselineskip10\p@
 \abovedisplayskip10\p@ plus2.4\p@ minus7.2\p@
 \belowdisplayskip\abovedisplayskip
 \abovedisplayshortskip\z@ plus2.4\p@
 \belowdisplayshortskip5.6\p@ plus2.4\p@ minus3.2\p@
 \textonlyfont@\rm\eightrm \textonlyfont@\it\eightit
 \textonlyfont@\sl\eightsl \textonlyfont@\bf\eightbf
 \textonlyfont@\smc\eightsmc \textonlyfont@\tt\eighttt
%Erg"nzung des fetten Small-Capitals-Fonts:
 \textonlyfont@\bsmc\eightbsmc
 \ifsyntax@\def\big##1{{\hbox{$\left##1\right.$}}}%
  \let\Big\big \let\bigg\big \let\Bigg\big
 \else
  \textfont\z@=\eightrm \scriptfont\z@=\sixrm \scriptscriptfont\z@=\fiverm
  \textfont\@ne=\eighti \scriptfont\@ne=\sixi \scriptscriptfont\@ne=\fivei
  \textfont\tw@=\eightsy \scriptfont\tw@=\sixsy \scriptscriptfont\tw@=\fivesy
  \textfont\thr@@=\eightex \scriptfont\thr@@=\sevenex
   \scriptscriptfont\thr@@=\sevenex
  \textfont\itfam=\eightit \scriptfont\itfam=\sevenit
   \scriptscriptfont\itfam=\sevenit
  \textfont\bffam=\eightbf \scriptfont\bffam=\sixbf
   \scriptscriptfont\bffam=\fivebf
 \setbox\strutbox\hbox{\vrule height7\p@ depth3\p@ width\z@}%
 \setbox\strutbox@\hbox{\raise.5\normallineskiplimit\vbox{%
   \kern-\normallineskiplimit\copy\strutbox}}%
 \setbox\z@\vbox{\hbox{$($}\kern\z@}\bigsize@=1.2\ht\z@
 \fi
 \normalbaselines\eightrm\ex@.2326ex\jot3\ex@\the\eightpoint@}

%Definition von 12-point, 14-point und 17-point Fonts
\font@\twelverm=cmr10 scaled\magstep1
\font@\twelveit=cmti10 scaled\magstep1
\font@\twelvesl=cmsl10 scaled\magstep1
\font@\twelvesmc=cmcsc10 scaled\magstep1
\font@\twelvett=cmtt10 scaled\magstep1
\font@\twelvebf=cmbx10 scaled\magstep1
\font@\twelvei=cmmi10 scaled\magstep1
\font@\twelvesy=cmsy10 scaled\magstep1
\font@\twelveex=cmex10 scaled\magstep1
\font@\twelvemsa=msam10 scaled\magstep1
\font@\twelveeufm=eufm10 scaled\magstep1
\font@\twelvemsb=msbm10 scaled\magstep1
\newtoks\twelvepoint@
\def\twelvepoint{\normalbaselineskip15\p@
 \abovedisplayskip15\p@ plus3.6\p@ minus10.8\p@
 \belowdisplayskip\abovedisplayskip
 \abovedisplayshortskip\z@ plus3.6\p@
 \belowdisplayshortskip8.4\p@ plus3.6\p@ minus4.8\p@
 \textonlyfont@\rm\twelverm \textonlyfont@\it\twelveit
 \textonlyfont@\sl\twelvesl \textonlyfont@\bf\twelvebf
 \textonlyfont@\smc\twelvesmc \textonlyfont@\tt\twelvett
%Erg"nzung des fetten Small-Capitals-Fonts:
 \textonlyfont@\bsmc\twelvebsmc
 \ifsyntax@ \def\big##1{{\hbox{$\left##1\right.$}}}%
  \let\Big\big \let\bigg\big \let\Bigg\big
 \else
  \textfont\z@=\twelverm  \scriptfont\z@=\tenrm  \scriptscriptfont\z@=\sevenrm
  \textfont\@ne=\twelvei  \scriptfont\@ne=\teni  \scriptscriptfont\@ne=\seveni
  \textfont\tw@=\twelvesy \scriptfont\tw@=\tensy \scriptscriptfont\tw@=\sevensy
  \textfont\thr@@=\twelveex \scriptfont\thr@@=\tenex
        \scriptscriptfont\thr@@=\tenex
  \textfont\itfam=\twelveit \scriptfont\itfam=\tenit
        \scriptscriptfont\itfam=\tenit
  \textfont\bffam=\twelvebf \scriptfont\bffam=\tenbf
        \scriptscriptfont\bffam=\sevenbf
  \setbox\strutbox\hbox{\vrule height10.2\p@ depth4.2\p@ width\z@}%
  \setbox\strutbox@\hbox{\lower.6\normallineskiplimit\vbox{%
        \kern-\normallineskiplimit\copy\strutbox}}%
 \setbox\z@\vbox{\hbox{$($}\kern\z@}\bigsize@=1.4\ht\z@
 \fi
 \normalbaselines\rm\ex@.2326ex\jot3.6\ex@\the\twelvepoint@}

\def\headfonts{\twelvepoint\bf}

\font@\fourteenrm=cmr10 scaled\magstep2
\font@\fourteenit=cmti10 scaled\magstep2
\font@\fourteensl=cmsl10 scaled\magstep2
\font@\fourteensmc=cmcsc10 scaled\magstep2
\font@\fourteentt=cmtt10 scaled\magstep2
\font@\fourteenbf=cmbx10 scaled\magstep2
\font@\fourteeni=cmmi10 scaled\magstep2
\font@\fourteensy=cmsy10 scaled\magstep2
\font@\fourteenex=cmex10 scaled\magstep2
\font@\fourteenmsa=msam10 scaled\magstep2
\font@\fourteeneufm=eufm10 scaled\magstep2
\font@\fourteenmsb=msbm10 scaled\magstep2
\newtoks\fourteenpoint@
\def\fourteenpoint{\normalbaselineskip15\p@
 \abovedisplayskip18\p@ plus4.3\p@ minus12.9\p@
 \belowdisplayskip\abovedisplayskip
 \abovedisplayshortskip\z@ plus4.3\p@
 \belowdisplayshortskip10.1\p@ plus4.3\p@ minus5.8\p@
 \textonlyfont@\rm\fourteenrm \textonlyfont@\it\fourteenit
 \textonlyfont@\sl\fourteensl \textonlyfont@\bf\fourteenbf
 \textonlyfont@\smc\fourteensmc \textonlyfont@\tt\fourteentt
%Erg"nzung des fetten Small-Capitals-Fonts:
 \textonlyfont@\bsmc\fourteenbsmc
 \ifsyntax@ \def\big##1{{\hbox{$\left##1\right.$}}}%
  \let\Big\big \let\bigg\big \let\Bigg\big
 \else
  \textfont\z@=\fourteenrm  \scriptfont\z@=\twelverm  \scriptscriptfont\z@=\tenrm
  \textfont\@ne=\fourteeni  \scriptfont\@ne=\twelvei  \scriptscriptfont\@ne=\teni
  \textfont\tw@=\fourteensy \scriptfont\tw@=\twelvesy \scriptscriptfont\tw@=\tensy
  \textfont\thr@@=\fourteenex \scriptfont\thr@@=\twelveex
        \scriptscriptfont\thr@@=\twelveex
  \textfont\itfam=\fourteenit \scriptfont\itfam=\twelveit
        \scriptscriptfont\itfam=\twelveit
  \textfont\bffam=\fourteenbf \scriptfont\bffam=\twelvebf
        \scriptscriptfont\bffam=\tenbf
  \setbox\strutbox\hbox{\vrule height12.2\p@ depth5\p@ width\z@}%
  \setbox\strutbox@\hbox{\lower.72\normallineskiplimit\vbox{%
        \kern-\normallineskiplimit\copy\strutbox}}%
 \setbox\z@\vbox{\hbox{$($}\kern\z@}\bigsize@=1.7\ht\z@
 \fi
 \normalbaselines\rm\ex@.2326ex\jot4.3\ex@\the\fourteenpoint@}

\def\chapheadfonts{\fourteenpoint\bf}

\font@\seventeenrm=cmr10 scaled\magstep3
\font@\seventeenit=cmti10 scaled\magstep3
\font@\seventeensl=cmsl10 scaled\magstep3
\font@\seventeensmc=cmcsc10 scaled\magstep3
\font@\seventeentt=cmtt10 scaled\magstep3
\font@\seventeenbf=cmbx10 scaled\magstep3
\font@\seventeeni=cmmi10 scaled\magstep3
\font@\seventeensy=cmsy10 scaled\magstep3
\font@\seventeenex=cmex10 scaled\magstep3
\font@\seventeenmsa=msam10 scaled\magstep3
\font@\seventeeneufm=eufm10 scaled\magstep3
\font@\seventeenmsb=msbm10 scaled\magstep3
\newtoks\seventeenpoint@
\def\seventeenpoint{\normalbaselineskip18\p@
 \abovedisplayskip21.6\p@ plus5.2\p@ minus15.4\p@
 \belowdisplayskip\abovedisplayskip
 \abovedisplayshortskip\z@ plus5.2\p@
 \belowdisplayshortskip12.1\p@ plus5.2\p@ minus7\p@
 \textonlyfont@\rm\seventeenrm \textonlyfont@\it\seventeenit
 \textonlyfont@\sl\seventeensl \textonlyfont@\bf\seventeenbf
 \textonlyfont@\smc\seventeensmc \textonlyfont@\tt\seventeentt
%Erg"nzung des fetten Small-Capitals-Fonts:
 \textonlyfont@\bsmc\seventeenbsmc
 \ifsyntax@ \def\big##1{{\hbox{$\left##1\right.$}}}%
  \let\Big\big \let\bigg\big \let\Bigg\big
 \else
  \textfont\z@=\seventeenrm  \scriptfont\z@=\fourteenrm  \scriptscriptfont\z@=\twelverm
  \textfont\@ne=\seventeeni  \scriptfont\@ne=\fourteeni  \scriptscriptfont\@ne=\twelvei
  \textfont\tw@=\seventeensy \scriptfont\tw@=\fourteensy \scriptscriptfont\tw@=\twelvesy
  \textfont\thr@@=\seventeenex \scriptfont\thr@@=\fourteenex
        \scriptscriptfont\thr@@=\fourteenex
  \textfont\itfam=\seventeenit \scriptfont\itfam=\fourteenit
        \scriptscriptfont\itfam=\fourteenit
  \textfont\bffam=\seventeenbf \scriptfont\bffam=\fourteenbf
        \scriptscriptfont\bffam=\twelvebf
  \setbox\strutbox\hbox{\vrule height14.6\p@ depth6\p@ width\z@}%
  \setbox\strutbox@\hbox{\lower.86\normallineskiplimit\vbox{%
        \kern-\normallineskiplimit\copy\strutbox}}%
 \setbox\z@\vbox{\hbox{$($}\kern\z@}\bigsize@=2\ht\z@
 \fi
 \normalbaselines\rm\ex@.2326ex\jot5.2\ex@\the\seventeenpoint@}

\font@\rrrrrm=cmr10 scaled\magstep4
\font@\bigtitlefont=cmbx10 scaled\magstep4

\parindent1pc
\normallineskiplimit\p@
\newdimen\indenti \indenti=2pc
\def\pageheight#1{\vsize#1}
\def\pagewidth#1{\hsize#1%
   \captionwidth@\hsize \advance\captionwidth@-2\indenti}
\pagewidth{30pc} \pageheight{47pc}
\def\topmatter{%
 \ifx\undefined\msafam
 \else\font@\eightmsa=msam8 \font@\sixmsa=msam6
   \ifsyntax@\else \addto\tenpoint{\textfont\msafam=\tenmsa
              \scriptfont\msafam=\sevenmsa \scriptscriptfont\msafam=\fivemsa}%
     \addto\eightpoint{\textfont\msafam=\eightmsa \scriptfont\msafam=\sixmsa
              \scriptscriptfont\msafam=\fivemsa}%
   \fi
 \fi
 \ifx\undefined\msbfam
 \else\font@\eightmsb=msbm8 \font@\sixmsb=msbm6
   \ifsyntax@\else \addto\tenpoint{\textfont\msbfam=\tenmsb
         \scriptfont\msbfam=\sevenmsb \scriptscriptfont\msbfam=\fivemsb}%
     \addto\eightpoint{\textfont\msbfam=\eightmsb \scriptfont\msbfam=\sixmsb
         \scriptscriptfont\msbfam=\fivemsb}%
   \fi
 \fi
 \ifx\undefined\eufmfam
 \else \font@\eighteufm=eufm8 \font@\sixeufm=eufm6
   \ifsyntax@\else \addto\tenpoint{\textfont\eufmfam=\teneufm
       \scriptfont\eufmfam=\seveneufm \scriptscriptfont\eufmfam=\fiveeufm}%
     \addto\eightpoint{\textfont\eufmfam=\eighteufm
       \scriptfont\eufmfam=\sixeufm \scriptscriptfont\eufmfam=\fiveeufm}%
   \fi
 \fi
 \ifx\undefined\eufbfam
 \else \font@\eighteufb=eufb8 \font@\sixeufb=eufb6
   \ifsyntax@\else \addto\tenpoint{\textfont\eufbfam=\teneufb
      \scriptfont\eufbfam=\seveneufb \scriptscriptfont\eufbfam=\fiveeufb}%
    \addto\eightpoint{\textfont\eufbfam=\eighteufb
      \scriptfont\eufbfam=\sixeufb \scriptscriptfont\eufbfam=\fiveeufb}%
   \fi
 \fi
 \ifx\undefined\eusmfam
 \else \font@\eighteusm=eusm8 \font@\sixeusm=eusm6
   \ifsyntax@\else \addto\tenpoint{\textfont\eusmfam=\teneusm
       \scriptfont\eusmfam=\seveneusm \scriptscriptfont\eusmfam=\fiveeusm}%
     \addto\eightpoint{\textfont\eusmfam=\eighteusm
       \scriptfont\eusmfam=\sixeusm \scriptscriptfont\eusmfam=\fiveeusm}%
   \fi
 \fi
 \ifx\undefined\eusbfam
 \else \font@\eighteusb=eusb8 \font@\sixeusb=eusb6
   \ifsyntax@\else \addto\tenpoint{\textfont\eusbfam=\teneusb
       \scriptfont\eusbfam=\seveneusb \scriptscriptfont\eusbfam=\fiveeusb}%
     \addto\eightpoint{\textfont\eusbfam=\eighteusb
       \scriptfont\eusbfam=\sixeusb \scriptscriptfont\eusbfam=\fiveeusb}%
   \fi
 \fi
 \ifx\undefined\eurmfam
 \else \font@\eighteurm=eurm8 \font@\sixeurm=eurm6
   \ifsyntax@\else \addto\tenpoint{\textfont\eurmfam=\teneurm
       \scriptfont\eurmfam=\seveneurm \scriptscriptfont\eurmfam=\fiveeurm}%
     \addto\eightpoint{\textfont\eurmfam=\eighteurm
       \scriptfont\eurmfam=\sixeurm \scriptscriptfont\eurmfam=\fiveeurm}%
   \fi
 \fi
 \ifx\undefined\eurbfam
 \else \font@\eighteurb=eurb8 \font@\sixeurb=eurb6
   \ifsyntax@\else \addto\tenpoint{\textfont\eurbfam=\teneurb
       \scriptfont\eurbfam=\seveneurb \scriptscriptfont\eurbfam=\fiveeurb}%
    \addto\eightpoint{\textfont\eurbfam=\eighteurb
       \scriptfont\eurbfam=\sixeurb \scriptscriptfont\eurbfam=\fiveeurb}%
   \fi
 \fi
 \ifx\undefined\cmmibfam
 \else \font@\eightcmmib=cmmib8 \font@\sixcmmib=cmmib6
   \ifsyntax@\else \addto\tenpoint{\textfont\cmmibfam=\tencmmib
       \scriptfont\cmmibfam=\sevencmmib \scriptscriptfont\cmmibfam=\fivecmmib}%
    \addto\eightpoint{\textfont\cmmibfam=\eightcmmib
       \scriptfont\cmmibfam=\sixcmmib \scriptscriptfont\cmmibfam=\fivecmmib}%
   \fi
 \fi
 \ifx\undefined\cmbsyfam
 \else \font@\eightcmbsy=cmbsy8 \font@\sixcmbsy=cmbsy6
   \ifsyntax@\else \addto\tenpoint{\textfont\cmbsyfam=\tencmbsy
      \scriptfont\cmbsyfam=\sevencmbsy \scriptscriptfont\cmbsyfam=\fivecmbsy}%
    \addto\eightpoint{\textfont\cmbsyfam=\eightcmbsy
      \scriptfont\cmbsyfam=\sixcmbsy \scriptscriptfont\cmbsyfam=\fivecmbsy}%
   \fi
 \fi
 \let\topmatter\relax}
\def\chapterno@{\uppercase\expandafter{\romannumeral\chaptercount@}}
\newcount\chaptercount@
\def\chapter{\nofrills@{\afterassignment\chapterno@
                        CHAPTER \global\chaptercount@=}\chapter@
 \DNii@##1{\leavevmode\hskip-\leftskip
   \rlap{\vbox to\z@{\vss\centerline{\eightpoint
   \chapter@##1\unskip}\baselineskip2pc\null}}\hskip\leftskip
   \nofrills@false}%
 \FN@\next@}
\newbox\titlebox@

%Uppercase ist abgestellt.
\def\title{\nofrills@{\relax}\title@%
 \DNii@##1\endtitle{\global\setbox\titlebox@\vtop{\tenpoint\bf
 \raggedcenter@\ignorespaces
 \baselineskip1.3\baselineskip\title@{##1}\endgraf}%
 \ifmonograph@ \edef\next{\the\leftheadtoks}\ifx\next\empty
    \leftheadtext{##1}\fi
 \fi
 \edef\next{\the\rightheadtoks}\ifx\next\empty \rightheadtext{##1}\fi
 }\FN@\next@}
\newbox\authorbox@
\def\author#1\endauthor{\global\setbox\authorbox@
 \vbox{\tenpoint\smc\raggedcenter@\ignorespaces
 #1\endgraf}\relaxnext@ \edef\next{\the\leftheadtoks}%
 \ifx\next\empty\leftheadtext{#1}\fi}
\newbox\affilbox@
\def\affil#1\endaffil{\global\setbox\affilbox@
 \vbox{\tenpoint\raggedcenter@\ignorespaces#1\endgraf}}
\newcount\addresscount@
\addresscount@\z@
\def\address#1\endaddress{\global\advance\addresscount@\@ne
  \expandafter\gdef\csname address\number\addresscount@\endcsname
  {\vskip12\p@ minus6\p@\noindent\eightpoint\smc\ignorespaces#1\par}}
\def\email{\nofrills@{\eightpoint{\it E-mail\/}:\enspace}\email@
  \DNii@##1\endemail{%
  \expandafter\gdef\csname email\number\addresscount@\endcsname
  {\def\usualspace{{\it\enspace}}\smallskip\noindent\eightpoint\email@
  \ignorespaces##1\par}}%
 \FN@\next@}
\def\thedate@{}
\def\date#1\enddate{\gdef\thedate@{\tenpoint\ignorespaces#1\unskip}}
\def\thethanks@{}
\def\thanks#1\endthanks{\gdef\thethanks@{\eightpoint\ignorespaces#1.\unskip}}
\def\thekeywords@{}
\def\keywords{\nofrills@{{\it Key words and phrases.\enspace}}\keywords@
 \DNii@##1\endkeywords{\def\thekeywords@{\def\usualspace{{\it\enspace}}%
 \eightpoint\keywords@\ignorespaces##1\unskip.}}%
 \FN@\next@}
\def\thesubjclass@{}
\def\subjclass{\nofrills@{{\rm2010 {\it Mathematics Subject
   Classification\/}.\enspace}}\subjclass@
 \DNii@##1\endsubjclass{\def\thesubjclass@{\def\usualspace
  {{\rm\enspace}}\eightpoint\subjclass@\ignorespaces##1\unskip.}}%
 \FN@\next@}
\newbox\abstractbox@
\def\abstract{\nofrills@{{\smc Abstract.\enspace}}\abstract@
 \DNii@{\setbox\abstractbox@\vbox\bgroup\noindent$$\vbox\bgroup
  \def\envir@{abstract}\advance\hsize-2\indenti
  \usualspace@{{\enspace}}\eightpoint \noindent\abstract@\ignorespaces}%
 \FN@\next@}
\def\endabstract{\par\unskip\egroup$$\egroup}
\def\widestnumber#1#2{\begingroup\let\head\null\let\subhead\empty
   \let\subsubhead\subhead
   \ifx#1\head\global\setbox\tocheadbox@\hbox{#2.\enspace}%
   \else\ifx#1\subhead\global\setbox\tocsubheadbox@\hbox{#2.\enspace}%
   \else\ifx#1\key\bgroup\let\endrefitem@\egroup
        \key#2\endrefitem@\global\refindentwd\wd\keybox@
   \else\ifx#1\no\bgroup\let\endrefitem@\egroup
        \no#2\endrefitem@\global\refindentwd\wd\nobox@
   \else\ifx#1\page\global\setbox\pagesbox@\hbox{\quad\bf#2}%
   \else\ifx#1\item\setboxz@h{#2}\global\rosteritemwd\wdz@
        \global\advance\rosteritemwd by.5\parindent
   \else\message{\string\widestnumber is not defined for this option
   (\string#1)}%
\fi\fi\fi\fi\fi\fi\endgroup}
\newif\ifmonograph@
\def\Monograph{\monograph@true \let\headmark\rightheadtext
  \let\varindent@\indent \def\headfont@{\bf}\def\proclaimheadfont@{\smc}%
  \def\demofont@{\smc}}
%Bei Proclaim,...: Einrcken.
\let\varindent@\indent

\newbox\tocheadbox@    \newbox\tocsubheadbox@
\newbox\tocbox@
\def\toc{\toc@{Contents}}
\def\newtocdefs{%
   \def \title##1\endtitle
       {\penaltyandskip@\z@\smallskipamount
        \hangindent\wd\tocheadbox@\noindent{\bf##1}}%
   \def \chapter##1{%
        Chapter \uppercase\expandafter{\romannumeral##1.\unskip}\enspace}%
   \def \specialhead##1\endspecialhead
       {\par\hangindent\wd\tocheadbox@ \noindent##1\par}%
   \def \head##1 ##2\endhead
       {\par\hangindent\wd\tocheadbox@ \noindent
        \if\notempty{##1}\hbox to\wd\tocheadbox@{\hfil##1\enspace}\fi
        ##2\par}%
   \def \subhead##1 ##2\endsubhead
       {\par\vskip-\parskip {\normalbaselines
        \advance\leftskip\wd\tocheadbox@
        \hangindent\wd\tocsubheadbox@ \noindent
        \if\notempty{##1}\hbox to\wd\tocsubheadbox@{##1\unskip\hfil}\fi
         ##2\par}}%
   \def \subsubhead##1 ##2\endsubsubhead
       {\par\vskip-\parskip {\normalbaselines
        \advance\leftskip\wd\tocheadbox@
        \hangindent\wd\tocsubheadbox@ \noindent
        \if\notempty{##1}\hbox to\wd\tocsubheadbox@{##1\unskip\hfil}\fi
        ##2\par}}}
\def\toc@#1{\relaxnext@
   \def\page##1%
       {\unskip\penalty0\null\hfil
        \rlap{\hbox to\wd\pagesbox@{\quad\hfil##1}}\hfilneg\penalty\@M}%
 \DN@{\ifx\next\nofrills\DN@\nofrills{\nextii@}%
      \else\DN@{\nextii@{{#1}}}\fi
      \next@}%
 \DNii@##1{%
\ifmonograph@\bgroup\else\setbox\tocbox@\vbox\bgroup
   \centerline{\headfont@\ignorespaces##1\unskip}\nobreak
   \vskip\belowheadskip \fi
   \setbox\tocheadbox@\hbox{0.\enspace}%
   \setbox\tocsubheadbox@\hbox{0.0.\enspace}%
   \leftskip\indenti \rightskip\leftskip
   \setbox\pagesbox@\hbox{\bf\quad000}\advance\rightskip\wd\pagesbox@
   \newtocdefs
 }%
 \FN@\next@}
\def\endtoc{\par\egroup}
\let\pretitle\relax
\let\preauthor\relax
\let\preaffil\relax
\let\predate\relax
\let\preabstract\relax
\let\prepaper\relax
\def\dedicatory #1\enddedicatory{\def\preabstract{{\medskip
  \eightpoint\it \raggedcenter@#1\endgraf}}}
\def\thetranslator@{}
\def\translator#1\endtranslator{\def\thetranslator@{\nobreak\medskip
 \line{\eightpoint\hfil Translated by \uppercase{#1}\qquad\qquad}\nobreak}}
\outer\def\endtopmatter{\runaway@{abstract}%
 \edef\next{\the\leftheadtoks}\ifx\next\empty
  \expandafter\leftheadtext\expandafter{\the\rightheadtoks}\fi
 \ifmonograph@\else
   \ifx\thesubjclass@\empty\else \makefootnote@{}{\thesubjclass@}\fi
   \ifx\thekeywords@\empty\else \makefootnote@{}{\thekeywords@}\fi
   \ifx\thethanks@\empty\else \makefootnote@{}{\thethanks@}\fi
 \fi
  \pretitle
  \ifmonograph@ \topskip7pc \else \topskip4pc \fi
  \box\titlebox@
  \topskip10pt% reset to normal value
  \preauthor
  \ifvoid\authorbox@\else \vskip2.5pc plus1pc \unvbox\authorbox@\fi
  \preaffil
  \ifvoid\affilbox@\else \vskip1pc plus.5pc \unvbox\affilbox@\fi
  \predate
  \ifx\thedate@\empty\else \vskip1pc plus.5pc \line{\hfil\thedate@\hfil}\fi
  \preabstract
  \ifvoid\abstractbox@\else \vskip1.5pc plus.5pc \unvbox\abstractbox@ \fi
  \ifvoid\tocbox@\else\vskip1.5pc plus.5pc \unvbox\tocbox@\fi
  \prepaper
  \vskip2pc plus1pc
}
\def\document{\let\fontlist@\relax\let\alloclist@\relax
  \tenpoint}

%Modifizierte Head-Skips
\newskip\aboveheadskip       \aboveheadskip1.8\bigskipamount
\newdimen\belowheadskip      \belowheadskip1.8\medskipamount

\def\headfont@{\smc}
\def\penaltyandskip@#1#2{\relax\ifdim\lastskip<#2\relax\removelastskip
      \ifnum#1=\z@\else\penalty@#1\relax\fi\vskip#2%
  \else\ifnum#1=\z@\else\penalty@#1\relax\fi\fi}
\def\nobreak{\penalty\@M
  \ifvmode\def\penalty@{\let\penalty@\penalty\count@@@}%
  \everypar{\let\penalty@\penalty\everypar{}}\fi}
\let\penalty@\penalty
\def\heading#1\endheading{\head#1\endhead}

\def\specialheadfont@{\bf}
\outer\def\specialhead{\par\penaltyandskip@{-200}\aboveheadskip
  \begingroup\interlinepenalty\@M\rightskip\z@ plus\hsize \let\\\linebreak
  \specialheadfont@\noindent\ignorespaces}
\def\endspecialhead{\par\endgroup\nobreak\vskip\belowheadskip}
%\outer\def\head#1\endhead{\par\penaltyandskip@{-200}\aboveheadskip
% {\headfont@\raggedcenter@\interlinepenalty\@M
% \ignorespaces#1\endgraf}\nobreak
% \vskip\belowheadskip
% \headmark{#1}}
\let\headmark\eat@
\newskip\subheadskip       \subheadskip\medskipamount
\def\subheadfont@{\bf}
\outer\def\subhead{\nofrills@{.\enspace}\subhead@
 \DNii@##1\endsubhead{\par\penaltyandskip@{-100}\subheadskip
  \varindent@{\usualspace@{{\subheadfont@\enspace}}%
 \subheadfont@\ignorespaces##1\unskip\subhead@}\ignorespaces}%
 \FN@\next@}
\outer\def\subsubhead{\nofrills@{.\enspace}\subsubhead@
 \DNii@##1\endsubsubhead{\par\penaltyandskip@{-50}\medskipamount
      {\usualspace@{{\it\enspace}}%
  \it\ignorespaces##1\unskip\subsubhead@}\ignorespaces}%
 \FN@\next@}
\def\proclaimheadfont@{\bf}
\outer\def\proclaim{\runaway@{proclaim}\def\envir@{proclaim}%
  \nofrills@{.\enspace}\proclaim@
 \DNii@##1{\penaltyandskip@{-100}\medskipamount\varindent@
   \usualspace@{{\proclaimheadfont@\enspace}}\proclaimheadfont@
   \ignorespaces##1\unskip\proclaim@
  \sl\ignorespaces}% 
 \FN@\next@}
\outer\def\endproclaim{\let\envir@\relax\par\rm
  \penaltyandskip@{55}\medskipamount}
\def\demoheadfont@{\it}
\def\demo{\runaway@{proclaim}\nofrills@{.\enspace}\demo@
     \DNii@##1{\par\penaltyandskip@\z@\medskipamount
  {\usualspace@{{\demoheadfont@\enspace}}%
  \varindent@\demoheadfont@\ignorespaces##1\unskip\demo@}\rm
  \ignorespaces}\FN@\next@}
\def\enddemo{\par\medskip}
\def\qed{\ifhmode\unskip\nobreak\fi\quad\ifmmode\square\else$\m@th\square$\fi}
\let\remark\demo
%\endremark=\enddemo
\let\endremark\enddemo
%Headfont fr Definition wie bei Beweisen.
\def\definition{\runaway@{proclaim}%
  \nofrills@{.\demoheadfont@\enspace}\definition@
        \DNii@##1{\penaltyandskip@{-100}\medskipamount
        {\usualspace@{{\demoheadfont@\enspace}}%
        \varindent@\demoheadfont@\ignorespaces##1\unskip\definition@}%
        \rm \ignorespaces}\FN@\next@}

%\let\example\definition
%\let\endexample\enddefinition
%Modifikation:

\newdimen\rosteritemwd
\newcount\rostercount@
\newif\iffirstitem@
\let\plainitem@\item
\newtoks\everypartoks@
\def\par@{\everypartoks@\expandafter{\the\everypar}\everypar{}}
\def\roster{\edef\leftskip@{\leftskip\the\leftskip}%
 \relaxnext@
 \rostercount@\z@  
 \def\item{\FN@\rosteritem@}% 
 \DN@{\ifx\next\runinitem\let\next@\nextii@\else
  \let\next@\nextiii@\fi\next@}%
 \DNii@\runinitem  
  {\unskip  
   \DN@{\ifx\next[\let\next@\nextii@\else
    \ifx\next"\let\next@\nextiii@\else\let\next@\nextiv@\fi\fi\next@}%
   \DNii@[####1]{\rostercount@####1\relax
    \enspace{\rm(\number\rostercount@)}~\ignorespaces}%
   \def\nextiii@"####1"{\enspace{\rm####1}~\ignorespaces}%
   \def\nextiv@{\enspace{\rm(1)}\rostercount@\@ne~}%
   \par@\firstitem@false  
   \FN@\next@}% 
 \def\nextiii@{\par\par@  
  \penalty\@m\smallskip\vskip-\parskip
  \firstitem@true}%  
 \FN@\next@}
\def\rosteritem@{\iffirstitem@\firstitem@false\else\par\vskip-\parskip\fi
 \leftskip3\parindent\noindent  
 \DNii@[##1]{\rostercount@##1\relax
  \llap{\hbox to2.5\parindent{\hss\rm(\number\rostercount@)}%
   \hskip.5\parindent}\ignorespaces}%
 \def\nextiii@"##1"{%
  \llap{\hbox to2.5\parindent{\hss\rm##1}\hskip.5\parindent}\ignorespaces}%
 \def\nextiv@{\advance\rostercount@\@ne
  \llap{\hbox to2.5\parindent{\hss\rm(\number\rostercount@)}%
   \hskip.5\parindent}}%
 \ifx\next[\let\next@\nextii@\else\ifx\next"\let\next@\nextiii@\else
  \let\next@\nextiv@\fi\fi\next@}

\newif\ifnextRunin@
\def\endroster{\relaxnext@
 \par\leftskip@  
 \penalty-50 \vskip-\parskip\smallskip  
 \DN@{\ifx\next\Runinitem\let\next@\relax
  \else\nextRunin@false\let\item\plainitem@  
   \ifx\next\par 
    \DN@\par{\everypar\expandafter{\the\everypartoks@}}%
   \else  
    \DN@{\noindent\everypar\expandafter{\the\everypartoks@}}%
  \fi\fi\next@}%
 \FN@\next@}
\newcount\rosterhangafter@
\def\Runinitem#1\roster\runinitem{\relaxnext@
 \rostercount@\z@ 
 \def\item{\FN@\rosteritem@}%  
 \def\runinitem@{#1}%  
 \DN@{\ifx\next[\let\next\nextii@\else\ifx\next"\let\next\nextiii@
  \else\let\next\nextiv@\fi\fi\next}%
 \DNii@[##1]{\rostercount@##1\relax
  \def\item@{{\rm(\number\rostercount@)}}\nextv@}%
 \def\nextiii@"##1"{\def\item@{{\rm##1}}\nextv@}%
 \def\nextiv@{\advance\rostercount@\@ne
  \def\item@{{\rm(\number\rostercount@)}}\nextv@}%
 \def\nextv@{\setbox\z@\vbox  
  {\ifnextRunin@\noindent\fi  
  \runinitem@\unskip\enspace\item@~\par  
  \global\rosterhangafter@\prevgraf}% 
  \firstitem@false  
  \ifnextRunin@\else\par\fi  
  \hangafter\rosterhangafter@\hangindent3\parindent
  \ifnextRunin@\noindent\fi  
  \runinitem@\unskip\enspace 
  \item@~\ifnextRunin@\else\par@\fi  
  \nextRunin@true\ignorespaces}%  
 \FN@\next@}
\def\footmarkform@#1{$\m@th^{#1}$}
\let\thefootnotemark\footmarkform@
\def\makefootnote@#1#2{\insert\footins
 {\interlinepenalty\interfootnotelinepenalty
 \eightpoint\splittopskip\ht\strutbox\splitmaxdepth\dp\strutbox
 \floatingpenalty\@MM\leftskip\z@\rightskip\z@\spaceskip\z@\xspaceskip\z@
 \leavevmode{#1}\footstrut\ignorespaces#2\unskip\lower\dp\strutbox
 \vbox to\dp\strutbox{}}}
\newcount\footmarkcount@
\footmarkcount@\z@
\def\footnotemark{\let\@sf\empty\relaxnext@
 \ifhmode\edef\@sf{\spacefactor\the\spacefactor}\/\fi
 \DN@{\ifx[\next\let\next@\nextii@\else
  \ifx"\next\let\next@\nextiii@\else
  \let\next@\nextiv@\fi\fi\next@}%
 \DNii@[##1]{\footmarkform@{##1}\@sf}%
 \def\nextiii@"##1"{{##1}\@sf}%
 \def\nextiv@{\iffirstchoice@\global\advance\footmarkcount@\@ne\fi
  \footmarkform@{\number\footmarkcount@}\@sf}%
 \FN@\next@}
\def\footnotetext{\relaxnext@
 \DN@{\ifx[\next\let\next@\nextii@\else
  \ifx"\next\let\next@\nextiii@\else
  \let\next@\nextiv@\fi\fi\next@}%
 \DNii@[##1]##2{\makefootnote@{\footmarkform@{##1}}{##2}}%
 \def\nextiii@"##1"##2{\makefootnote@{##1}{##2}}%
 \def\nextiv@##1{\makefootnote@{\footmarkform@{\number\footmarkcount@}}{##1}}%
 \FN@\next@}
\def\footnote{\let\@sf\empty\relaxnext@
 \ifhmode\edef\@sf{\spacefactor\the\spacefactor}\/\fi
 \DN@{\ifx[\next\let\next@\nextii@\else
  \ifx"\next\let\next@\nextiii@\else
  \let\next@\nextiv@\fi\fi\next@}%
 \DNii@[##1]##2{\footnotemark[##1]\footnotetext[##1]{##2}}%
 \def\nextiii@"##1"##2{\footnotemark"##1"\footnotetext"##1"{##2}}%
 \def\nextiv@##1{\footnotemark\footnotetext{##1}}%
 \FN@\next@}
\def\adjustfootnotemark#1{\advance\footmarkcount@#1\relax}
\def\footnoterule{\kern-3\p@
  \hrule width 5pc\kern 2.6\p@} 
\def\captionfont@{\smc}
\def\topcaption#1#2\endcaption{%
  {\dimen@\hsize \advance\dimen@-\captionwidth@
   \rm\raggedcenter@ \advance\leftskip.5\dimen@ \rightskip\leftskip
  {\captionfont@#1}%
  \if\notempty{#2}.\enspace\ignorespaces#2\fi
  \endgraf}\nobreak\bigskip}
\def\botcaption#1#2\endcaption{%
  \nobreak\bigskip
  \setboxz@h{\captionfont@#1\if\notempty{#2}.\enspace\rm#2\fi}%
  {\dimen@\hsize \advance\dimen@-\captionwidth@
   \leftskip.5\dimen@ \rightskip\leftskip
   \noindent \ifdim\wdz@>\captionwidth@ 
   \else\hfil\fi 
  {\captionfont@#1}\if\notempty{#2}.\enspace\rm#2\fi\endgraf}}
\def\@ins{\par\begingroup\def\vspace##1{\vskip##1\relax}%
  \def\captionwidth##1{\captionwidth@##1\relax}%
  \setbox\z@\vbox\bgroup} % start a \vbox
\def\block{\RIfMIfI@\nondmatherr@\block\fi
       \else\ifvmode\vskip\abovedisplayskip\noindent\fi
        $$\def\endblock{\par\egroup$$}\fi
  \vbox\bgroup\advance\hsize-2\indenti\noindent}
\def\endblock{\par\egroup}
\def\cite#1{{\rm[{\citefont@\m@th#1}]}}
\def\citefont@{\rm}
\def\refsfont@{\eightpoint}
\outer\def\Refs{\runaway@{proclaim}%
 \relaxnext@ \DN@{\ifx\next\nofrills\DN@\nofrills{\nextii@}\else
  \DN@{\nextii@{References}}\fi\next@}%
 \DNii@##1{\penaltyandskip@{-200}\aboveheadskip
  \line{\hfil\headfont@\ignorespaces##1\unskip\hfil}\nobreak
  \vskip\belowheadskip
  \begingroup\refsfont@\sfcode`.=\@m}%
 \FN@\next@}
\def\endRefs{\par\endgroup}
\newbox\nobox@            \newbox\keybox@           \newbox\bybox@
\newbox\paperbox@         \newbox\paperinfobox@     \newbox\jourbox@
\newbox\volbox@           \newbox\issuebox@         \newbox\yrbox@
\newbox\pagesbox@         \newbox\bookbox@          \newbox\bookinfobox@
\newbox\publbox@          \newbox\publaddrbox@      \newbox\finalinfobox@
\newbox\edsbox@           \newbox\langbox@
\newif\iffirstref@        \newif\iflastref@
\newif\ifprevjour@        \newif\ifbook@            \newif\ifprevinbook@
\newif\ifquotes@          \newif\ifbookquotes@      \newif\ifpaperquotes@
\newdimen\bysamerulewd@
\setboxz@h{\refsfont@\kern3em}
\bysamerulewd@\wdz@
\newdimen\refindentwd
\setboxz@h{\refsfont@ 00. }
\refindentwd\wdz@
\outer\def\ref{\begingroup \noindent\hangindent\refindentwd
 \firstref@true \def\nofrills{\def\refkern@{\kern3sp}}%
 \ref@}
\def\ref@{\book@false \bgroup\let\endrefitem@\egroup \ignorespaces}
\def\moreref{\endrefitem@\endref@\firstref@false\ref@}%
\def\transl{\endrefitem@\endref@\firstref@false
  \book@false
  \prepunct@
  \setboxz@h\bgroup \aftergroup\unhbox\aftergroup\z@
    \def\endrefitem@{\unskip\refkern@\egroup}\ignorespaces}%
\def\emptyifempty@{\dimen@\wd\currbox@
  \advance\dimen@-\wd\z@ \advance\dimen@-.1\p@
  \ifdim\dimen@<\z@ \setbox\currbox@\copy\voidb@x \fi}
\let\refkern@\relax
\def\endrefitem@{\unskip\refkern@\egroup
  \setboxz@h{\refkern@}\emptyifempty@}\ignorespaces
\def\refdef@#1#2#3{\edef\next@{\noexpand\endrefitem@
  \let\noexpand\currbox@\csname\expandafter\eat@\string#1box@\endcsname
    \noexpand\setbox\noexpand\currbox@\hbox\bgroup}%
  \toks@\expandafter{\next@}%
  \if\notempty{#2#3}\toks@\expandafter{\the\toks@
  \def\endrefitem@{\unskip#3\refkern@\egroup
  \setboxz@h{#2#3\refkern@}\emptyifempty@}#2}\fi
  \toks@\expandafter{\the\toks@\ignorespaces}%
  \edef#1{\the\toks@}}
\refdef@\no{}{. }
\refdef@\key{[\m@th}{] }
\refdef@\by{}{}
\def\bysame{\by\hbox to\bysamerulewd@{\hrulefill}\thinspace
   \kern0sp}
\def\manyby{\message{\string\manyby is no longer necessary; \string\by
  can be used instead, starting with version 2.0 of \styname.STY}\by}
\refdef@\paper{\ifpaperquotes@``\fi\it}{}
\refdef@\paperinfo{}{}
\def\jour{\endrefitem@\let\currbox@\jourbox@
  \setbox\currbox@\hbox\bgroup
  \def\endrefitem@{\unskip\refkern@\egroup
    \setboxz@h{\refkern@}\emptyifempty@
    \ifvoid\jourbox@\else\prevjour@true\fi}%
\ignorespaces}
\refdef@\vol{\ifbook@\else\bf\fi}{}
\refdef@\issue{no. }{}
\refdef@\yr{}{}
\refdef@\pages{}{}
\def\page{\endrefitem@\def\pp@{\def\pp@{pp.~}p.~}\let\currbox@\pagesbox@
  \setbox\currbox@\hbox\bgroup\ignorespaces}
\def\pp@{pp.~}
\def\book{\endrefitem@ \let\currbox@\bookbox@
 \setbox\currbox@\hbox\bgroup\def\endrefitem@{\unskip\refkern@\egroup
  \setboxz@h{\ifbookquotes@``\fi}\emptyifempty@
  \ifvoid\bookbox@\else\book@true\fi}%
  \ifbookquotes@``\fi\it\ignorespaces}
\def\inbook{\endrefitem@
  \let\currbox@\bookbox@\setbox\currbox@\hbox\bgroup
  \def\endrefitem@{\unskip\refkern@\egroup
  \setboxz@h{\ifbookquotes@``\fi}\emptyifempty@
  \ifvoid\bookbox@\else\book@true\previnbook@true\fi}%
  \ifbookquotes@``\fi\ignorespaces}
\refdef@\eds{(}{, eds.)}
\def\ed{\endrefitem@\let\currbox@\edsbox@
 \setbox\currbox@\hbox\bgroup
 \def\endrefitem@{\unskip, ed.)\refkern@\egroup
  \setboxz@h{(, ed.)}\emptyifempty@}(\ignorespaces}
\refdef@\bookinfo{}{}
\refdef@\publ{}{}
\refdef@\publaddr{}{}
\refdef@\finalinfo{}{}
\refdef@\lang{(}{)}

\let\refdef@\relax 
\def\ppunbox@#1{\ifvoid#1\else\prepunct@\unhbox#1\fi}
\def\nocomma@#1{\ifvoid#1\else\changepunct@3\prepunct@\unhbox#1\fi}
\def\changepunct@#1{\ifnum\lastkern<3 \unkern\kern#1sp\fi}
\def\prepunct@{\count@\lastkern\unkern
  \ifnum\lastpenalty=0
    \let\penalty@\relax
  \else
    \edef\penalty@{\penalty\the\lastpenalty\relax}%
  \fi
  \unpenalty
  \let\refspace@\ \ifcase\count@,% usual case, do a comma
\or;\or.\or % do nothing; this case is from nofrills.
  \or\let\refspace@\relax
  \else,\fi
  \ifquotes@''\quotes@false\fi \penalty@ \refspace@
}
\def\transferpenalty@#1{\dimen@\lastkern\unkern
  \ifnum\lastpenalty=0\unpenalty\let\penalty@\relax
  \else\edef\penalty@{\penalty\the\lastpenalty\relax}\unpenalty\fi
  #1\penalty@\kern\dimen@}
\def\endref{\endrefitem@\lastref@true\endref@
  \par\endgroup \prevjour@false \previnbook@false }
\def\endref@{%
\iffirstref@
  \ifvoid\nobox@\ifvoid\keybox@\indent\fi
  \else\hbox to\refindentwd{\hss\unhbox\nobox@}\fi
  \ifvoid\keybox@
  \else\ifdim\wd\keybox@>\refindentwd
         \box\keybox@
       \else\hbox to\refindentwd{\unhbox\keybox@\hfil}\fi\fi
  \kern4sp\ppunbox@\bybox@
\fi 
  \ifvoid\paperbox@
  \else\prepunct@\unhbox\paperbox@
    \ifpaperquotes@\quotes@true\fi\fi
  \ppunbox@\paperinfobox@
  \ifvoid\jourbox@
    \ifprevjour@ \nocomma@\volbox@
      \nocomma@\issuebox@
      \ifvoid\yrbox@\else\changepunct@3\prepunct@(\unhbox\yrbox@
        \transferpenalty@)\fi
      \ppunbox@\pagesbox@
    \fi 
  \else \prepunct@\unhbox\jourbox@
    \nocomma@\volbox@
    \nocomma@\issuebox@
    \ifvoid\yrbox@\else\changepunct@3\prepunct@(\unhbox\yrbox@
      \transferpenalty@)\fi
    \ppunbox@\pagesbox@
  \fi 
  \ifbook@\prepunct@\unhbox\bookbox@ \ifbookquotes@\quotes@true\fi \fi
  \nocomma@\edsbox@
  \ppunbox@\bookinfobox@
  \ifbook@\ifvoid\volbox@\else\prepunct@ vol.~\unhbox\volbox@
  \fi\fi
  \ppunbox@\publbox@ \ppunbox@\publaddrbox@
  \ifbook@ \ppunbox@\yrbox@
    \ifvoid\pagesbox@
    \else\prepunct@\pp@\unhbox\pagesbox@\fi
  \else
    \ifprevinbook@ \ppunbox@\yrbox@
      \ifvoid\pagesbox@\else\prepunct@\pp@\unhbox\pagesbox@\fi
    \fi \fi
  \ppunbox@\finalinfobox@
  \iflastref@
    \ifvoid\langbox@.\ifquotes@''\fi
    \else\changepunct@2\prepunct@\unhbox\langbox@\fi
  \else
    \ifvoid\langbox@\changepunct@1%
    \else\changepunct@3\prepunct@\unhbox\langbox@
      \changepunct@1\fi
  \fi
}
\outer\def\enddocument{%
 \runaway@{proclaim}%
\ifmonograph@ % do nothing
\else
 \nobreak
 \thetranslator@
 \count@\z@ \loop\ifnum\count@<\addresscount@\advance\count@\@ne
 \csname address\number\count@\endcsname
 \csname email\number\count@\endcsname
 \repeat
\fi
 \vfill\supereject\end}

%Modifizierte Fonts fr Proclaim,...
\def\headfont@{\headfonts}
\def\proclaimheadfont@{\bf}
\def\specialheadfont@{\bf}
\def\subheadfont@{\bf}
\def\demoheadfont@{\smc}

%Kontrollsequenzen fr Inhaltsverzeichnis und Index:
\newif\ifThisToToc \ThisToTocfalse
\newif\iftocloaded \tocloadedfalse

\def\C@L{\noexpand\Cal}\def\B@B{\noexpand\Bbb}\def\fR@K{\noexpand\frak}
\def\S@{\noexpand\S}\def\P@P{\noexpand\"}
\def\xpar{\\}

\def\writetoc#1{\iftocloaded\ifThisToToc\begingroup\def\totoc{}
  \def\Cal{\noexpand\C@L}\def\Bbb{\noexpand\B@B}
  \def\frak{\noexpand\fR@K}\def\goth{\frak}\def\S{\noexpand\S@}
  \def\"{\noexpand\P@P}
  \def\xpar{\par\penalty100000 }\def\idx##1{##1}\def\\{\xpar}
  \edef\next@{\write\toc{\noindent#1\leaderfill\noexpand\folio\par}}%
  \next@\endgroup\global\ThisToTocfalse\fi\fi}
\def\leaderfill{\leaders\hbox to 1em{\hss.\hss}\hfill}

\newif\ifindexloaded \indexloadedfalse
\def\idx#1{\ifindexloaded\begingroup\def\ign{}\def\it{}\def\/{}%
 \def\smc{}\def\bf{}\def\tt{}%
 \def\Cal{\noexpand\C@L}\def\Bbb{\noexpand\B@B}%
 \def\frak{\noexpand\fR@K}\def\goth{\frak}\def\S{\noexpand\S@}%
  \def\"{\noexpand\P@P}%
 {\edef\next@{\write\index{#1, \noexpand\folio}}\next@}%
 \endgroup\fi{#1}}
\def\ign#1{}

\def\totoc{\global\ThisToToctrue}

\outer\def\head#1\endhead{\par\penaltyandskip@{-200}\aboveheadskip
 {\headfont@\raggedcenter@\interlinepenalty\@M
 \ignorespaces#1\endgraf}\nobreak
 \vskip\belowheadskip
 \headmark{#1}\writetoc{#1}}

\outer\def\chaphead#1\endchaphead{\par\penaltyandskip@{-200}\aboveheadskip
 {\chapheadfonts\raggedcenter@\interlinepenalty\@M
 \ignorespaces#1\endgraf}\nobreak
 \vskip3\belowheadskip
 \headmark{#1}\writetoc{#1}}

\def\folio{{\foliofont@\ifnum\pageno<\z@ \romannumeral-\pageno
 \else\number\pageno \fi}}
\newtoks\leftheadtoks
\newtoks\rightheadtoks

%Uppercase ist abgestellt:
\def\leftheadtext{\nofrills@{\relax}\lht@
  \DNii@##1{\leftheadtoks\expandafter{\lht@{##1}}%
    \mark{\the\leftheadtoks\noexpand\else\the\rightheadtoks}
    \ifsyntax@\setboxz@h{\def\\{\unskip\space\ignorespaces}%
        \headlinefont@##1}\fi}%
  \FN@\next@}
%Uppercase ist abgestellt:
\def\rightheadtext{\nofrills@{\relax}\rht@
  \DNii@##1{\rightheadtoks\expandafter{\rht@{##1}}%
    \mark{\the\leftheadtoks\noexpand\else\the\rightheadtoks}%
    \ifsyntax@\setboxz@h{\def\\{\unskip\space\ignorespaces}%
        \headlinefont@##1}\fi}%
  \FN@\next@}
%\headline={\def\chapter#1{\chapterno@. }%
%  \def\\{\unskip\space\ignorespaces}\headlinefont@
%  \ifodd\pageno \rightheadline \else \leftheadline\fi}
\def\NoRunningHeads{\global\runheads@false\global\let\headmark\eat@}

\newif\iffirstpage@     \firstpage@true
\newif\ifrunheads@      \runheads@true

%Erg"nzungen zu Runningheads und Pagenumbers:
\newdimen\fullhsize \fullhsize=\hsize
\newdimen\fullvsize \fullvsize=\vsize
\def\fullline{\hbox to\fullhsize}

\def\pagenumbers{\gdef\folio{\folio@}}

\let\norunningheads\NoRunningHeads
\def\userunningheads{\global\runheads@true}
%Default: Seitennumerierung unten.
\norunningheads

\headline={\def\chapter#1{\chapterno@. }%
  \def\\{\unskip\space\ignorespaces}\ifrunheads@\headlinefont@
    \ifodd\pageno\rightheadline \else\leftheadline\fi
   \else\hfil\fi\ifNoRunHeadline\global\NoRunHeadlinefalse\fi}
\let\folio@\folio
\def\foliofont@{\foliofont}
\def\foliofont{\eightrm}
\def\headlinefont@{\headlinefont}
\def\headlinefont{\eightpoint\smc}
\def\leftheadline{\rlap{\folio}\hfill
   \ifNoRunHeadline\else\iftrue\topmark\fi\fi \hfill}
\def\rightheadline{\hfill\ifNoRunHeadline
   \else \expandafter\iffalse\botmark\fi\fi
  \hfill \llap{\folio}}
\footline={{\eightpoint\bottremark}%
   \ifrunheads@\else\hfil{\let\foliofont\tenrm\folio}\fi\hfil}
\def\bottremark{}
 
%Definition von \norunninghead:
\newif\ifNoRunHeadline      
\def\norunninghead{\global\NoRunHeadlinetrue}
\norunninghead

\output={\output@}
%\def\output@{\shipout\vbox{%
% \iffirstpage@ \global\firstpage@false
%  \pagebody \logo@ \makefootline%
% \else \ifrunheads@ \makeheadline \pagebody
%       \else \pagebody \makefootline \fi
% \fi}%
% \advancepageno \ifnum\outputpenalty>-\@MM\else\dosupereject\fi}
%
%Modifizierter Output + Index-Output:
\newif\ifoffset\offsetfalse
\output={\output@}
\def\output@{%
 \ifoffset 
  \ifodd\count0\advance\hoffset by0.5truecm
   \else\advance\hoffset by-0.5truecm\fi\fi
 \shipout\vbox{%
  \makeheadline \pagebody \makefootline }%
 \advancepageno \ifnum\outputpenalty>-\@MM\else\dosupereject\fi}

\def\indexoutput#1{%
  \ifoffset 
   \ifodd\count0\advance\hoffset by0.5truecm
    \else\advance\hoffset by-0.5truecm\fi\fi
  \shipout\vbox{\makeheadline
  \vbox to\fullvsize{\boxmaxdepth\maxdepth%
  \ifvoid\topins\else\unvbox\topins\fi% 
  #1 %
  \ifvoid\footins\else % footnote info is present
    \vskip\skip\footins
    \footnoterule
    \unvbox\footins\fi
  \ifr@ggedbottom \kern-\dimen@ \vfil \fi}%
  \baselineskip2pc
  \makefootline}%
 \global\advance\pageno\@ne
 \ifnum\outputpenalty>-\@MM\else\dosupereject\fi}
 
 \newbox\partialpage \newdimen\halfsize \halfsize=0.5\fullhsize
 \advance\halfsize by-0.5em

 \def\begindoublecolumns{\output={\indexoutput{\unvbox255}}%
   \begingroup \def\line{\fullline}
   \output={\global\setbox\partialpage=\vbox{\unvbox255\bigskip}}\eject
   \output={\doublecolumnout}\hsize=\halfsize \vsize=2\fullvsize}
 \def\enddoublecolumns{\output={\balancecolumns}\eject
  \endgroup \pagegoal=\fullvsize%
  \output={\output@}}
\def\doublecolumnout{\splittopskip=\topskip \splitmaxdepth=\maxdepth
  \dimen@=\fullvsize \advance\dimen@ by-\ht\partialpage
  \setbox0=\vsplit255 to \dimen@ \setbox2=\vsplit255 to \dimen@
  \indexoutput{\pagesofar} \unvbox255 \penalty\outputpenalty}
\def\pagesofar{\unvbox\partialpage
  \wd0=\hsize \wd2=\hsize \hbox to\fullhsize{\box0\hfil\box2}}
\def\balancecolumns{\setbox0=\vbox{\unvbox255} \dimen@=\ht0
  \advance\dimen@ by\topskip \advance\dimen@ by-\baselineskip
  \divide\dimen@ by2 \splittopskip=\topskip
  {\vbadness=10000 \loop \global\setbox3=\copy0
    \global\setbox1=\vsplit3 to\dimen@
    \ifdim\ht3>\dimen@ \global\advance\dimen@ by1pt \repeat}
  \setbox0=\vbox to\dimen@{\unvbox1} \setbox2=\vbox to\dimen@{\unvbox3}
  \pagesofar}

\tenpoint
\catcode`\@=\active

\def\smallheadings{\let\chapheadfonts\tenpoint\let\headfonts\tenpoint}

\tenpoint
\catcode`\@=\active

\def\LL{\leavevmode\setbox0=\hbox{L}\hbox to\wd0{\hss\char'40L}}
\def\al{\alpha}

\def\om{\omega}

            %used for crossreferencing, Tex should ignore.
             %used for refencing (section-numbers)
          %used for new-section numbers
\def\P{{\Bbb P}}

\def\today{\ifcase\month\or
 January\or February\or March\or April\or May\or June\or
 July\or August\or September\or October\or November\or December\fi
 \space\number\day, \number\year}
 %zum Nummerieren

\def\({\left(}
\def\){\right)}
\def\[{\left[}
\def\]{\right]}

\def\3{\ss}
\catcode`\@=11
\def\dddot#1{\vbox{\ialign{##\crcr
      .\hskip-.5pt.\hskip-.5pt.\crcr\noalign{\kern1.5\p@\nointerlineskip}
      $\hfil\displaystyle{#1}\hfil$\crcr}}}

\newif\iftab@\tab@false
\newif\ifvtab@\vtab@false
\def\tab{\bgroup\tab@true\vtab@false\vst@bfalse\Strich@false%
   \def\\{\global\hline@@false%
     \ifhline@\global\hline@false\global\hline@@true\fi\cr}
   \edef\l@{\the\leftskip}\ialign\bgroup\hskip\l@##\hfil&&##\hfil\cr}
\def\endtab{\cr\egroup\egroup}
\def\vtab{\vtop\bgroup\vst@bfalse\vtab@true\tab@true\Strich@false%
   \bgroup\def\\{\cr}\ialign\bgroup&##\hfil\cr}
\def\endvtab{\cr\egroup\egroup\egroup}
\def\stab{\D@cke0.5pt\null 
 \bgroup\tab@true\vtab@false\vst@bfalse\Strich@true\Let@@\vspace@
 \normalbaselines\offinterlineskip
  \openup\spreadmlines@
 \edef\l@{\the\leftskip}\ialign
 \bgroup\hskip\l@##\hfil&&##\hfil\crcr}
\def\endstab{\crcr\egroup
 \egroup}
\newif\ifvst@b\vst@bfalse
\def\vstab{\D@cke0.5pt\null
 \vtop\bgroup\tab@true\vtab@false\vst@btrue\Strich@true\bgroup\Let@@\vspace@
 \normalbaselines\offinterlineskip
  \openup\spreadmlines@\bgroup}
\def\endvstab{\crcr\egroup\egroup
 \egroup\tab@false\Strich@false}

\newdimen\htstrut@
\htstrut@8.5\p@
\newdimen\htStrut@
\htStrut@12\p@
\newdimen\dpstrut@
\dpstrut@3.5\p@
\newdimen\dpStrut@
\dpStrut@3.5\p@
\def\openup{\afterassignment\@penup\dimen@=}
\def\@penup{\advance\lineskip\dimen@
  \advance\baselineskip\dimen@
  \advance\lineskiplimit\dimen@
  \divide\dimen@ by2
  \advance\htstrut@\dimen@
  \advance\htStrut@\dimen@
  \advance\dpstrut@\dimen@
  \advance\dpStrut@\dimen@}
\def\Let@@{\relax\iffalse{\fi%
    \def\\{\global\hline@@false%
     \ifhline@\global\hline@false\global\hline@@true\fi\cr}%
    \iffalse}\fi}
\def\matrix{\null\,\vcenter\bgroup
 \tab@false\vtab@false\vst@bfalse\Strich@false\Let@@\vspace@
 \normalbaselines\openup\spreadmlines@\ialign
 \bgroup\hfil$\m@th##$\hfil&&\quad\hfil$\m@th##$\hfil\crcr
 \Mathstrut@\crcr\noalign{\kern-\baselineskip}}
\def\endmatrix{\crcr\Mathstrut@\crcr\noalign{\kern-\baselineskip}\egroup
 \egroup\,}
\def\smatrix{\D@cke0.5pt\null\,
 \vcenter\bgroup\tab@false\vtab@false\vst@bfalse\Strich@true\Let@@\vspace@
 \normalbaselines\offinterlineskip
  \openup\spreadmlines@\ialign
 \bgroup\hfil$\m@th##$\hfil&&\quad\hfil$\m@th##$\hfil\crcr}
\def\endsmatrix{\crcr\egroup
 \egroup\,\Strich@false}
\newdimen\D@cke
\def\Dicke#1{\global\D@cke#1}
\newtoks\tabs@\tabs@{&}
\newif\ifStrich@\Strich@false
\newif\iff@rst

\def\Stricherr@{\iftab@\ifvtab@\errmessage{\noexpand\s not allowed
     here. Use \noexpand\vstab!}%
  \else\errmessage{\noexpand\s not allowed here. Use \noexpand\stab!}%
  \fi\else\errmessage{\noexpand\s not allowed
     here. Use \noexpand\smatrix!}\fi}
\def\format{\ifvst@b\else\crcr\fi\egroup\iffalse{\fi\ifnum`}=0 \fi\format@}
\def\format@#1\\{\def\preamble@{#1}%
 \def\Str@chfehlt##1{\ifx##1\s\Stricherr@\fi\ifx##1\\\let\Next\relax%
   \else\let\Next\Str@chfehlt\fi\Next}%
 \def\c{\hfil\noexpand\ifhline@@\hbox{\vrule height\htStrut@%
   depth\dpstrut@ width\z@}\noexpand\fi%
   \ifStrich@\hbox{\vrule height\htstrut@ depth\dpstrut@ width\z@}%
   \fi\iftab@\else$\m@th\fi\the\hashtoks@\iftab@\else$\fi\hfil}%
 \def\r{\hfil\noexpand\ifhline@@\hbox{\vrule height\htStrut@%
   depth\dpstrut@ width\z@}\noexpand\fi%
   \ifStrich@\hbox{\vrule height\htstrut@ depth\dpstrut@ width\z@}%
   \fi\iftab@\else$\m@th\fi\the\hashtoks@\iftab@\else$\fi}%
 \def\l{\noexpand\ifhline@@\hbox{\vrule height\htStrut@%
   depth\dpstrut@ width\z@}\noexpand\fi%
   \ifStrich@\hbox{\vrule height\htstrut@ depth\dpstrut@ width\z@}%
   \fi\iftab@\else$\m@th\fi\the\hashtoks@\iftab@\else$\fi\hfil}%
 \def\s{\ifStrich@\ \the\tabs@\vrule width\D@cke\the\hashtoks@%
          \fi\the\tabs@\ }%
 \def\sa{\ifStrich@\vrule width\D@cke\the\hashtoks@%
            \the\tabs@\ %
            \fi}%
 \def\se{\ifStrich@\ \the\tabs@\vrule width\D@cke\the\hashtoks@\fi}%
 \def\cd{\hfil\noexpand\ifhline@@\hbox{\vrule height\htStrut@%
   depth\dpstrut@ width\z@}\noexpand\fi%
   \ifStrich@\hbox{\vrule height\htstrut@ depth\dpstrut@ width\z@}%
   \fi$\dsize\m@th\the\hashtoks@$\hfil}%
 \def\rd{\hfil\noexpand\ifhline@@\hbox{\vrule height\htStrut@%
   depth\dpstrut@ width\z@}\noexpand\fi%
   \ifStrich@\hbox{\vrule height\htstrut@ depth\dpstrut@ width\z@}%
   \fi$\dsize\m@th\the\hashtoks@$}%
 \def\ld{\noexpand\ifhline@@\hbox{\vrule height\htStrut@%
   depth\dpstrut@ width\z@}\noexpand\fi%
   \ifStrich@\hbox{\vrule height\htstrut@ depth\dpstrut@ width\z@}%
   \fi$\dsize\m@th\the\hashtoks@$\hfil}%
 \ifStrich@\else\Str@chfehlt#1\\\fi%
 \setbox\z@\hbox{\xdef\Preamble@{\preamble@}}\ifnum`{=0 \fi\iffalse}\fi
 \ialign\bgroup\span\Preamble@\crcr}
\newif\ifhline@\hline@false
\newif\ifhline@@\hline@@false
\def\hlinefor#1{\multispan@{\strip@#1 }\leaders\hrule height\D@cke\hfill%
    \global\hline@true\ignorespaces}
\def\Item "#1"{\par\noindent\hangindent2\parindent%
  \hangafter1\setbox0\hbox{\rm#1\enspace}\ifdim\wd0>2\parindent%
  \box0\else\hbox to 2\parindent{\rm#1\hfil}\fi\ignorespaces}
\def\ITEM #1"#2"{\par\noindent\hangafter1\hangindent#1%
  \setbox0\hbox{\rm#2\enspace}\ifdim\wd0>#1%
  \box0\else\hbox to 0pt{\rm#2\hss}\hskip#1\fi\ignorespaces}
\def\item"#1"{\par\noindent\hang%
  \setbox0=\hbox{\rm#1\enspace}\ifdim\wd0>\the\parindent%
  \box0\else\hbox to \parindent{\rm#1\hfil}\enspace\fi\ignorespaces}
\let\plainitem@\item
\catcode`\@=13

\catcode`\@=11
\font\tenln    = line10
\font\tenlnw   = linew10

\newskip\Einheit \Einheit=0.5cm
\newcount\xcoord \newcount\ycoord
\newdimen\xdim \newdimen\ydim \newdimen\PfadD@cke \newdimen\Pfadd@cke

%%%%%%%%%%%%%%%%%%%%%%%%%%%%%%%%%%%%%%%%%%%%%%%%%
%LaTeX counters, dimensions, variables for lines%
%%%%%%%%%%%%%%%%%%%%%%%%%%%%%%%%%%%%%%%%%%%%%%%%%
\newcount\@tempcnta
\newcount\@tempcntb

\newdimen\@tempdima
\newdimen\@tempdimb

\newdimen\@wholewidth
\newdimen\@halfwidth

\newcount\@xarg
\newcount\@yarg
\newcount\@yyarg
\newbox\@linechar
\newbox\@tempboxa
\newdimen\@linelen
\newdimen\@clnwd
\newdimen\@clnht

\newif\if@negarg

\def\@whilenoop#1{}
\def\@whiledim#1\do #2{\ifdim #1\relax#2\@iwhiledim{#1\relax#2}\fi}
\def\@iwhiledim#1{\ifdim #1\let\@nextwhile=\@iwhiledim
        \else\let\@nextwhile=\@whilenoop\fi\@nextwhile{#1}}

\def\@whileswnoop#1\fi{}
\def\@whilesw#1\fi#2{#1#2\@iwhilesw{#1#2}\fi\fi}
\def\@iwhilesw#1\fi{#1\let\@nextwhile=\@iwhilesw
         \else\let\@nextwhile=\@whileswnoop\fi\@nextwhile{#1}\fi}

\def\thinlines{\let\@linefnt\tenln \let\@circlefnt\tencirc
  \@wholewidth\fontdimen8\tenln \@halfwidth .5\@wholewidth}
\def\thicklines{\let\@linefnt\tenlnw \let\@circlefnt\tencircw
  \@wholewidth\fontdimen8\tenlnw \@halfwidth .5\@wholewidth}
\thinlines
%%%%%%%%%%%%%%%%%%%%%%%%%%%%%%%%%%%%%%%%%%%%%%%%%%%%%%%%%%%

\PfadD@cke1pt \Pfadd@cke0.5pt
\def\PfadDicke#1{\PfadD@cke#1 \divide\PfadD@cke by2 \Pfadd@cke\PfadD@cke \multiply\PfadD@cke by2}
\long\def\LOOP#1\REPEAT{\def\BODY{#1}\ITERATE}
\def\ITERATE{\BODY \let\next\ITERATE \else\let\next\relax\fi \next}
\let\REPEAT=\fi
\def\Punkt{\hbox{\raise-2pt\hbox to0pt{\hss$\ssize\bullet$\hss}}}
\def\DuennPunkt(#1,#2){\unskip
  \raise#2 \Einheit\hbox to0pt{\hskip#1 \Einheit
          \raise-2.5pt\hbox to0pt{\hss$\bullet$\hss}\hss}}
\def\NormalPunkt(#1,#2){\unskip
  \raise#2 \Einheit\hbox to0pt{\hskip#1 \Einheit
          \raise-3pt\hbox to0pt{\hss\twelvepoint$\bullet$\hss}\hss}}
\def\DickPunkt(#1,#2){\unskip
  \raise#2 \Einheit\hbox to0pt{\hskip#1 \Einheit
          \raise-4pt\hbox to0pt{\hss\fourteenpoint$\bullet$\hss}\hss}}
\def\Kreis(#1,#2){\unskip
  \raise#2 \Einheit\hbox to0pt{\hskip#1 \Einheit
          \raise-4pt\hbox to0pt{\hss\fourteenpoint$\circ$\hss}\hss}}

%%%%%%%%%%%%%%%%%%%%%
%LaTeX line macros%
%%%%%%%%%%%%%%%%%%%%%
\def\Line@(#1,#2)#3{\@xarg #1\relax \@yarg #2\relax
\@linelen=#3\Einheit
\ifnum\@xarg =0 \@vline
  \else \ifnum\@yarg =0 \@hline \else \@sline\fi
\fi}

\def\@sline{\ifnum\@xarg< 0 \@negargtrue \@xarg -\@xarg \@yyarg -\@yarg
  \else \@negargfalse \@yyarg \@yarg \fi
\ifnum \@yyarg >0 \@tempcnta\@yyarg \else \@tempcnta -\@yyarg \fi
\ifnum\@tempcnta>6 \@badlinearg\@tempcnta0 \fi
\ifnum\@xarg>6 \@badlinearg\@xarg 1 \fi
\setbox\@linechar\hbox{\@linefnt\@getlinechar(\@xarg,\@yyarg)}%
\ifnum \@yarg >0 \let\@upordown\raise \@clnht\z@
   \else\let\@upordown\lower \@clnht \ht\@linechar\fi
\@clnwd=\wd\@linechar
\if@negarg \hskip -\wd\@linechar \def\@tempa{\hskip -2\wd\@linechar}\else
     \let\@tempa\relax \fi
\@whiledim \@clnwd <\@linelen \do
  {\@upordown\@clnht\copy\@linechar
   \@tempa
   \advance\@clnht \ht\@linechar
   \advance\@clnwd \wd\@linechar}%
\advance\@clnht -\ht\@linechar
\advance\@clnwd -\wd\@linechar
\@tempdima\@linelen\advance\@tempdima -\@clnwd
\@tempdimb\@tempdima\advance\@tempdimb -\wd\@linechar
\if@negarg \hskip -\@tempdimb \else \hskip \@tempdimb \fi
\multiply\@tempdima \@m
\@tempcnta \@tempdima \@tempdima \wd\@linechar \divide\@tempcnta \@tempdima
\@tempdima \ht\@linechar \multiply\@tempdima \@tempcnta
\divide\@tempdima \@m
\advance\@clnht \@tempdima
\ifdim \@linelen <\wd\@linechar
   \hskip \wd\@linechar
  \else\@upordown\@clnht\copy\@linechar\fi}

\def\@hline{\ifnum \@xarg <0 \hskip -\@linelen \fi
\vrule height\Pfadd@cke width \@linelen depth\Pfadd@cke
\ifnum \@xarg <0 \hskip -\@linelen \fi}

\def\@getlinechar(#1,#2){\@tempcnta#1\relax\multiply\@tempcnta 8
\advance\@tempcnta -9 \ifnum #2>0 \advance\@tempcnta #2\relax\else
\advance\@tempcnta -#2\relax\advance\@tempcnta 64 \fi
\char\@tempcnta}

\def\Vektor(#1,#2)#3(#4,#5){\unskip\leavevmode
  \xcoord#4\relax \ycoord#5\relax
      \raise\ycoord \Einheit\hbox to0pt{\hskip\xcoord \Einheit
         \Vector@(#1,#2){#3}\hss}}

\def\Vector@(#1,#2)#3{\@xarg #1\relax \@yarg #2\relax
\@tempcnta \ifnum\@xarg<0 -\@xarg\else\@xarg\fi
\ifnum\@tempcnta<5\relax
\@linelen=#3\Einheit
\ifnum\@xarg =0 \@vvector
  \else \ifnum\@yarg =0 \@hvector \else \@svector\fi
\fi
\else\@badlinearg\fi}

\def\@hvector{\@hline\hbox to 0pt{\@linefnt
\ifnum \@xarg <0 \@getlarrow(1,0)\hss\else
    \hss\@getrarrow(1,0)\fi}}

\def\@vvector{\ifnum \@yarg <0 \@downvector \else \@upvector \fi}

\def\@svector{\@sline
\@tempcnta\@yarg \ifnum\@tempcnta <0 \@tempcnta=-\@tempcnta\fi
\ifnum\@tempcnta <5
  \hskip -\wd\@linechar
  \@upordown\@clnht \hbox{\@linefnt  \if@negarg
  \@getlarrow(\@xarg,\@yyarg) \else \@getrarrow(\@xarg,\@yyarg) \fi}%
\else\@badlinearg\fi}

\def\@upline{\hbox to \z@{\hskip -.5\Pfadd@cke \vrule width \Pfadd@cke
   height \@linelen depth \z@\hss}}

\def\@downline{\hbox to \z@{\hskip -.5\Pfadd@cke \vrule width \Pfadd@cke
   height \z@ depth \@linelen \hss}}

\def\@upvector{\@upline\setbox\@tempboxa\hbox{\@linefnt\char'66}\raise
     \@linelen \hbox to\z@{\lower \ht\@tempboxa\box\@tempboxa\hss}}

\def\@downvector{\@downline\lower \@linelen
      \hbox to \z@{\@linefnt\char'77\hss}}

\def\@getlarrow(#1,#2){\ifnum #2 =\z@ \@tempcnta='33\else
\@tempcnta=#1\relax\multiply\@tempcnta \sixt@@n \advance\@tempcnta
-9 \@tempcntb=#2\relax\multiply\@tempcntb \tw@
\ifnum \@tempcntb >0 \advance\@tempcnta \@tempcntb\relax
\else\advance\@tempcnta -\@tempcntb\advance\@tempcnta 64
\fi\fi\char\@tempcnta}

\def\@getrarrow(#1,#2){\@tempcntb=#2\relax
\ifnum\@tempcntb < 0 \@tempcntb=-\@tempcntb\relax\fi
\ifcase \@tempcntb\relax \@tempcnta='55 \or
\ifnum #1<3 \@tempcnta=#1\relax\multiply\@tempcnta
24 \advance\@tempcnta -6 \else \ifnum #1=3 \@tempcnta=49
\else\@tempcnta=58 \fi\fi\or
\ifnum #1<3 \@tempcnta=#1\relax\multiply\@tempcnta
24 \advance\@tempcnta -3 \else \@tempcnta=51\fi\or
\@tempcnta=#1\relax\multiply\@tempcnta
\sixt@@n \advance\@tempcnta -\tw@ \else
\@tempcnta=#1\relax\multiply\@tempcnta
\sixt@@n \advance\@tempcnta 7 \fi\ifnum #2<0 \advance\@tempcnta 64 \fi
\char\@tempcnta}
%%%%%%%%%%%%%%%%%%%%%%%%%%%%%%%%%%%%%%%%%%%%%%%%%%%%%%%%%%%%%

\def\Diagonale(#1,#2)#3{\unskip\leavevmode
  \xcoord#1\relax \ycoord#2\relax
      \raise\ycoord \Einheit\hbox to0pt{\hskip\xcoord \Einheit
         \Line@(1,1){#3}\hss}}
\def\AntiDiagonale(#1,#2)#3{\unskip\leavevmode
  \xcoord#1\relax \ycoord#2\relax %\advance\xcoord by -0.05\relax
      \raise\ycoord \Einheit\hbox to0pt{\hskip\xcoord \Einheit
         \Line@(1,-1){#3}\hss}}
\def\Pfad(#1,#2),#3\endPfad{\unskip\leavevmode
  \xcoord#1 \ycoord#2 \thicklines\ZeichnePfad#3\endPfad\thinlines}
\def\ZeichnePfad#1{\ifx#1\endPfad\let\next\relax
  \else\let\next\ZeichnePfad
    \ifnum#1=1
      \raise\ycoord \Einheit\hbox to0pt{\hskip\xcoord \Einheit
         \vrule height\Pfadd@cke width1 \Einheit depth\Pfadd@cke\hss}%
      \advance\xcoord by 1
    \else\ifnum#1=2
      \raise\ycoord \Einheit\hbox to0pt{\hskip\xcoord \Einheit
        \hbox{\hskip-\PfadD@cke\vrule height1 \Einheit width\PfadD@cke depth0pt}\hss}%
      \advance\ycoord by 1
    \else\ifnum#1=3
      \raise\ycoord \Einheit\hbox to0pt{\hskip\xcoord \Einheit
         \Line@(1,1){1}\hss}
      \advance\xcoord by 1
      \advance\ycoord by 1
    \else\ifnum#1=4
      \raise\ycoord \Einheit\hbox to0pt{\hskip\xcoord \Einheit
         \Line@(1,-1){1}\hss}
      \advance\xcoord by 1
      \advance\ycoord by -1
    \else\ifnum#1=5
      \advance\xcoord by -1
      \raise\ycoord \Einheit\hbox to0pt{\hskip\xcoord \Einheit
         \vrule height\Pfadd@cke width1 \Einheit depth\Pfadd@cke\hss}%
    \else\ifnum#1=6
      \advance\ycoord by -1
      \raise\ycoord \Einheit\hbox to0pt{\hskip\xcoord \Einheit
        \hbox{\hskip-\PfadD@cke\vrule height1 \Einheit width\PfadD@cke depth0pt}\hss}%
    \else\ifnum#1=7
      \advance\xcoord by -1
      \advance\ycoord by -1
      \raise\ycoord \Einheit\hbox to0pt{\hskip\xcoord \Einheit
         \Line@(1,1){1}\hss}
    \else\ifnum#1=8
      \advance\xcoord by -1
      \advance\ycoord by +1
      \raise\ycoord \Einheit\hbox to0pt{\hskip\xcoord \Einheit
         \Line@(1,-1){1}\hss}
    \fi\fi\fi\fi
    \fi\fi\fi\fi
  \fi\next}
\def\hSSchritt{\leavevmode\raise-.4pt\hbox to0pt{\hss.\hss}\hskip.2\Einheit
  \raise-.4pt\hbox to0pt{\hss.\hss}\hskip.2\Einheit
  \raise-.4pt\hbox to0pt{\hss.\hss}\hskip.2\Einheit
  \raise-.4pt\hbox to0pt{\hss.\hss}\hskip.2\Einheit
  \raise-.4pt\hbox to0pt{\hss.\hss}\hskip.2\Einheit}
\def\vSSchritt{\vbox{\baselineskip.2\Einheit\lineskiplimit0pt
\hbox{.}\hbox{.}\hbox{.}\hbox{.}\hbox{.}}}
\def\DSSchritt{\leavevmode\raise-.4pt\hbox to0pt{%
  \hbox to0pt{\hss.\hss}\hskip.2\Einheit
  \raise.2\Einheit\hbox to0pt{\hss.\hss}\hskip.2\Einheit
  \raise.4\Einheit\hbox to0pt{\hss.\hss}\hskip.2\Einheit
  \raise.6\Einheit\hbox to0pt{\hss.\hss}\hskip.2\Einheit
  \raise.8\Einheit\hbox to0pt{\hss.\hss}\hss}}
\def\dSSchritt{\leavevmode\raise-.4pt\hbox to0pt{%
  \hbox to0pt{\hss.\hss}\hskip.2\Einheit
  \raise-.2\Einheit\hbox to0pt{\hss.\hss}\hskip.2\Einheit
  \raise-.4\Einheit\hbox to0pt{\hss.\hss}\hskip.2\Einheit
  \raise-.6\Einheit\hbox to0pt{\hss.\hss}\hskip.2\Einheit
  \raise-.8\Einheit\hbox to0pt{\hss.\hss}\hss}}
\def\SPfad(#1,#2),#3\endSPfad{\unskip\leavevmode
  \xcoord#1 \ycoord#2 \ZeichneSPfad#3\endSPfad}
\def\ZeichneSPfad#1{\ifx#1\endSPfad\let\next\relax
  \else\let\next\ZeichneSPfad
    \ifnum#1=1
      \raise\ycoord \Einheit\hbox to0pt{\hskip\xcoord \Einheit
         \hSSchritt\hss}%
      \advance\xcoord by 1
    \else\ifnum#1=2
      \raise\ycoord \Einheit\hbox to0pt{\hskip\xcoord \Einheit
        \hbox{\hskip-2pt \vSSchritt}\hss}%
      \advance\ycoord by 1
    \else\ifnum#1=3
      \raise\ycoord \Einheit\hbox to0pt{\hskip\xcoord \Einheit
         \DSSchritt\hss}
      \advance\xcoord by 1
      \advance\ycoord by 1
    \else\ifnum#1=4
      \raise\ycoord \Einheit\hbox to0pt{\hskip\xcoord \Einheit
         \dSSchritt\hss}
      \advance\xcoord by 1
      \advance\ycoord by -1
    \else\ifnum#1=5
      \advance\xcoord by -1
      \raise\ycoord \Einheit\hbox to0pt{\hskip\xcoord \Einheit
         \hSSchritt\hss}%
    \else\ifnum#1=6
      \advance\ycoord by -1
      \raise\ycoord \Einheit\hbox to0pt{\hskip\xcoord \Einheit
        \hbox{\hskip-2pt \vSSchritt}\hss}%
    \else\ifnum#1=7
      \advance\xcoord by -1
      \advance\ycoord by -1
      \raise\ycoord \Einheit\hbox to0pt{\hskip\xcoord \Einheit
         \DSSchritt\hss}
    \else\ifnum#1=8
      \advance\xcoord by -1
      \advance\ycoord by 1
      \raise\ycoord \Einheit\hbox to0pt{\hskip\xcoord \Einheit
         \dSSchritt\hss}
    \fi\fi\fi\fi
    \fi\fi\fi\fi
  \fi\next}
\def\Koordinatenachsen(#1,#2){\unskip
 \hbox to0pt{\hskip-.5pt\vrule height#2 \Einheit width.5pt depth1 \Einheit}%
 \hbox to0pt{\hskip-1 \Einheit \xcoord#1 \advance\xcoord by1
    \vrule height0.25pt width\xcoord \Einheit depth0.25pt\hss}}
\def\Koordinatenachsen(#1,#2)(#3,#4){\unskip
 \hbox to0pt{\hskip-.5pt \ycoord-#4 \advance\ycoord by1
    \vrule height#2 \Einheit width.5pt depth\ycoord \Einheit}%
 \hbox to0pt{\hskip-1 \Einheit \hskip#3\Einheit 
    \xcoord#1 \advance\xcoord by1 \advance\xcoord by-#3 
    \vrule height0.25pt width\xcoord \Einheit depth0.25pt\hss}}
\def\Gitter(#1,#2){\unskip \xcoord0 \ycoord0 \leavevmode
  \LOOP\ifnum\ycoord<#2
    \loop\ifnum\xcoord<#1
      \raise\ycoord \Einheit\hbox to0pt{\hskip\xcoord \Einheit\Punkt\hss}%
      \advance\xcoord by1
    \repeat
    \xcoord0
    \advance\ycoord by1
  \REPEAT}
\def\Gitter(#1,#2)(#3,#4){\unskip \xcoord#3 \ycoord#4 \leavevmode
  \LOOP\ifnum\ycoord<#2
    \loop\ifnum\xcoord<#1
      \raise\ycoord \Einheit\hbox to0pt{\hskip\xcoord \Einheit\Punkt\hss}%
      \advance\xcoord by1
    \repeat
    \xcoord#3
    \advance\ycoord by1
  \REPEAT}
\def\Label#1#2(#3,#4){\unskip \xdim#3 \Einheit \ydim#4 \Einheit
  \def\lo{\advance\xdim by-.5 \Einheit \advance\ydim by.5 \Einheit}%
  \def\llo{\advance\xdim by-.25cm \advance\ydim by.5 \Einheit}%
  \def\loo{\advance\xdim by-.5 \Einheit \advance\ydim by.25cm}%
  \def\o{\advance\ydim by.25cm}%
  \def\ro{\advance\xdim by.5 \Einheit \advance\ydim by.5 \Einheit}%
  \def\rro{\advance\xdim by.25cm \advance\ydim by.5 \Einheit}%
  \def\roo{\advance\xdim by.5 \Einheit \advance\ydim by.25cm}%
  \def\l{\advance\xdim by-.30cm}%
  \def\r{\advance\xdim by.30cm}%
  \def\lu{\advance\xdim by-.5 \Einheit \advance\ydim by-.6 \Einheit}%
  \def\llu{\advance\xdim by-.25cm \advance\ydim by-.6 \Einheit}%
  \def\luu{\advance\xdim by-.5 \Einheit \advance\ydim by-.30cm}%
  \def\u{\advance\ydim by-.30cm}%
  \def\ru{\advance\xdim by.5 \Einheit \advance\ydim by-.6 \Einheit}%
  \def\rru{\advance\xdim by.25cm \advance\ydim by-.6 \Einheit}%
  \def\ruu{\advance\xdim by.5 \Einheit \advance\ydim by-.30cm}%
  #1\raise\ydim\hbox to0pt{\hskip\xdim
     \vbox to0pt{\vss\hbox to0pt{\hss$#2$\hss}\vss}\hss}%
}
\catcode`\@=13

\hsize13cm
\vsize19cm
\newdimen\fullhsize
\newdimen\fullvsize
\newdimen\halfsize
\fullhsize13cm
\fullvsize19cm
\halfsize=0.5\fullhsize
\advance\halfsize by-0.5em

\magnification1200

\TagsOnRight

\def\AignAB{1}
\def\CJKrAA{2}
\def\ComtAA{3}
\def\KratBD{4}
\def\SlatAC{5}
\def\OEIS{6}
\def\StanBI{7}
\def\VienAE{8}

%equation numbers
\def\IA{1.1}
\def\IB{1.2}
\def\IC{1.3}
\def\ID{1.4}
\def\CKa{1.5}
\def\CKb{1.6}
\def\CKc{1.7}
\def\AA{1.8}
\def\AB{1.9}
\def\AAa{2.1}
\def\AG{2.2}
\def\AAb{2.3}
\def\AAc{2.4}
\def\SA{2.5}
\def\SB{2.6}
\def\SC{2.7}
\def\SD{2.8}
\def\SE{2.9}
\def\SF{2.10}
\def\SK{2.11}
\def\SL{2.12}
\def\SM{2.13}
\def\SN{2.14}
\def\SO{2.15}
\def\SR{2.16}

\def\AC{3.1}
\def\AK{3.2}
\def\AH{3.3}
\def\AHa{3.4}
\def\AD{3.5}
\def\AF{3.6}
\def\AI{3.7}
\def\AJ{3.8}
\def\AKa{3.9}
\def\AL{3.10}
\def\ALa{4.1}
\def\AM{4.2}
\def\AN{4.3}
\def\AO{4.4}
\def\AOa{4.5}
\def\AOb{4.6}
\def\AP{4.7}
\def\APa{4.8}
\def\AQ{4.9}
\def\AQa{4.10}
\def\AQb{4.11}
\def\AR{5.1}
\def\AS{5.2}
\def\ASa{5.3}
\def\AT{5.4}
\def\AU{5.5}
\def\AV{5.6}
\def\AW{5.7}
\def\AX{5.8}
\def\AY{5.9}
\def\BA{5.10}
\def\BJ{5.11}
\def\BK{5.12}
\def\BL{5.13}
\def\BM{5.14}
\def\BN{5.15}
\def\BO{5.16}
\def\CA{6.1}
\def\CAa{6.2}
\def\CAb{6.3}
\def\CB{6.4}
\def\CAc{6.5}
\def\CAd{6.6}
\def\CC{6.7}
\def\CAe{6.8}
\def\CAf{6.9}
\def\CD{6.10}
\def\CAg{6.11}
\def\CAh{6.12}
\def\CAi{6.13}
\def\CE{6.14}
\def\CF{6.15}
\def\CG{6.16}
\def\CH{6.17}
\def\CI{6.18}

%theorem numbers
\def\TCKA{1}
\def\TCKB{2}
\def\TA{3}
\def\TB{4}
\def\TF{5}
\def\TC{6}
\def\TEa{7}
\def\TE{8}
\def\TG{9}
\def\TH{10}
\def\TD{11}
\def\TI{12}
\def\TJ{13}
\def\TK{14}
\def\TL{15}
\def\TM{16}
\def\TN{17}
\def\TO{18}
\def\TP{19}

%figure numbers
\def\FA{1}

\def\P{{\Cal P}}

\def\MP{\operatorname{\text{\it MP\/}}}

\def\fl#1{\lfloor#1\rfloor}
\def\cl#1{\lceil#1\rceil}

\def\po#1#2{(#1)_#2} 

\topmatter 
\title Some determinants of path generating functions, II
\endtitle 
\author C.~Krattenthaler$^{\dagger*}$ and D. Yaqubi$^{\ddagger**}$
\endauthor 
\affil 
$^\dagger$Fakult\"at f\"ur Mathematik, Universit\"at Wien,\\
Oskar-Norgenstern-Platz~1, A-1090 Vienna, Austria.\\
WWW: {\tt http://www.mat.univie.ac.at/\~{}kratt}\\\vskip6pt
$^\ddagger$Department of Mathematics, Ferdowsi University of Mashhad,\\
 Mashhad, Iran
\endaffil
\address Fakult\"at f\"ur Mathematik, Universit\"at Wien,
Oskar-Norgenstern-Platz~1, A-1090 Vienna, Austria.
WWW: {\tt http://www.mat.univie.ac.at/\~{}kratt}
\endaddress
\address
Department of Mathematics, Ferdowsi University Of Mashhad, Mashhad, Iran
\endaddress
\thanks $^*$Research partially supported 
by the Austrian Science Foundation FWF, grant S50-N15,
in the framework of the Special Research Program
``Algorithmic and Enumerative Combinatorics".%
\newline
\indent $^{**}$Current address: Faculty of Agriculture and Animal Science, 
University of Torbat-e Jam, Iran%
\endthanks

\subjclass Primary 05A19;
 Secondary 05A10 05A15 11C20 15A15
\endsubjclass
\keywords Hankel determinants, Catalan numbers, ballot numbers, 
Motzkin numbers, Motzkin prefix numbers, Motzkin paths 
\endkeywords
\abstract 
We evaluate Hankel determinants of matrices in which the entries are
generating functions for paths
consisting of up-steps, down-steps and level steps with a fixed
starting point but variable end point. By specialisation,
these determinant evaluations have numerous corollaries. 
In particular, one consequence is that the Hankel determinant
of Motzkin prefix numbers equals~1, regardless of the size of
the Hankel matrix.
\endabstract
\endtopmatter

\document

\subhead 1. Introduction\endsubhead
Determinants (and Hankel determinants in particular) of path counting numbers
(respectively, more generally, of path generating functions) are
ubiquitous in the literature. Their ``popularity" stems in part from
the fact that, frequently, such determinants can be 
evaluated into attractive, compact closed formulae. 
This article contributes further to this body of results.

\midinsert
%\vskip10pt
$$
\PfadDicke{.5pt}
\Gitter(12,5)(-1,0)
\Koordinatenachsen(12,5)(-1,0)
\Pfad(0,0),33141433144\endPfad
\DickPunkt(0,0)
\DickPunkt(11,0)
\PfadDicke{.5pt}
\hbox{\hskip7.5cm}
\Gitter(9,5)(-1,0)
\Koordinatenachsen(9,5)(-1,0)
\Pfad(0,0),33344344\endPfad
\DickPunkt(0,0)
\DickPunkt(8,0)
\hbox{\hskip3.5cm}
$$
\centerline{\eightpoint a. A Motzkin path\quad \kern4cm
b. A Catalan path}
\vskip7pt
\centerline{\eightpoint Figure \FA}
\endinsert

The determinants that we consider here involve matrices formed from
numbers and generating functions of {\it three-step paths}. More
precisely, our paths consist of up-steps $(1,1)$, level steps $(1,0)$,
and down-steps $(1,-1)$. The number of such paths from $(0,0)$ to
$(n,0)$ that never run below the $x$-axis is known as the {\it Motzkin
number} $M_n$ (cf\. \cite{\StanBI, Exercise~6.38}; Figure~\FA.a shows
an example of a path contributing to $M_{11}$). On the other hand,
the number of paths from $(0,0)$ to $(2n,0)$ that consist only of
up-steps $(1,1)$ and down-steps $(1,-1)$ and do not run below the
$x$-axis is known as the {\it Catalan number} $C_n=\frac {1}
{n+1}\binom {2n}n$ (cf\. \cite{\StanBI, Exercise~6.19};
Figure~\FA.b shows an example of a path contributing to $C_{8}$).
It is well-known that
$$\align 
\det_{0\le i,j\le n-1}(C_{i+j})&=1,
\tag\IA\\
\det_{0\le i,j\le n-1}(C_{i+j+1})&=1,
\tag\IB\\
\det_{0\le i,j\le n-1}(M_{i+j})&=1,
\tag\IC\\
\det_{0\le i,j\le n-1}(M_{i+j+1})&=\cases 
(-1)^{n/3}&\text{if }n\equiv 0\pmod 3,\\
(-1)^{(n-1)/3}&\text{if }n\equiv 1\pmod 3,\\
0&\text{if }n\equiv 2\pmod 3,\\
\endcases
\tag\ID
\endalign$$
see e.g\. \cite{\VienAE, \AignAB}).

In \cite{\CJKrAA}, these Hankel determinant evaluations were
generalised to Hankel determinant evaluations of path generating
functions as follows (among others).
We define $\P_n(k ,l )$ as the generating function $\sum _{P} ^{}w(P),$
where $P$ runs over all three-step paths from $(0,k )$ to $(n,l )$, 
and where $w(P)$ is the
product of all weights of the steps of $P$, where the weights of the
steps are defined by $w((1,0))=x+y$, $w((1,1))=1$, and $w((1,-1))=xy$.
Furthermore, let $\P^+_n(k ,l )$ be the analogous generating function
$\sum _{P} 
^{}w(P)$, where $P$ runs over the subset of the set of the above
three-step paths which never run below the $x$-axis.
Clearly, if we specialise $x=-y=\sqrt{-1}$, then
$\P^+_{2n}(0,0)$ reduces to $C_n$ (and $\P^+_{2n+1}(0,0)=0$ for all~$n$), 
while, if we specialise $x=\frac {1} {2}(1+\sqrt{-3})$,
$y=\frac {1} {2}(1-\sqrt{-3})$, then
$\P^+_{n}(0,0)$ reduces to $M_n$. 
The somewhat unusual parametrisation that we have chosen
here turns out to be useful in the context of the Hankel determinant
evaluations of \cite{\CJKrAA} and of the present article, in that
the evaluations can be much more elegantly presented than it would be
possible under more straightforward parametrisations.

The theorems from \cite{\CJKrAA} that generalise (\IA)--(\ID) to the
weighted setting are the following two.

\proclaim{Theorem \TCKA\ {\smc (\cite{\CJKrAA, Theorem~1})}} 
For all positive integers $n$ and non-negative integers $k$, we have
$$\det_{0\le i,j\le n-1}\(\P^+_{i+j}(0,k)\)=\cases (-1)^{n_1\binom
{k+1}2}(xy)^{(k+1)^2\binom {n_1}2}&n=n_1(k+1),\\
0&n\not\equiv 0~(\text {\rm mod }k+1).
\endcases
\tag\CKa$$
\endproclaim

\proclaim{Theorem \TCKB\ {\smc (\cite{\CJKrAA, Theorem~2})}}
For all positive integers $n$ and non-negative integers $k$, we have
$$\multline 
\det_{0\le i,j\le n-1}\(\P^+_{i+j+1}(0,k)\)\\
=\cases (-1)^{n_1\binom
{k+1}2}(xy)^{(k+1)^2\binom {n_1}2}
\frac {y^{(k+1)(n_1+1)}-x^{(k+1)(n_1+1)}} {y^{k+1}-x^{k+1}}&n=n_1(k+1),\\
(-1)^{n_1\binom
{k+1}2+\binom k2}(xy)^{(k+1)^2\binom {n_1}2+n_1k(k+1)}\\
\kern2cm\times\frac {y^{(k+1)(n_1+1)}-x^{(k+1)(n_1+1)}}
{y^{k+1}-x^{k+1}}&n=n_1(k+1)+k,\\
0&n\not\equiv 0,k~(\text {\rm mod }k+1).
\endcases
\endmultline
\tag\CKb$$
\endproclaim

\remark{Remark}
If $k=0$, the formulae in Theorems~\TCKA\ and \TCKB\ have to be read
according to the convention that only the first line on the right-hand
sides of (\CKa) and (\CKb) applies; that is,
$$\det_{0\le i,j\le n-1}\(\P^+_{i+j}(0,0)\)=
(xy)^{\binom {n}2}
$$
and
$$
\det_{0\le i,j\le n-1}\(\P^+_{i+j+1}(0,0)\)
=(xy)^{\binom {n}2}
\frac {y^{n+1}-x^{n+1}} {y-x}.
$$
\endremark

The work on the present article began with computer experiments of the
second author on Hankel determinants of {\it Motzkin prefix numbers}. 
By definition, the $n$-th Motzkin prefix number is
the number of three-step paths consisting of $n$~steps, starting at the
origin, and not running below the $x$-axis (with any end point).
We denote this number by  $\MP_n$. The aforementioned 
computer experiments seemed to indicate that 
$$
\det_{0\le i,j\le n-1}(\MP_{i+j})=1
\tag\CKc$$
for all $n$. Subsequent consultation of the {\it On-Line Encyclopedia of
Integer Sequences} \cite{\OEIS, sequence~{\tt A005773}} 
revealed that this same observation
had already been made earlier by Philippe Del\'eham in 2007.
On the other hand, in view of the earlier work \cite{\CJKrAA}, the
first author was obviously led to look at the weighted generalisation
of the Hankel determinant in (\CKc), namely 
$$
\det_{0\le i,j\le n-1}\(\sum_{l \ge0}^{}\P^+_{i+j}(0 ,l )\),$$
or even more generally at
$$
\det_{0\le i,j\le n-1}\(\sum_{l \ge0}^{}\P^+_{i+j}(k ,l )\).
$$
The entries here are generating functions for three-step paths of a
given length that start at height~$k$ and never run below the
$x$-axis.

From here, it did not take very long to discover the closed form
evaluations of the Hankel determinants of path generating
functions in (\AA) and (\AB) below. Of course, these were at this
point only conjectures. 

\proclaim{Theorem \TA}
For all positive integers $n$ and non-negative integers $k$, we have
$$\det_{0\le i,j\le n-1}\(\sum_{l \ge0}\P^+_{i+j}(k ,l )\)=
\cases 
(-1)^{n_1\binom {k +1}2}(xy)^{(k +1)^2\binom {n_1+1}2-n},
&\text{if }n=(k +1)n_1,\\
(-1)^{n_1\binom {k +1}2}(xy)^{(k +1)^2\binom {n_1+1}2},
&\text{if }n=(k +1)n_1+1,\\
0,&\text{if }n\not\equiv0,1~(\text{\rm mod }k +1).
\endcases
\tag\AA$$
\endproclaim

\proclaim{Theorem \TB}
For all positive integers $n$ and non-negative integers $k$, we have
$$\multline
\det_{0\le i,j\le n-1}\(\sum_{l \ge0}\P^+_{i+j+1}(k,l )\)\\
=
%(xy)^{\binom {n}2}
%\frac {y^{n+1}+y^n-x^{n+1}-x^n} {y-x}.
\cases 
(-1)^{n_1\binom {k +1}2}(xy)^{(k +1)^2\binom {n_1+1}2-n}\\
\kern.3cm
\times
   \left(\frac {y^{(k+1)(n_1+1)}-x^{(k+1)(n_1+1)}} {y^{k+1}-x^{k+1}}
+(-1)^k\frac {y^{(k+1)n_1}-x^{(k+1)n_1}}
   {y^{k+1}-x^{k+1}}\right),
&\kern-5pt\text{if }n=(k +1)n_1,\\
(-1)^{n_1\binom {k +1}2}(xy)^{(k +1)^2\binom {n_1+1}2}\\
\kern.3cm
\times
\frac {(1+x)(1+y)\left(y^{(k+1)(n_1+1)}-x^{(k+1)(n_1+1)}\right)}
       {y^{k+1}-x^{k+1}}
&\kern-5pt\text{if }n=(k +1)n_1+1,\\
(-1)^{(n_1+1)\binom {k +1}2+k}(xy)^{(k +1)^2\binom {n_1+1}2+(k^2-1)(n_1+1)}\\
\kern.3cm
\times
\frac {(1+x)(1+y)\left(y^{(k+1)(n_1+1)}-x^{(k+1)(n_1+1)}\right)}
       {y^{k+1}-x^{k+1}}
&\kern-5pt\text{if }n=(k +1)n_1+k,\\
0,&\kern-5pt\text{if }n\not\equiv0,1,k~(\text{\rm mod }k +1).
\endcases\\
\endmultline
\tag\AB
$$
\endproclaim

\remark{Remarks}
(1) Also here, if $k=0$, the formulae in Theorems~\TA\ and \TB\ have to be read
according to the convention that only the first line on the right-hand
sides of (\AA) and (\AB) applies; that is,
$$\det_{0\le i,j\le n-1}\(\sum_{l \ge0}\P^+_{i+j}(0 ,l )\)=
(xy)^{\binom {n}2}
$$
and
$$
\det_{0\le i,j\le n-1}\(\sum_{l \ge0}\P^+_{i+j+1}(0,l )\)
=(xy)^{\binom {n}2}
   \left(\frac {y^{n+1}-x^{n+1}} {y-x}
+\frac {y^{n}-x^{n}}
   {y-x}\right).
$$

Similarly, if $k=1$ only the first two lines in Theorem~\TB\ apply;
that is,
$$\multline
\det_{0\le i,j\le n-1}\(\sum_{l \ge0}\P^+_{i+j+1}(1,l )\)\\
=
\cases 
(-1)^{n_1}(xy)^{4\binom {n_1+1}2-n}
   \left(\frac {y^{2(n_1+1)}-x^{2(n_1+1)}} {y^{2}-x^{2}}
-\frac {y^{2n_1}-x^{2n_1}}
   {y^{2}-x^{2}}\right),
&\kern-5pt\text{if }n=2n_1,\\
(-1)^{n_1}(xy)^{4\binom {n_1+1}2}
\frac {(1+x)(1+y)\left(y^{2(n_1+1)}-x^{2(n_1+1)}\right)}
       {y^2-x^2}
&\kern-5pt\text{if }n=2n_1+1.
\endcases\\
\endmultline
$$

\medskip
(2) Clearly, the specialisation $x=\frac {1} {2}(1+\sqrt{-3})$,
$y=\frac {1} {2}(1-\sqrt{-3})$, and $k=0$ of Theorem~\TA\ establishes
(\CKc). Many more interesting specialisations are possible, see Section~6.
\endremark

In the present article, we provide proofs for these
determinant evaluations. As it turns out, there is a ``connection
matrix" (see Section~3) which, upon multiplication on the left,
transforms the Hankel matrices on the left-hand sides of (\AA) and
(\AB) into the matrices on the left-hand sides of (\CKa) and (\CKb),
respectively, up to some ``correction" in the last row, see
Lemmas~\TD\ and \TI\ in Sections~4 and~5. This fact then allows us
to complete the evaluation of the determinants in Theorems~\TA\ and \TB\
by adapting arguments from the proofs of Theorems~\TCKA\ and \TCKB\
in \cite{\CJKrAA} to the new situation here, see the proofs of
Theorems~\TA\ and~\TB\ in Sections~4 and~5. Auxiliary results for
these proofs are collected in Section~3, which themselves depend on
elementary properties of our path generating functions that are
recalled in Section~2. We conclude our article with a list of
interesting specialisations of our two main theorems in Section~6.

\subhead 2. Elementary facts about three-step paths\endsubhead
In the proofs of our theorems, we need three elementary properties that
our path generating functions satisfy. We list them here as
(\AAa)--(\AAb). In the remainder of this section, we discuss four
types of specialisations of the path generating functions, which will
then be considered in Section~6 in the context of Theorems~\TA\ and \TB.

By retracing paths from the back to the beginning, one sees
that 
$$
\P^+_n(k ,l )=
(xy)^{k -l }\P^+_n(l ,k ),
\tag\AAa$$
and the same symmetry relation holds for $\P_n(k,l)$, but we shall not
have any need for the latter.

The reflection principle (see e.g\. \cite{\ComtAA, p.~22}) allows us to
express the generating functions $\P^+_n(k,l)$ for {\it
restricted\/} paths in terms of the generating functions $\P_n(k,l)$
for {\it unrestricted\/} paths, in terms of the relation
$$\P^+_n(k,l)=\P_n(k,l)-(xy)^{k+1}\P_n(-k-2,l).\tag\AG$$

In fact, elementary combinatorial reasoning shows that there is an
explicit formula for the path generating
function $\P_n(k,l)$, namely
$$
\P_n(k,l)=\sum _{s\ge0} ^{}\frac {n!} {s!\,(s+l-k)!\,(n-l+k-2s)!}
(x+y)^{n-l+k-2s}(xy)^{s}.
\tag\AAb
$$
In combination with (\AG), this also yields an explicit formula
for $\P^+_n(k,l)$. To be precise, we have
$$
\P^+_n(k ,l )=\sum_{s\ge0}\binom n{l -k +2s}
\left(\binom {l -k +2s}s - \binom {l -k +2s}{s-k -1}\right)
(x+y)^{n-l +k -2s}(x y)^s.
\tag\AAc
$$

\medskip
In Section~6, we shall discuss specialisations of Theorems~\TA\
and \TB. The relevant 
specialisations of our path generating functions from \cite{\CJKrAA,
Eqs.~(2.5)--(2.10)} are the following: with $\om$ denoting a primitive
sixth root of unity, we have
$$\align 
\P_n(k,l)\Big\vert_{x=-y=\sqrt{-1}}&=\chi(n+k+l\text{ even})
\binom n{\frac {1} {2}(n+l-k)},
\tag\SA\\
\P^+_n(k,l)\Big\vert_{x=-y=\sqrt{-1}}&=\chi(n+k+l\text{ even})
\(\binom n{\frac {1} {2}(n+l-k)}-\binom n{\frac {1} {2}(n+l+k+2)}\),
\tag\SB\\ 
\P_n(k,l)\Big\vert_{x=y^{-1}=\om}&=\sum _{\ell\ge0} ^{}
\binom n {\ell,\ell+l-k},
\tag\SC\\ 
\P^+_n(k,l)\Big\vert_{x=y^{-1}=\om}&=\sum _{\ell\ge0} ^{}\(
\binom n {\ell,\ell+l-k}-\binom n {\ell,\ell+l+k+2}\),
\tag\SD\\ 
\P_n(k,l)\Big\vert_{x=y=1}&=
\binom {2n} {n+l-k},
\tag\SE\\ 
\P^+_n(k,l)\Big\vert_{x=y=1}&=
\binom {2n} {n+l-k}-\binom {2n} {n+l+k+2},
\tag\SF
\endalign$$
where $\chi(\Cal A)=1$ if $\Cal A$ is true and  $\chi(\Cal A)$=0
otherwise, and 
$$\binom n{k_1,k_2}=\frac {n!} {k_1!\,k_2!\,(n-k_1-k_2)!}$$
is a trinomial coefficient. We add one more such specialisation,
$$\align
%\P_n(k,l)\Big\vert_{x=\om^2,\,y=\om^{-2}}&=
%(-1)^{n+k+l}\sum _{\ell\ge0} ^{}
%\binom n {\ell,\ell+l-k},
%\tag\SG\\ 
%\P^+_n(k,l)\Big\vert_{x=\om^2,\,y=\om^{-2}}&=
%(-1)^{n+k+l}\sum _{\ell\ge0} ^{}
%\(\binom n {\ell,\ell+l-k}-\binom n {\ell,\ell+l+k+2}\),
%\tag\SH\\ 
%\P_n(k,l)\Big\vert_{x=-y=1}
%&=(-1)^{(n+k-l)/2}\chi(n+k+l\text{ even})\binom n{\frac {1} {2}(n+l-k)},
%\tag\SI\\
%\P^+_n(k,l)\Big\vert_{x=-y=1}
%&=(-1)^{(n+k-l)/2}\chi(n+k+l\text{ even})\\
%&\kern2cm
%\times
%\(\binom n{\frac {1} {2}(n+l-k)}-\binom n{\frac {1} {2}(n+l+k+2)}\),
%\tag\SJ\\
\P_n(k,l)\Big\vert_{x=y=-1}&=
(-1)^{n+k+l}\binom {2n} {n+l-k},
\tag\SK\\ 
\P^+_n(k,l)\Big\vert_{x=y=-1}&=
(-1)^{n+k+l}\(\binom {2n} {n+l-k}-\binom {2n} {n+l+k+2}\).
\tag\SL
\endalign$$
This is easy to derive from \cite{\CJKrAA, Eq.~(2.4)}.
(The specialisations~(\SK) and (\SL) were 
not given in \cite{\CJKrAA} since they do not lead to anything
new in the context of \cite{\CJKrAA}. In our context they do.)
By summation over~$l$ on both sides of (\SB), (\SD), (\SF), and (\SL), we obtain
$$\align 
\sum_{l\ge0}\P^+_n(k,l)\Big\vert_{x=-y=\sqrt{-1}}&=
%\cases
%\dsize
%\sum_{l=0}^k\binom n{\frac {1} {2}(n-k)+l},
%&\text{if }n+k\text{ is even,}\\
%\dsize
%\sum_{l=0}^k\binom n{\frac {1} {2}(n+1-k)+l},
%&\text{if }n+k\text{ is odd,}
%\endcases
\sum_{l=0}^k\binom n{\fl{\frac {1} {2}(n+1-k)}+l},
\tag\SM\\ 
\sum_{l\ge0}\P^+_n(k,l)\Big\vert_{x=y^{-1}=\om}&=\sum _{\ell\ge0} ^{}
\sum_{l=-k}^{k+1}\binom n {\ell,\ell+l},
\tag\SN\\ 
\sum_{l\ge0}\P^+_n(k,l)\Big\vert_{x=y=1}&=
\sum_{l=-k}^{k+1}\binom {2n} {n+l},
\tag\SO\\
%\sum_{l\ge0}\P^+_n(k,l)\Big\vert_{x=\om^2,\,y=\om^{-2}}&=
%\sum _{\ell\ge0} ^{}\sum_{l=-k}^{k+1}(-1)^{n+l}
%\binom n {\ell,\ell+l},\\
%\sum_{l\ge0}\P^+_n(k,l)\Big\vert_{x=-y=1}&=
%\sum_{l\ge0}
%(-1)^{\fl{(n+k)/2}+l}\left(
%\binom n{\fl{\frac {1} {2}(n-k+1)}+l}\right.\\
%&\kern4cm
%\left.-\binom n{\fl{\frac {1} {2}(n+k+3)}+l}\right),
%\\
\sum_{l\ge0}\P^+_n(k,l)\Big\vert_{x=y=-1}&=
\sum_{l=-k}^{k+1}(-1)^{n+l}\binom {2n} {n+l}.
\endalign$$
The last identity can in fact be simplified, in view of the
elementary summation formula
$$
\sum_{s=0}^M(-1)^s\binom Ns=(-1)^M\binom {N-1}M.
$$
Namely, we have
$$\align 
%\sum_{l\ge0}\P^+_n(k,l)\Big\vert_{x=\om^2,\,y=\om^{-2}}&=
%(-1)^{n+k+1}\sum _{\ell\ge0} ^{}
%\binom n\ell
%\(\binom {n-\ell-1}{k+\ell+1}-\binom {n-\ell-1}{\ell-k-1}\)
%\tag\SP\\
%\sum_{l\ge0}\P^+_n(k,l)\Big\vert_{x=-y=1}&=
%(-1)^{\fl{(n+k)/2}-\fl{\frac {1} {2}(n+k+1)}+k}
%\sum_{l\ge\fl{\frac {1} {2}(n-k+1)}}
%(-1)^l\binom n{l}\\
%&-(-1)^{\fl{(n+k)/2}-\fl{\frac {1} {2}(n+k+1)}+1}
%\sum_{l\ge\fl{\frac {1} {2}(n+k+3)}}(-1)^l\binom n{l},
%\\
%&=
%-(-1)^{-\chi(n+k\text{ odd})+k}
%(-1)^{\fl{\frac {1} {2}(n-k-1)}}\binom {n-1}{\fl{\frac {1} {2}(n-k-1)}}\\
%&+(-1)^{-\chi(n+k\text{ odd})+1}
%(-1)^{\fl{\frac {1} {2}(n+k+1)}}\binom {n-1}{\fl{\frac {1} {2}(n+k+1)}},
%\\
%&=\chi(n+k\text{ odd})(-1)^{\frac {1} {2}(n+k+1)}
%\(-\binom {n-1}{\frac {1} {2}(n-k-1)}
%+\binom {n-1}{\frac {1} {2}(n+k+1)}\),
%\\
%&=\chi(n+k\text{ odd})(-1)^{\frac {1} {2}(n+k-1)}
%\frac {k+1} {n}\binom {n}{\frac {1} {2}(n+k+1)},
%\\
%\chi(n+k\text{ odd})(-1)^{\frac {1} {2}(n+k-1)}
%\frac {k+1} {n}\binom {n}{\frac {1} {2}(n+k+1)},
%\tag\SQ\\
\sum_{l\ge0}\P^+_n(k,l)\Big\vert_{x=y=-1}&=
%\sum_{l=n-k}^{n+k+1}(-1)^{l}\binom {2n} {l}
%=
%(-1)^{n+k+1}\binom {2n-1} {n+k+1}
%-(-1)^{n-k-1}\binom {2n-1} {n-k-1}\\
%&=
(-1)^{n+k}\frac {k+1} {n}\binom {2n}{n+k+1}.
\tag\SR
\endalign$$
Care must be applied of how to interpret
the expression on the right-hand side for $n=0$: in (\SR),
the value for $n=0$ must be taken as~$1$, regardless of the choice of~$k$.

\subhead 3. The connection matrix $A(n)$\endsubhead
In this section, we define the {\it connection matrix} $A(n)$ announced in the
introduction, see (\AC). 
Multiplication of our matrices of Motzkin prefix generating
functions on the left by $A(n)$
allows us to connect them --- via Lemmas~\TD\ and \TI\ --- 
to the matrices of Motzkin generating functions
in \cite{\CJKrAA}, presented here in Theorems~\TCKA\ and \TCKB\
in the previous section.
The connection is not completely direct, it is only up to {\it correction
matrices} (the matrices $C_0(n,k)$ and $C_1(n,k)$ in the lemmas).
Nevertheless, since $A(n)$ has determinant~$1$ (see Lemma~\TC), 
multiplication on the left by $A(n)$ does not change the determinant,
and some further work makes it possible to deduce Theorems~\TA\
and \TB\ on the basis of results from Theorems~\TCKA\ and \TCKB.

\medskip
We define the matrix $A(n):=(A_{n,i,j})_{0\le i,j\le n-1}$ by
$$
A_{n,i,j}=\cases 
\frac {(1 + x) (1 + y)}{xy},&\text{if }i=j<n-1,\\
-\frac {1} {xy},&\text{if }i=j-1<n-1,\\
\frac {(-1)^{n + j}} {xy}
\sum_{l =j}^{n}\left(\binom l  j \binom {n + j - 1 - l } { j} x^{l  - j} y^{n - 1 - l } \right.\\
\kern3cm\left.
+ 
         \binom {l } { j} \binom {n + j - l } { j} x^{l  - j} y^{n - l }\right),
&\text{if }i=n-1\text{ and }j<n-1,\\
\frac { x y - (n - 1) (x + y)} {xy},&\text{if }i=j=n-1,\\
0,&
\text{otherwise.}\endcases
\tag\AC
$$
Here, binomial coefficients have to be interpreted as~$0$ as soon as a lower
parameter becomes negative or an upper parameter is less than the
lower parameter. 
For example, the matrix $A(4)$ has the form
$$\pmatrix 
\frac {(1 + x) (1 + y)}{xy}&-\frac {1} {xy}&0&0\\
0&\frac {(1 + x) (1 + y)}{xy}&-\frac {1} {xy}&0\\
0&0&\frac {(1 + x) (1 + y)}{xy}&-\frac {1} {xy}\\
A_{4,3,0}&A_{4,3,1}&A_{4,3,2}&    \frac {x y - 3(x + y)} {xy}
\endpmatrix,
$$
where the entries $A_{4,3,0},A_{4,3,1},A_{4,3,2}$ are the polynomials in
$x$ and $y$ divided by $xy$ given by the next-to-last line in (\AC).

In the proof of Lemma~\TG, we shall need
alternative formulae for the matrix entries in
the last row of $A(n)$, which are presented in the lemma below.

\proclaim{Lemma \TF}
For all non-negative integers $n$ and $j$ with $0\le j\le n-2$, we have
$$\multline
A_{n,n-1,j}=\frac {(-1)^{n + j}} {xy}
\sum_{r\ge0}(-1)^r\left( \binom {n-r-1}r \binom {n-2r-1}j
(xy)^r (x+y)^{n -1- j - 2r}\right.\\
\left. +\binom {n-r}r \binom {n-2r}j
(xy)^r (x+y)^{n - j - 2r}\right).
\endmultline
\tag\AK$$
\endproclaim

\demo{Proof}
Since the two terms in the summand on the right-hand side of (\AK)
arise from each other by a shift of $n$ by~$1$, it suffices to
concentrate on one of them:
$$\align
\sum_{r\ge0}(-1)^r&
 \binom {n-r}r \binom {n-2r}j
(xy)^r (x+y)^{n - j - 2r}\\
&=
\sum_{r\ge0}(-1)^r
 \binom {n-r}r \binom {n-2r}j
(xy)^r \sum_{\ell\ge0}\binom {n - j - 2r}\ell x^\ell y^{n-j-2r-\ell}\\
&=
\sum_{l \ge0}x^{l } y^{n-j-l }\sum_{r\ge0}(-1)^r
 \binom {n-r}r \binom {n-2r}j
 \binom {n - j - 2r}{l -r}\\
&=
\sum_{l \ge0}x^{l } y^{n-j-l }
 \binom {n}j
 \binom {n - j }{l }
{} _{2} F _{1} \!\left [ \matrix 
j+l -n,-l \\ -n\endmatrix ; {\displaystyle
   1}\right ].
\endalign$$
Here, we used the standard hypergeometric notation
$$
{}_r F_s\!\left[\matrix a_1,\dots,a_r\\ b_1,\dots,b_s\endmatrix;  
z\right]=\sum _{l =0} ^{\infty}\frac {\po{a_1}{l }\cdots\po{a_r}{l }} 
{l !\,\po{b_1}{l }\cdots\po{b_s}{l }} z^l \ , 
$$
where the {\it Pochhammer symbol\/} $(\al)_m$ is defined by
$(\alpha)_m  =  \alpha(\hbox{$\alpha+1$})(\alpha+2)\cdots (\alpha+m-1)$ for
$m>0$, and $(\alpha)_0  = 1$.
The above $_2F_1$-series can be evaluated by means of the
Chu--Vandermonde summation formula
(see \cite{\SlatAC, (1.7.7); Appendix (III.4)}),
$$
{} _{2} F _{1} \!\left [ \matrix { a, -N}\\ { c}\endmatrix ; {\displaystyle
   1}\right ]  = {{({ \textstyle c-a}) _{N} }\over
    {({ \textstyle c}) _{N} }},
\tag\AH$$
where $N$ is a nonnegative integer.
Thus, we obtain
$$\align
\sum_{l \ge0}x^{l } y^{n-j-l }
 \binom {n}j
 \binom {n - j }{l }
\frac {(-j-l )_l } {(-n)_l }
&=
\sum_{l \ge0}x^{l } y^{n-j-l }
 \binom {n-l }j
 \binom {l +j }{j}\\
&=
\sum_{l =j}^{n}x^{l -j} y^{n-l }
 \binom {n-l +j}j
 \binom {l  }{j}.
\endalign$$
This matches exactly with the original definition of
$A_{n,n-1,j}$.\quad \quad \qed
\enddemo

Next we show that the determinant of the connection matrix $A(n)$ is
equal to~$1$. The proof requires an identity that is established
separately in Lemma~\TEa.

\proclaim{Lemma \TC}
For all positive integers $n$, we have $\det A(n)=1$.
\endproclaim

\demo{Proof}
We replace the last column of $A(n)$ by
$$
\sum_{j=0}^{n-1}\frac {1} {(1+x)^{n-1-j}(1+y)^{n-1-j}}\cdot
(\text{column }j).
$$
Clearly, this does not change the value of the determinant of $A(n)$.
Moreover, in the resulting matrix, all entries in the last column
will become~$0$, except for the entry in the last row, which equals
$$
\sum_{j=0}^{n-2}\frac {1} {(1+x)^{n-1-j}(1+y)^{n-1-j}}A_{n,n-1,j}
+\frac { x y - (n - 1) (x + y)} {xy}.
$$
By Lemma~\TEa, this expression equals
$$
\left(\frac {xy} {(1+x)(1+y)}\right)^{n-1}.
\tag\AHa
$$

We may now expand the determinant with respect to the last column.
Since in this column only the entry in the last row is non-zero,
the result is that entry --- that is, the expression in (\AHa) ---
multiplied by the remaining minor $\det_{0\le i,j\le n-2}(A_{n,i,j})$.
This latter determinant is the determinant of a lower triangular matrix. 
Obviously, it equals the product of the diagonal entries. Altogether
we obtain~$1$ as is straightforward to see. This establishes the 
lemma.\quad \quad \qed
\enddemo

\proclaim{Lemma \TEa}
For all non-negative integers $m$, we have
$$
\sum_{j=0}^{n-2}\frac {1} {(1+x)^{n-1-j}(1+y)^{n-1-j}}A_{n,n-1,j}
+\frac { x y - (n - 1) (x + y)} {xy}
=
\left(\frac {xy} {(1+x)(1+y)}\right)^{n-1}.
\tag\AD
$$
\endproclaim

\demo{Proof}
We start with the precise form of the left-hand side of (\AD),
$$\multline 
\frac {1} {xy}
\sum_{j=0}^{n-2}\frac {(-1)^{n + j}} {(1+x)^{n-1-j}(1+y)^{n-1-j}}
\sum_{l =j}^{n}\left(\binom l  j \binom {n + j - 1 - l } { j} x^{l  - j} y^{n - 1 - l } \right.\\
\kern3cm\left.
+ 
         \binom {l } { j} \binom {n + j - l } { j} x^{l  - j} y^{n - l }\right)
+\frac { x y - (n - 1) (x + y)} {xy}.
\endmultline$$
It is straightforward to see that, by extending the sum over $j$ 
to the range $0\le j\le n$, the last term in the above expression gets
``swallowed". In other terms, the expression can be rewritten as
$$\align
\frac {1} {xy}&
\sum_{j=0}^{n}\frac {(-1)^{n + j}} {(1+x)^{n-1-j}(1+y)^{n-1-j}}
\sum_{l =j}^{n}\left(\binom l  j \binom {n + j - 1 - l } { j} x^{l  - j} y^{n - 1 - l } \right.\\
&\kern3cm\left.
+ 
         \binom {l } { j} \binom {n + j - l } { j} x^{l  - j} y^{n - l }\right)\\
&=
\frac {(-1)^n} {xy(1+x)^{n-1}(1+y)^{n-1}}\\
&\kern.5cm
\times\left(
\sum_{j=0}^{n-1}(-1)^{j}
\sum_{l =j}^{n-1}\binom l  j \binom {n + j - 1 - l } { j} x^{l  - j} y^{n
- 1 - l } (1+x)^{j}(1+y)^{j}\right.\\
&\kern2cm
\left.
+ 
\sum_{j=0}^{n}(-1)^{j}
\sum_{l =j}^{n}
         \binom {l } { j} \binom {n + j - l } { j} x^{l  - j} y^{n - l }
(1+x)^{j}(1+y)^{j}\right).
\tag\AF
\endalign$$
Here, to lower the upper bounds on the summation indices of the sums
over~$j$ and~$l $ in the first double sum is allowed due to the 
vanishing properties of the binomial coefficients $\binom l j$
and $\binom {n+j-1-l }{j}$. The purpose of this ``exercise" is to
make it visible that the first double sum arises from the
second by replacing $n$ by $n-1$.

Hence it suffices to concentrate on the second double sum.
The coefficient of $x^Ay^B$, $0\le A,B\le n$, in this sum is given by
$$\align
\sum_{j=0}^{n}(-1)^{j}&
\sum_{l =j}^{n}
         \binom {l } { j} \binom {n + j - l } { j}
\binom j{A-l +j}\binom j{B-n+l }\\
&=
\sum_{j=0}^{n}(-1)^{j}
\sum_{l =0}^{n-j}
         \binom {l +j} { j} \binom {n - l } { j} 
\binom j{A-l }\binom j{B-n+l +j}\\
&=
\sum_{l =0}^{n}
\sum_{j=A-l }^{n}(-1)^{j}
         \binom {l +j} { j} \binom {n - l } { j} 
\binom j{A-l }\binom j{B-n+l +j}.
\endalign$$
We write the inner sum over~$j$ in hypergeometric notation.
Thereby we obtain
$$
\sum_{l =0}^{n}
(-1)^{A-l }
         \binom {A} { l } \binom {n - l } { A-l } 
\binom j{A-l }\binom j{A+B-n}
{}_2 F_1\!\left[\matrix A+1,A-n\\ A+B-n+1\endmatrix;  
1\right].
$$
The $_2F_1$-series can again be evaluated by means of the Chu--Vandermonde
summation formula (\AH).
We substitute the result and now write the remaining sum over~$l $
in hypergeometric notation. Thus, we arrive at
$$
(-1)^n\binom {n-B}{n-A}\binom nB
{} _{2} F _{1} \!\left [ \matrix B-n,-A\\ -n\endmatrix ; {\displaystyle
   1}\right ].
$$
By applying (\AH) once again and some simplification, we finally get
$$
(-1)^{A+B+n}\frac {(-B)_A\,(-A)_B} {A!\,B!}
=\cases (-1)^n,&\text{if }A=B,\\
0,&\text{otherwise,}\endcases
\tag\AI
$$
for the second double sum in (\AF).
As we discussed earlier, the first double arises from the second by
replacing $n$ by $n-1$. Thus, these two double sums either both
vanish or cancel each other, except for $A=B=n$; in that
latter case we are asking for the coefficient of $x^Ay^B=x^ny^n$,
which is necessarily zero in the first double sum (because only
monomials of lower degree can appear) while it is the $(-1)^n$
from (\AI) for the second double sum that survives. If this is substituted in
(\AF), the assertion (\AD) follows immediately, which completes
the proof of the lemma.\quad \quad \qed
\enddemo

The next two lemmas provide the identities that are crucial for 
establishing the link between the matrices in Theorems~\TA\ and \TB\
and the matrices in \cite{\CJKrAA}, made explicit in Lemmas~\TD\
and \TI. The first is a relatively simple combinatorial identity.

\proclaim{Lemma \TE}
For all non-negative integers $m$, we have
$$
(1+x)(1+y)\sum_{l =0}^{m}\P^+_{m}(k ,l )
-\sum_{l =0}^{m+1}\P^+_{m+1}(k ,l )=
xy\P^+_{m}(k ,0).
\tag\AJ
$$
\endproclaim

\demo{Proof}
We have
$$
(1+x)(1+y)=1+(x+y)+xy.
$$
This is exactly the sum of the weights of an up-step, of a horizontal
step, and of a down-step. Thus, the first term in (\AJ) is the
generating function for paths consisting of $m+1$ (horizontal, up-
and down-)steps that start at $(0,k )$ and do not run below the
$x$-axis {\it for the first $m$ steps}. On the other hand, the second term
is the negative of the generating function for the same paths, except
that one requires the stronger condition that they do not run below
the $x$-axis {\it for all of their $m+1$ steps}. Hence, the
difference on the left-hand side of (\AJ) equals the generating
function for all those paths which reach $(m,0)$ without having passed
below the $x$-axis, but then continue with a down-step. Since a
down-step has weight $xy$, this is exactly the expression on the
right-hand side of (\AJ).\quad \quad \qed
\enddemo

The second identity is not combinatorial (at least, the authors do
not have a combinatorial interpretation for it). Instead, its proof
requires a certain summation formula for hypergeometric series which
is stated separately in Lemma~\TH.

\proclaim{Lemma \TG}
For non-negative integers $m$ and positive integers $n$ 
with $0\le m\le n$, we have
$$
\sum_{j=0}^{n-1}A_{n,n-1,j}
\sum_{l =0}^{m+j}\P^+_{m+j}(k ,l )
=
\P^+_{m+n-1}(k ,0)+(xy)^{n-1}\sum_{l\ge0}^{}\P_m(0,n-k+l),
%\P^+_{m+n-1}(k ,0)\\
%+
%\sum_{t\ge0}\sum_{l =0}^{m+k }
%\frac {m!} {t!\,(l -k +t+n)!\,(m-n-l +k -2t)!}
%%\binom m{l +n-k +2t}
%%\binom {l -k +2t+n} {t}
%(xy)^{t+n-1} (x+y)^{m-n-l +k - 2t},
\tag\AKa$$
where the coefficients $A_{n,n-1,j}$ are given in {\rm(\AC)}.
\endproclaim

\demo{Proof}
Using the expression for $A_{n,n-1,j}$ for $0\le j\le n-2$ from Lemma~{\TF} and
the expression for $\P^+_{m+j}(k ,l )$ in (\AAc), we have
$$\align 
\sum_{j=0}^{n-1}&A_{n,n-1,j}
\sum_{l =0}^{m+j+k }\P^+_{m+j}(k ,l )
\\
&=
\sum_{j=0}^{n-2}\sum_{r,s\ge0}\sum_{l =0}^{m+j+k }
\frac {(-1)^{n + j+r}} {xy}
\binom {m+j}{l -k +2s}
\left(\binom {l -k +2s}s - \binom {l -k +2s}{s-k -1}\right)\\
&\kern1cm
\cdot
\left( \binom {n-r-1}r \binom {n-2r-1}j
(xy)^{r+s} (x+y)^{n -1 +m-l +k - 2r-2s}\right.\\
&\kern3cm
\left. +\binom {n-r}r \binom {n-2r}j
(xy)^{r+s} (x+y)^{n +m-l +k - 2r-2s}\right)\\
&\kern.5cm
+\frac { x y - (n - 1) (x + y)} {xy}\sum_{l =0}^{m+n+k-1}\P^+_{m+n-1}(k ,l ).
\endalign$$
We would like to extend the sum over~$j$ to range over $0\le j\le n$.
Because of the binomial coefficients $\binom {n-2r-1}j$ and $\binom
{n-2r}j$, this extension is indeed without any harm except if $r=0$.
Taking the corresponding corrections into account,
we see that the above expression is equal to
$$\align 
\sum_{j=0}^{n}&\sum_{r,s\ge0}\sum_{l =0}^{m+j+k }
\frac {(-1)^{n + j+r}} {xy}
\binom {m+j}{l -k +2s}
\left(\binom {l -k +2s}s - \binom {l -k +2s}{s-k -1}\right)\\
&\kern1cm
\cdot
\left( \binom {n-r-1}r \binom {n-2r-1}j
(xy)^{r+s} (x+y)^{n -1 +m-l +k - 2r-2s}\right.\\
&\kern3cm
\left. +\binom {n-r}r \binom {n-2r}j
(xy)^{r+s} (x+y)^{n +m-l +k - 2r-2s}\right)\\
&
+\sum_{s\ge0}\sum_{l =0}^{m+n+k -1}
\frac {1} {xy}
\binom {m+n-1}{l -k +2s}
\left(\binom {l -k +2s}s - \binom {l -k +2s}{s-k -1}\right)\\
&\kern1cm
\cdot
\left(  (xy)^{s} (x+y)^{n -1 +m-l +k -2s}
 + n
(xy)^{s} (x+y)^{n +m-l +k -2s}\right)\\
&
-\sum_{s\ge0}\sum_{l =0}^{m+n+k }
\frac {1} {xy}
\binom {m+n}{l -k +2s}
\left(\binom {l -k +2s}s - \binom {l -k +2s}{s-k -1}\right)\\
&\kern2cm
\cdot
 (xy)^{s} (x+y)^{n +m-l +k -2s}\\
&
+\frac { x y - (n - 1) (x + y)} {xy}\sum_{l =0}^{m+n+k-1}\P^+_{m+n-1}(k ,l ).
\endalign$$
We have
$$\align
\sum_{j=0}^n(-1)^j\binom {m+j}{l -k +2s}\binom {n-2r}j
&=\binom {m}{l -k +2s}
{} _{2} F _{1} \!\left [ \matrix 
m+1,2r-n\\ m-l +k -2s+1\endmatrix ; {\displaystyle
   1}\right ]\\
&=
\binom {m}{l -k +2s}\frac {(-l +k -2s)_{n-2r}} {(m-l +k -2s+1)_{n-2r}}\\
&=
(-1)^n\binom m{l -n-k +2s+2r},
\endalign$$
again by the Chu--Vandermonde summation formula (\AH). If this is
substituted (twice --- once with $n$ replaced by~$n-1$) in the 
quadruple sum of our earlier obtained expression, then we get
$$\align 
\sum_{r,s\ge0}&\sum_{l =0}^{m+k }
\frac {(-1)^{r}} {xy}
\left(\binom {l -k +2s}s - \binom {l -k +2s}{s-k -1}\right)\\
&\kern1cm
\cdot
\left( -\binom {n-r-1}r \binom m{l -n-k +2s+2r+1}
(xy)^{r+s} (x+y)^{n -1 +m-l +k - 2r-2s}\right.\\
&\kern3cm
\left. +\binom {n-r}r \binom m{l -n-k +2s+2r}
(xy)^{r+s} (x+y)^{n +m-l +k - 2r-2s}\right)\\
&
+\frac {\left(  1 + n (x+y)\right)} {xy}
\sum_{l =0}^{m+n+k -1}\P^+_{m+n-1}(k ,l )
-\frac {1} {xy}\sum_{l =0}^{m+n+k }\P^+_{m+n}(k ,l )\\
&
+\frac { x y - (n - 1) (x + y)} {xy}\sum_{l =0}^{m+n+k-1}\P^+_{m+n-1}(k ,l ).
\tag\AL
\endalign$$
We observe that
$$
\big(1+n(x+y)\big)+\big(x y - (n - 1) (x + y)\big)=(1+x)(1+y).
$$
Hence, by Lemma~\TE,
the last two lines of (\AL) simplify to $\P^+_{m+n-1}(k ,0)$,
which is the first term on the right-hand side of (\AKa).

The remaining task is therefore to simplify the triple sum in (\AL).
In order to do so, we split it into two parts,
$$\align 
S_1&=\sum_{r,s\ge0}\sum_{l =0}^{m+k }
\frac {(-1)^{r}} {xy}
\binom {l -k +2s}s \\
&\kern1cm
\cdot
\left( -\binom {n-r-1}r \binom m{l -n-k +2s+2r+1}
(xy)^{r+s} (x+y)^{n -1 +m-l +k - 2r-2s}\right.\\
&\kern3cm
\left. +\binom {n-r}r \binom m{l -n-k +2s+2r}
(xy)^{r+s} (x+y)^{n +m-l +k - 2r-2s}\right)\\
&=\sum_{r,s\ge0}\sum_{l =0}^{m+k }
\frac {(-1)^{r}} {xy}
\binom m{l -n-k +2s+2r}
(xy)^{r+s} (x+y)^{n +m-l +k - 2r-2s}
\\
&\kern2cm
\cdot
\left(\binom {l -k +2s}s \binom {n-r}r 
 -\binom {l -k +2s-1}s \binom {n-r-1}r \right),
\endalign$$
and
$$\align 
S_2&=-\sum_{r,s\ge0}\sum_{l =0}^{m+k }
\frac {(-1)^{r}} {xy}
\binom {l -k +2s}{s-k -1} \\
&\kern1cm
\cdot
\left( -\binom {n-r-1}r \binom m{l -n-k +2s+2r+1}
(xy)^{r+s} (x+y)^{n -1 +m-l +k - 2r-2s}\right.\\
&\kern3cm
\left. +\binom {n-r}r \binom m{l -n-k +2s+2r}
(xy)^{r+s} (x+y)^{n +m-l +k - 2r-2s}\right)\\
&=-\sum_{r,s\ge0}\sum_{l =0}^{m+k }
\frac {(-1)^{r}} {xy}
\binom m{l -n-k +2s+2r}
(xy)^{r+s} (x+y)^{n +m-l +k - 2r-2s}
\\
&\kern2cm
\cdot
\left(\binom {l -k +2s}{s-k -1} \binom {n-r}r 
 -\binom {l -k +2s-1}{s-k -1} \binom {n-r-1}r \right).
\endalign$$

We start with the computation of $S_1$. We let $t=r+s$ and write the
sum over~$r$ in hypergeometric notation. This leads to
$$\align 
S_1
&=\sum_{t\ge0}\sum_{l =0}^{m+k }
\frac {1} {xy}
\binom m{l -n-k +2t}
(xy)^{t} (x+y)^{n +m-l +k - 2t}
\\
&\kern2cm
\cdot
\binom {l -k +2t-1}{t-1}
{} _{5} F _{4} \!\left [ \matrix 
1+\frac {n t} { l  - k  - n} ,
 -\frac n2, \frac 12 - \frac n2 ,  - t ,  -l  + k  - t
\\ 
\frac {n t} { l  - k  - n} , 1 - n, \frac 12 - \frac l 2 + \frac k 2 - t ,
1 - \frac l 2 + \frac k 2 - t 
\endmatrix ; {\displaystyle
   1}\right ].
\endalign$$
Next, we apply the contiguous relation
$$\multline
{} _{5} F _{4} \!\left [ \matrix { a, b, c, A_1,A_2}\\ { B_1,B_2,B_3,B_4}\endmatrix ;
   {\displaystyle z}\right ] =
  {{b \left( c - a - 1  \right)  
        }\over 
     {\left( b - a \right)  \left( c - 1 \right) }}
  {} _{5} F _{4} \!\left [ \matrix { a, b + 1, c - 1, A_1,A_2}\\ {
        B_1,B_2,B_3,B_4}\endmatrix ; {\displaystyle z}\right ] \\
+ 
   {{a \left( c - b - 1  \right)  
        }\over 
     {\left( a - b \right)  \left( c - 1 \right) }}
  {} _{5} F _{4} \!\left [ \matrix { a + 1, b, c - 1, A_1,A_2}\\ {
        B_1,B_2,B_3,B_4}\endmatrix ; {\displaystyle z}\right ]
\endmultline
$$
with $a=-t$, $b=-l +k -t$, $c=1+\frac {nt} {l -k -n}$, $A_1=-\frac {n}
{2}$, $A_2=\frac {1} {2}-\frac {n} {2}$, $B_1=\frac {nt} {l -k -n}$,
$B_2=1-n$, $B_3= \frac 12 - \frac l 2 + \frac k 2 - t $, and
$B_4= 1 - \frac l 2 + \frac k 2 - t $. Since, with our choice, we have
$c-1=B_1$, the effect is that, on the right-hand side, the
$_5F_4$-series reduce to $_4F_3$-series. Thus, we obtain
$$\align 
S_1
&=\sum_{t\ge0}\sum_{l =0}^{m+k }
\frac {1} {xy}
\binom m{l -n-k +2t}
\binom {l -k +2t-1}{t-1}
(xy)^{t} (x+y)^{n +m-l +k - 2t}
\\
&\kern1cm
\cdot
\left(
\frac {l -k +t} {n}
{} _{4} F _{3} \!\left [ \matrix 
 -\frac n2, \frac 12 - \frac n2 ,  - t ,  1-l  + k  - t
\\ 
1 - n, \frac 12 - \frac l 2 + \frac k 2 - t ,
1 - \frac l 2 + \frac k 2 - t 
\endmatrix ; {\displaystyle
   1}\right ]\right.\\
&\kern2cm
\left.
-\frac {l -k +t-n} {n}
{} _{4} F _{3} \!\left [ \matrix 
 -\frac n2, \frac 12 - \frac n2 ,  1- t ,  -l  + k  - t
\\ 
1 - n, \frac 12 - \frac l 2 + \frac k 2 - t ,
1 - \frac l 2 + \frac k 2 - t 
\endmatrix ; {\displaystyle
   1}\right ]\right).
\endalign$$
Both $_4F_3$-series can be evaluated by means of Lemma~\TH.
After some simplification, the result is
$$\align 
S_1
&=\frac {1} {xy}
\sum_{t\ge0}\sum_{l =0}^{m+k }
\binom m{l -n-k +2t}
\binom {l -k +2t-n} {t-n}
(xy)^{t} (x+y)^{n +m-l +k - 2t}\\
&=
\sum_{t\ge0}\sum_{l =0}^{m+k }
\frac {m!} {t!\,(l -k +t+n)!\,(m-n-l +k -2t)!}
(xy)^{t+n-1} (x+y)^{m-n-l +k - 2t}\\
&=
(xy)^{n-1}\sum_{l\ge0}^{}\P_m(0,n-k+l).
\endalign$$
Here, we replaced $t$ by $t+n$ to go from the first to the second line,
and subsequently we used (\AAb) to arrive at the last line.
Clearly, this is the second term on the right-hand side of (\AKa).

A similar computation yields
$$\align 
S_2
&=\frac {1} {xy}
\sum_{t\ge0}\sum_{l =0}^{m+k }
\binom m{l -n-k +2t}
\binom {l -k +2t-n} {t-n-k -1}
(xy)^{t} (x+y)^{n +m-l +k - 2t}.
\endalign$$
Due to the binomial coefficient $\binom {l -k +2t-n} {t-n-k -1}$, the
summation index $t$ must be at least $n+k+1$ in order to generate
non-vanishing summands. However, in that case we have 
$$l-n-k+2t\ge
l+n+k+2\ge l+m+k+2>m,
$$
which makes the binomial coefficient $\binom m{l -n-k +2t}$ vanish.
In other words, we have $S_2=0$. This completes the proof of the
lemma.\quad \quad \qed 
\enddemo

The following is Lemma~A3 from \cite{\KratBD}.

\proclaim{Lemma \TH}Let $n$ be a positive integer. Then
$${}_4F_3\!\[\matrix -\frac {n} {2},\frac {1} {2}-\frac {n} {2},-A,A+B\\
1-n,\frac {B} {2},\frac {1} {2}+\frac {B} {2}\endmatrix; 1\]=\frac {(A+B)_n} {(B)_n}+
\frac {(-A)_n} {(B)_n}.$$
\endproclaim

\subhead 4. Proof of Theorem \TA\endsubhead
We first use the results from the previous section to connect the
matrix of Motzkin prefix generating functions on the left-hand side of
(\AA) to a matrix of Motzkin generating functions that appeared
in \cite{\CJKrAA}.

\proclaim{Lemma \TD}
Define matrices $M_0(n,k )$ and $\MP_0(n,k )$ by
$$
M_0(n,k ):=\left(\P^+_{i+j}(k ,0)\right)_{0\le i,j\le n-1}
$$
and
$$
\MP_0(n,k):=\left(\sum_{l \ge0}\P^+_{i+j}(k ,l )\right)_{0\le i,j\le n-1}.
$$
Then
$$
A(n)\cdot \MP_0(n,k )= M_0(n,k )+C_0(n,k ),
$$
where $A(n)$ is given by {\rm(\AC)} and the ``correction matrix" 
$C_0(n,k):=(C^{(0)}_{n,k,i,j})_{0\le i,j\le n-1}$ is defined via
$$
C^{(0)}_{n,k,i,j}=\cases 
0,&\text{if }i\le n-2,\\
(xy)^{n-1}\sum_{l\ge0}^{}\P_j(0,n-k+l),&\text{if }i=n-1.
\endcases
$$
\endproclaim

\demo{Proof}
This is a direct consequence of Lemmas~\TE\ and \TG.\quad \quad \qed
\enddemo

We are now in the position to prove Theorem~\TA.

\demo{Proof of Theorem \TA}
We start with 
$$
\MP_0(n,k)=\left(\sum_{l \ge0}\P^+_{i+j}(k ,l )\right).
$$
We multiply on the left by $A(n)$. According to Lemma~\TD, we get
$$\multline
M_0(n,k)+C_0(n,k)\\
=\left(
\cases 
\P^+_{i+j}(k ,0),&\text{for }i\le n-2\\
\P^+_{n+j-1}(k ,0)+
(xy)^{n-1}\sum_{l\ge0}^{}\P_j(0,n-k+l),&\text{for }i=n-1
\endcases
\right).
\endmultline$$
Since $\det A(n)=1$ by Lemma~\TC, the determinant of $M_0(n,k)+C_0(n,k)$
is the same as the determinant of $\MP_0(n,k)$. If we apply relation
(\AAa) with $l=0$, then we see that we have transformed our problem
into the problem of evaluation of the determinant of
$$(xy)^kM'_0(n,k)+C_0(n,k),$$ 
where
$$
M'_0(n,k ):=\left(\P^+_{i+j}(0,k)\right)_{0\le i,j\le n-1}.
$$

The determinant of $M'_0(n,k)$ is evaluated in Theorem~\TCKA. 
As it turns out, for a while we may now follow the
arguments of the proof of this evaluation in \cite{\CJKrAA}. For the
convenience of the reader, we summarise the main steps here.

\medskip
The first step in \cite{\CJKrAA} makes use of the combinatorics of
non-intersecting lattice paths, see Section~4 there. One may however
do equally well without combinatorics, as we now explain. 
By cutting paths after $i$ steps, it is easy to see that the equation
$$\P^+_{i+j}(0,k)=
\sum _{\ell=0} ^{i}\P^+_i(0,\ell)\P^+_j(\ell,k)
\tag\ALa$$
holds. Thus, we see that $(xy)^kM'_0(n,k)+C_0(n,k)$ is equal to the product
of the matrices
$$
\left(\P^+_i(0,\ell)\right)_{0\le i,\ell\le n-1}
\tag\AM$$
and
$$
(xy)^k\left(\P^+_j(\ell,k)\right)_{0\le\ell,j\le n-1} + C_0(n,k).
\tag\AN$$
Indeed, since $\P^+_i(0,\ell)=0$ for $i<\ell$ and $\P^+_i(0,i)=1$,
the matrix in (\AM) is lower triangular with $1$'s on the main diagonal,
and thus we have
$$
\left(\P^+_i(0,\ell)\right)_{0\le i,\ell\le n-1}\cdot C_0(n,k)=C_0(n,k).
$$
Moreover, for the same reason the determinant of the matrix in (\AM)
is~$1$. Hence, the determinant of (\AN) still equals $\det \MP_0(n,k)$.

\medskip
The second step in \cite{\CJKrAA} (see the first paragraph of
Section~5 there) consists in the use of (\AG) in
order to rewrite $\P^+_j(\ell,k)$. In our case, we are led to the
problem of evaluating the determinant of
$$
(xy)^k\left(
\P_j(i,k)-(xy)^{\ell+1}\P_j(-i-2,k)\right)_{0\le i,j\le n-1} + C_0(n,k).
\tag\AO$$

\medskip
In the third step in \cite{\CJKrAA} (see Eqs.~(5.6) and (5.7) there
with $t=1$),
certain row operations are applied. To be precise, 
row $(h(2k+2)+b)$ of the matrix obtained so far gets replaced by
$$\multline 
\sum _{\ell=0} ^{h}(xy)^{(h-\ell)(k+1)}\cdot
\big(\text {row $(\ell(2k+2)+b)$}\big)\\
-\sum _{\ell=1} ^{h}(xy)^{(h-\ell)(k+1)+b+1}\cdot
\big(\text {row $(\ell(2k+2)-b-2)$}\big)
\endmultline\tag\AOa$$
if $0\le b\le k-1$, and by
$$\multline 
\sum _{\ell=0} ^{h}(xy)^{(h-\ell)(k+1)}\cdot
\big(\text {row $(\ell(2k+2)+b)$}\big)\\
-\sum _{\ell=1} ^{h+1}(xy)^{(h-\ell)(k+1)+b+1}\cdot
\big(\text {row $(\ell(2k+2)-b-2)$}\big)
\endmultline\tag\AOb$$
if $k+1\le b\le 2k$. 
We apply these same operations to the matrix in (\AO).
The important feature of these operations is that, to obtain a row
of the new matrix, only this row and {\it earlier} rows are involved.
Therefore, by the computations performed in \cite{\CJKrAA, proof of
Theorem~8 with $t=1$} that finally lead to (5.8) there and the subsequent two
displays, these operations transform the matrix in (\AO) into the matrix
$$
(xy)^kN_0(n,k)+C_0(n,k),
\tag\AP$$
where $N_0(n,k)=(N^{(0)}_{n,k,i,j})_{0\le i,j\le n-1}$ with
$$
N^{(0)}_{n,k,h(2k+2)+b,j}=\cases 
-(xy)^{(2h+1)(k+1)+b-k}\P_j(0,(h+1)(2k+2)+b-k)\\
\quad +(xy)^{2h(k+1)}\P_j(0,h(2k+2)-b+k),&\text{if }0\le b\le k,\\
-(xy)^{(2h+1)(k+1)+b-k}\P_j(0,(h+1)(2k+2)+b-k)\\
\quad +(xy)^{(2h+1)(k+1)}\P_j(0,(h+1)(2k+2)-b+k),\\
&\kern-1.5cm\text{if }k+1\le b\le 2k+1.
\endcases
$$

\medskip
Close inspection of the new matrix in (\AP) reveals that its determinant can
now rather straightforwardly be deduced.

\medskip
{\smc Case 1:} \hbox{$n\not\equiv0,1$ (mod $k+1$).}
Let
$$n=H(2k+2)+B
\tag\APa$$ 
with $0\le B\le 2k+1$ but $B\ne 0,1,k+1,k+2$.
Then it is not difficult to see (see the paragraphs after (5.8)
in \cite{\CJKrAA}) that, if $1\le B\le k$, row $H(2k+2)$ of $N_0(n,k)$
consists entirely of zeroes, while, if $k+2\le B\le 2k+1$, 
row $H(2k+2)+k+1$ consists entirely of zeroes. In particular, this
implies that in our case there is a row of zeroes in the
matrix in (\AP), and hence its determinant vanishes. This establishes
the third case on the right-hand side of (\AA).

\medskip
{\smc Case 2:} $n\equiv1$ (mod $k+1$).
With the notation of (\APa), we have $B=1$ or $B=k+2$. 
By reusing the arguments in Case~1, we see that, here, it is the last row
of $N_0(n,k)$ (namely row $n-1=H(2k+2)+B-1$) which consists entirely of zeroes.
Since, in (\AP), the matrix $C_0(n,k)$ --- which is a matrix with
potentially non-zero entries in the last row --- is added to
$(xy)^kN_0(n,k)$, we cannot conclude that the determinant of (\AP)
vanishes, but rather further analysis is required.

From now on, let $n=(k+1)n_1+1$.
In order to get a clearer picture, it is convenient to reverse the
order of rows
$s(k+1),s(k+1)+1,\dots,s(k+1)+k$, for $s=0,1,\dots,n_1-1$.
It should be noticed that this leaves the last row, namely row
$n-1=(k+1)n_1$, in place. In this manner, we arrive at the matrix
$$
(xy)^k{\bar N}_0(n,k)+C_0(n,k),
\tag\AQ$$
where the matrix ${\bar N}_0(n,k)=({\bar N}^{(0)}_{n,k,i,j})_{0\le i,j\le n-1}$ is given
by
$$
{\bar N}^{(0)}_{n,k,h(k+1)+b,j}=\cases 
-(xy)^{(h+1)(k+1)-b}\P_j(0,(h+2)(k+1)-b)\\
\quad +(xy)^{h(k+1)}\P_j(0,h(k+1)+b),\\
&\kern-3.5cm\text{if }0\le b\le k\text{ and }h(k+1)+b<n-1,\\
%-(xy)^{(h+1)(k+1)-b}\P_j(0,(h+2)(k+1)-b)\\
%\quad +(xy)^{h(k+1)}\P_j(0,h(k+1)+b),\\
%&\kern-4.5cm\text{if }k+1\le b\le 2k+1\text{ and }h(2k+2)+b<n-1,\\
%-(xy)^{(2h+3)(k+1)-b}\P_j(0,(h+2)(2k+2)-b)\\
%\quad +(xy)^{(2h+1)(k+1)}\P_j(0,h(2k+2)+b),\\
%&\kern-4.5cm\text{if }k+1\le b\le 2k+1\text{ and }h(2k+2)+b<n-1,\\
0,&\kern-4.5cm\text{if }h(k+1)=n-1.
\endcases
$$
Since $\P_a(0,b)=0$ for $a<b$, we see that the new matrix (\AQ) is
upper triangular except for the last row.

The reader should observe that, because of the permutation of the
rows, the determinant of the matrix in (\AQ) is not necessarily
equal to the determinant of $\MP_0(n,k)$, but that they rather differ
by a sign of $(-1)^{n_1\binom{k+1}2}$.

Let us concentrate on the last $k+2$ rows. There, all entries in
columns $0,1,\dots,n-k-3$ are zero. In other words, the matrix
${\bar N}_0(n,k)+C_0(n,k)$ has a block form
$$\pmatrix 
A&*\\0&B
\endpmatrix,
\tag\AQa$$
where the $(k+2)\times(k+2)$ submatrix $B$ looks as follows:
$$B=\pmatrix 
\cases
(xy)^{n-2}\P_j(0,i),&\text{for }n-k-2\le i\le n-2\\
(xy)^{n-1}\sum_{l\ge0}^{}\P_j(0,n-k+l),&\text{for }i=n-1
\endcases
\endpmatrix.
$$
Here, the index $j$ ranges over $j=n-k-2,n-k-1,\dots,n-1$.
By subtracting $xy$ times row $i$ for $i=n-k,n-k+1,\dots,n-2$
from the last row, we may transform the matrix in (\AQa) into
$$\pmatrix 
A&*\\0&B'
\endpmatrix,
\tag\AQb$$
where $B'$ is defined by
$$\pmatrix 
\cases
(xy)^{n-2}\P_j(0,i),&\text{for }n-k-2\le i\le n-2\\
0,&\text{for }i=n-1\text{ and }j\le n-2\\
(xy)^{n-1},&\text{for }i=j=n-1
\endcases
\endpmatrix,
$$
with $j$ again ranging over $j=n-k-2,n-k-1,\dots,n-1$.
These row operations do not change the value of the determinant,
and consequently the determinant of the matrix
in (\AQ) equals the one in (\AQb).

The determinant of (\AQb) is easy to compute since it is in fact
an upper triangular matrix (including the last
row!). Reading along the diagonal of this matrix, we find
$$\align 
&(xy)^k,(xy)^k,\dots,(xy)^k,\\
&(xy)^{2k+1},(xy)^{2k+1}\dots,(xy)^{2k+1},\\
&(xy)^{3k+2},(xy)^{3k+2},\dots,(xy)^{3k+2},\\
&\hbox to 5cm{\leaders \hbox{.}\hfil}\\
&(xy)^{n-2},(xy)^{n-2},
\dots,(xy)^{n-2},\\
&(xy)^{n-1},
\endalign$$
where, when arranged as above, there are exactly $k+1$ entries in each
line (except for the last line, of course). The product of these
entries is 
$$
(xy)^{n_1k(k+1)+(k+1)^2\binom {n_1}2+n-1}
=
(xy)^{(k+1)^2\binom {n_1+1}2},
$$
which, together with the earlier found sign $(-1)^{n_1\binom {k+1}2}$,
establishes the second case in (\AA).

\medskip
{\smc Case 3:} $n\equiv0$ (mod $k+1$).
Let $n=(k+1)n_1$. Here, we also depart from (\AP).
Inspection of the definition of the correction
matrix $C_0(n,k)$ in Lemma~\TD\ shows that the entries in
columns $0,1,\dots,n-k-1$ in its last row all vanish.
Therefore, if we reverse the order of rows
$s(k+1),s(k+1)+1,\dots,s(k+1)+k$, for $s=0,1,\dots,n_1-1$,
then we obtain an upper triangular matrix whose entries along
the diagonal are
$$\align 
&(xy)^k,(xy)^k,\dots,(xy)^k,\\
&(xy)^{2k+1},(xy)^{2k+1}\dots,(xy)^{2k+1},\\
&(xy)^{3k+2},(xy)^{3k+2},\dots,(xy)^{3k+2},\\
&\hbox to 5cm{\leaders \hbox{.}\hfil}\\
&(xy)^{n-1},(xy)^{n-1},
\dots,(xy)^{n-1}.
\endalign$$
The product of these entries is 
$$
(xy)^{n_1k(k+1)+(k+1)^2\binom {n_1}2}
=
(xy)^{(k+1)^2\binom {n_1+1}2-n},
$$
which, together with the sign $(-1)^{n_1\binom {k+1}2}$ that results
from the row permutation that we performed,
establishes the first case in (\AA).

\medskip
This completes the proof of the theorem.\quad \quad \qed
\enddemo

\subhead 5. Proof of Theorem \TB\endsubhead
We first use the results from Section~3 to connect the
matrix of Motzkin prefix generating functions on the left-hand side of
(\AB) to another matrix of Motzkin generating functions that appeared
in \cite{\CJKrAA}.

\proclaim{Lemma \TI}
Define matrices $M_1(n,k )$ and $\MP_1(n,k )$ by
$$
M_1(n,k ):=\left(\P^+_{i+j+1}(k ,0)\right)_{0\le i,j\le n-1}
$$
and
$$
\MP_1(n,k):=\left(\sum_{l \ge0}\P^+_{i+j+1}(k ,l )\right)_{0\le i,j\le n-1}.
$$
Then
$$
A(n)\cdot \MP_1(n,k )= M_1(n,k )+C_1(n,k ),
$$
where $A(n)$ is given by {\rm(\AC)} and the ``correction matrix" 
$C_1(n,k):=(C^{(1)}_{n,k,i,j})_{0\le i,j\le n-1}$ is defined via
$$
C^{(1)}_{n,k,i,j}=\cases 
0,&\text{if }i\le n-2,\\
(xy)^{n-1}\sum_{l\ge0}^{}\P_{j+1}(0,n-k+l),&\text{if }i=n-1.
\endcases
$$
\endproclaim

\demo{Proof}
This is a direct consequence of Lemmas~\TE\ and \TG.\quad \quad \qed
\enddemo

We are now in the position to prove Theorem~\TB.

\demo{Proof of Theorem \TB}
We start with 
$$
\MP_1(n,k)=\left(\sum_{l \ge0}\P^+_{i+j+1}(k ,l )\right).
$$
We multiply on the left by $A(n)$. According to Lemma~\TI, we get
$$\multline
M_1(n,k)+C_1(n,k)\\
=\left(
\cases 
\P^+_{i+j+1}(k ,0),&\text{for }i\le n-2\\
\P^+_{n+j}(k ,0)+
(xy)^{n-1}\sum_{l\ge0}^{}\P_{j+1}(0,n-k+l),&\text{for }i=n-1
\endcases
\right).
\endmultline$$
Since $\det A(n)=1$ by Lemma~\TC, the determinant of $M_1(n,k)+C_1(n,k)$
is the same as the determinant of $\MP_1(n,k)$. If we apply relation
(\AAa) with $l=0$, then we see that we have transformed our problem
into the problem of evaluation of the determinant of
$$(xy)^kM'_1(n,k)+C_1(n,k),$$ 
where
$$
M'_1(n,k ):=\left(\P^+_{i+j+1}(0,k)\right)_{0\le i,j\le n-1}.
$$

The determinant of $M'_1(n,k)$ is evaluated in Theorem~\TCKB. 
We now follow the arguments
of the proof of this evaluation in \cite{\CJKrAA} for a while. 

\medskip
Similarly to the proof of Theorem~\TA\ in the previous section, 
the first step consists
in the use of the decomposition (\ALa) in order to convert our problem
into the problem of the evaluation of the determinant of the matrix
$$
(xy)^k\left(\P^+_{j+1}(i,k)\right)_{0\le i,j\le n-1} + C_1(n,k).
$$

\medskip
In the second step, we rewrite $\P^+_{j+1}(i,k)$ by using (\AG).
In this manner, we arrive at the problem of evaluating the determinant
of
$$
(xy)^k\left(
\P_{j+1}(i,k)-(xy)^{i+1}\P_{j+1}(-i-2,k)\right)_{0\le i,j\le n-1} + 
C_1(n,k).
\tag\AR$$

\medskip
For the third step, we apply again the row operations given by (\AOa)
and (\AOb). 
By the computations performed in \cite{\CJKrAA, proof of
Theorem~9 with $t=1$} that finally lead to (5.14)--(5.16) there,
these operations transform the matrix in (\AR) into the matrix
$$
(xy)^kN_1(n,k)+C_1(n,k),
\tag\AS$$
where $N_1(n,k)=(N^{(1)}_{n,k,i,j})_{0\le i,j\le n-1}$ with
$$
N^{(1)}_{n,k,h(2k+2)+b,j}=\cases 
-(xy)^{(2h+1)(k+1)+b-k}\P_{j+1}(0,(h+1)(2k+2)+b-k)\\
\quad +(xy)^{2h(k+1)}\P_{j+1}(0,h(2k+2)-b+k)
&\kern-.5cm\text{if }0\le b\le k,\\
-(xy)^{(2h+1)(k+1)+b-k}\P_{j+1}(0,(h+1)(2k+2)+b-k)\\
\quad +(xy)^{(2h+1)(k+1)}\P_{j+1}(0,(h+1)(2k+2)-b+k)\\
&\kern-1.8cm\text{if }k+1\le b\le 2k+1.
\endcases
$$

\medskip
Closer inspection of the new matrix in (\AS) will lead to the claimed
result in (\AB). This is less straightforward than in the proof of
Theorem~\TA\ though.

\medskip
{\smc Case 1:} \hbox{$n\not\equiv0,1,k$ (mod $k+1$).}
Let
$$n=H(2k+2)+B
\tag\ASa$$ 
with $0\le B\le 2k+1$ but $B\ne 0,1,k,k+1,k+2,2k+1$.
Then it is not difficult to see (see the paragraphs after (5.16)
in \cite{\CJKrAA}) that, if $1\le B\le k-1$, row $H(2k+2)$ of $N_1(n,k)$
consists entirely of zeroes, while, if $k+2\le B\le 2k$, 
row $H(2k+2)+k+1$ consists entirely of zeroes. In particular, this
implies that in our case there is a row of zeroes in the
matrix in (\AS), and hence its determinant vanishes. This establishes
the fourth case on the right-hand side of (\AB).

\medskip
{\smc Case 2:} $n\equiv0$ (mod $k+1$).
With the notation of (\ASa), we have $B=0$ or $B=k+1$. 
Similarly to Case~2 in the proof of Theorem~\TA\ in the previous
section, we subtract $xy$ times row~$i$ for $i=n-k-1,n-k,\dots,n-2$
from the last row of the matrix in (\AS). Thus, we obtain the matrix
$$
(xy)^kN_1(n,k)+\bar C_1(n,k),
\tag\AT$$
where $\bar C_1(n,k)$ is equal to the zero matrix except for the
$(n-1,n-1)$-entry (the bottom-right entry), which equals $(xy)^{n-1}$.
Obviously, the determinant did not change. By using linearity in the
last column, we may write the determinant of the matrix in (\AT) as
$$
(xy)^{nk}\det N_1(n,k)+(xy)^{(n-1)k+n-1}\det \big(N_1(n,k)\big)^{n-1}_{n-1},
\tag\AU$$
where $\big(N_1(n,k)\big)^{n-1}_{n-1}$ denotes the matrix arising from
$N_1(n,k)$ by omitting the last row and the last column.

Since the matrix $N_1(n,k)$ arose from $M_1(n,k)$ by row operations
that did not change the determinant, the determinant of $N_1(n,k)$
equals the expression for $\det M_1(n,k)$ given in Theorem~\TCKB, namely
$$
(-1)^{n_1\binom
{k+1}2}(xy)^{(k+1)^2\binom {n_1}2}
\frac {y^{(k+1)(n_1+1)}-x^{(k+1)(n_1+1)}} {y^{k+1}-x^{k+1}},
$$
with $n=n_1(k+1)$. Similarly, the determinant of
$\big(N_1(n,k)\big)^{n-1}_{n-1}$ equals the expression for
$\det M_1(n-1,k)$ given in Theorem~\TCKB, namely
$$
(-1)^{n_1\binom
{k+1}2+k}(xy)^{(k+1)^2\binom {n_1-1}2+(n_1-1)k(k+1)}
\frac {y^{(k+1)n_1}-x^{(k+1)n_1}} {y^{k+1}-x^{k+1}}.
$$
If these two expressions are substituted in (\AU), then the expression
given in the first case on the right-hand side of (\AB) is obtained
after minor modification.

\medskip
{\smc Case 3:} $n\equiv1$ (mod $k+1$).
With the notation of (\ASa), we have $B=1$ or $B=k+2$. 
Without loss of generality, we may assume $k\ge1$ (cf\. Remark~(1)
after Theorem~\TB).

By reusing the arguments in Case~1, we see that, here, it is the last row
of $N_1(n,k)$ (namely row $n-1=H(2k+2)+B-1$) which consists entirely of zeroes.
Since, in (\AS), the matrix $C_1(n,k)$ --- which is a matrix with
potentially non-zero entries in the last row --- is added to
$(xy)^kN_1(n,k)$, we cannot conclude that the determinant of (\AS)
vanishes, but rather further work is required.

From now on, let $n=(k+1)n_1+1$.
As in Case~2 of the proof of Theorem~\TA\ in the previous section, 
we reverse the order of rows
$s(k+1),s(k+1)+1,\dots,s(k+1)+k$, for $s=0,1,\dots,n_1-1$,
leaving the last row, row $n-1=(k+1)n_1$, in place. 
Furthermore, we factor
$(xy)^{h(k+1)}$ from all the entries in rows
$h(k+1),h(k+1)+1,\dots,h(k+1)+k$, $h=0,1,\dots,n_1-1$. This yields an
overall factor of
$$
(xy)^{(k+1)^2\binom {n_1}2}
\tag\AV
$$ 
by which we have to multiply the determinant of the remaining matrix
in the end. We must as well multiply by the sign
$$(-1)^{n_1\binom {k+1}2}
\tag\AW$$ 
in order to take into account the permutation of the rows that we
performed.

In this manner, we arrive at the matrix
$$
{\bar N}_1(n,k)+C_1(n,k),
\tag\AX$$
where the matrix 
${\bar N}_1(n,k)=({\bar N}^{(1)}_{n,k,i,j})_{0\le i,j\le n-1}$ is given
by
$$
{\bar N}^{(1)}_{n,k,h(k+1)+b,j}=\cases 
\P_{j+1}(0,h(k+1)+b)
-
(xy)^{k-b+1}\P_{j+1}(0,(h+2)(k+1)-b)\\
&\kern-5.5cm\text{if }0\le b\le k\text{ and }h(k+1)+b<n-1,\\
0,&\kern-5.5cm\text{if }h(k+1)=n-1.
\endcases
\tag\AY
$$
We should observe that, for $1\le i\le n-2$, the first
non-zero entry in row $i$ (which is to be found in column $i-1$)
equals $1$.

In the matrix in (\AX), we replace the $0$-th row by
$$\sum _{h=0} ^{n_1-1}\sum _{b=0} ^{k} (-1)^{h(k+1)+b}
\sum _{s=0} ^{h}c(h,b,s)\,x^{s(k+1)}y^{(h-s)(k+1)}
\cdot
\big(\text {row $(h(k+1)+b)$}\big),\tag\BA$$
where the coefficients $c(h,b,s)$ are given by
$$
c(h,b,s)=\cases 
x^b+y^b,
&\text{if }b\ne0,\\
1,
&\text{if }b=0.
\endcases
$$
Since the coefficient of the $0$-th row in the linear combination
(\BA) is $1$, this does not change the value of the determinant.
It should be noted that the last row, row $n-1=n_1(k+1)$, is not
involved in the linear combination (\BA).

Now we have to redo the computation in \cite{\CJKrAA} with $t=1$, starting with
(5.21) and leading to the result in the display in the centre
of p.~161 there, however with the relaxed condition that $j+1\le
n=n_1(k+1)+1$. (In \cite{\CJKrAA} we have $j+1\le n_1(k+1)$ at that
point.)
Taking also into account our assumption that $k\ge1$, the final result
is that the $(0,j)$-entry in the new matrix is given by
$$\multline
-(-1)^{n_1(k+1)}\sum _{s=0} ^{n_1} 
x^{s(k+1)}y^{(n_1-s)(k+1)}
\big(\P_{j+1}(0,n_1(k+1))-(x+y)\P_{j+1}(0,n_1(k+1)+1)\big).
\endmultline
%\tag\BI
$$
It is important to observe that, again because $\P_a(0,b)=0$ for $a<b$, 
this expression vanishes for
$j<n_1(k+1)-1$, so that only the right-most two entries in row~$0$ are
non-zero. 

In addition to the above modifications of the $0$-th row, in analogy
to similar operations in Case~2 in the proof of Theorem~\TA\ and in
Case~2 of the current proof, we also replace the last row,
row~$n-1=n_1(k+1)$, by
$$
\big(\text {row $n-1$}\big)-
(xy)^{n-1}\sum _{i=n-k} ^{n-2}
\big(\text {row $i$}\big).
\tag\BJ$$
When doing this operation it is important to observe that all entries
in column~$j$ with $j\le n-k-1$ are actually zero, that the entries 
${\bar N}^{(1)}_{n,k,i,j}$ given by (\AY) for $i\le n-2$ and 
$j\ge n-k$ are given by
$\P_{j+1}(0,i)$, except for the $(n-2,n-1)$-entry, which is
equal to $\P_{n}(0,n-2)
-
(xy)\P_{n}(0,n)$.
Again, this operation does not change the determinant.

Altogether, the new matrix obtained is
${\widetilde M}_1(n,k)=({\widetilde M}^{(1)}_{n,k,i,j})_{0\le i,j\le
n-1}$, where
$$
{\widetilde M}^{(1)}_{n,k,h(k+1)+b,j}=\cases 
(-1)^{n_1(k+1)+1}\sum _{s=0} ^{n_1} 
x^{s(k+1)}y^{(n_1-s)(k+1)}\\
\quad \times
\big(\P_{j+1}(0,n_1(k+1))-(x+y)\P_{j+1}(0,n_1(k+1)+1)\big),\\
&\kern-6cm\text{if }h=b=0,\\
\P_{j+1}(0,h(k+1)+b)
-
(xy)^{k-b+1}\P_{j+1}(0,(h+2)(k+1)-b),\\
&\kern-6cm\text{if }0\le b\le k\text{ and }0<h(k+1)+b<n-1,\\
(xy)^{n-1}\big(\P_{j+1}(0,n_1(k+1))+(1+xy)\P_{j+1}(0,n_1(k+1)+1)\big),\\
&\kern-6cm\text{if }h(k+1)=n-1.
\endcases
$$
The determinant of this matrix is the same as that of the matrix
in (\AX). 

%Achtung: $(n_1(k+1)-1,j)$-Eintragung ($h=n_1-1$, $b=k$)
%von $C_1(n,k)$ involviert fuer $j=n_1(k+1)$ den zweiten
%Term!! Das erklaert den $-xy\P_{n_1(k+1)+1}(0,n_1(k+1)+1)$ Term.
It is helpful to display the schematic form of this matrix:
$$
{\widetilde M}_1(n,k)=\pmatrix 
0&\dots&0&a&b\\
&M&&\vdots&\vdots\\
0&\dots&0&c&d\\
\endpmatrix,
\tag\BK$$
where 
$$\align 
a&=(-1)^{n_1(k+1)+1}\frac {x^{(k+1)(n_1+1)}-y^{(k+1)(n_1+1)}} {x^{k+1}-y^{k+1}},\\
b&=(-1)^{n_1(k+1)+1}\frac {x^{(k+1)(n_1+1)}-y^{(k+1)(n_1+1)}} {x^{k+1}-y^{k+1}}
\big((n_1(k+1)+1)(x+y)-(x+y)\big),\\
c&=(xy)^{n-1},\\
d&=(xy)^{n-1}(n_1(k+1)+1)(x+y)+(1+xy),
\tag\BL
\endalign$$
and $M=(M_{i,j})_{1\le i\le n_1(k+1)-1,\,0\le j\le n_1(k+1)-2}$ with
$$
M_{h(k+1)+b,j}=\P_{j+1}(0,h(k+1)+b)
-
(xy)^{k-b+1}\P_{j+1}(0,(h+2)(k+1)-b)
$$
for $0\le b\le k$.

In order to evaluate the determinant of ${\widetilde M}_1(n,k)$, we do
a Laplace expansion simultaneously with respect to the top and the
bottom row. Thereby, we obtain
$$
\det {\widetilde M}_1(n,k)=(-1)^n\det\pmatrix a&b\\c&d\endpmatrix
\cdot \det M.
\tag\BM$$
Straightforward calculation shows that
$$
\det\pmatrix a&b\\c&d\endpmatrix
=
(-1)^n(xy)^{n-1}\frac {x^{(k+1)(n_1+1)}-y^{(k+1)(n_1+1)}} {x^{k+1}-y^{k+1}}
(1+x)(1+y),
$$
while $M$ is an upper triangular matrix so that its determinant is
equal to the product of its diagonal entries, all of which are~$1$. 
If everything is put
together with the earlier obtained factors (\AV) and (\AW), then
we arrive at the expression given in the second case on the right-hand
side of (\AB).

\medskip
{\smc Case 4:} $n\equiv k$ (mod $k+1$).
With the notation of (\ASa), we have $B=k$ or $B=2k+1$. 
Again, without loss of generality, we may assume $k\ge1$ (cf\. Remark~(1)
after Theorem~\TB).

Let $n=n_1(k+1)+k$. We do the same row operations as in Case~3, except
the one in (\BJ). This produces a matrix of the form
$$\pmatrix 
0&\dots&0&a&b&\dots\\\\
&\hbox{\seventeenpoint $M$}&&&\hbox{\seventeenpoint $*$}\\\\
&\hbox{\seventeenpoint $0$}&&&&\hbox{\seventeenpoint $D$}\\\\
%0&\dots&0&0&0&\dots\\
0&\dots&0&c&d&\dots
\endpmatrix,
\tag\BN$$
where $a,b,c,d,M$ are as in (\BK) and (\BL) (with the meaning of $n$ in the
definitions of $c$ and $d$ being the current one), 
and $D$ is a $(k-1)\times (k-1)$ ``reflected upper triangular" matrix. 
(By ``reflected upper triangular" we mean a matrix where all entries above 
the {\it anti}-diagonal of the matrix are equal to $0$.)
The entries $a$ and $c$ are located in column\break $n_1(k+1)-1$, so that the
submatrix~$M$ is located strictly to the left of this column while the
submatrix $D$ is located strictly to the right of column~$n_1(k+1)$
(the indexing of rows and columns starting at~$0$ as usual).
To the left of $D$ --- in the rows covered by~$D$ --- there are only zeroes.

By performing a Laplace expansion simultaneously with respect to the
top and the bottom row, one sees that the determinant of the above
matrix equals
$$
(-1)^{n}\det\pmatrix a&b\\c&d\endpmatrix
\cdot\det M\cdot \det D.
\tag\BO$$
Comparison with (\BM) shows that
 the determinant of the matrix in (\BN) differs from 
$\det\widetilde M(n,k)$ (with $\widetilde M(n,k)$ given in (\BK))
by a factor of
$$
(-1)^{k-1}(xy)^{k-1}\det D.
$$
Here, the factor of $(-1)^{k-1}$ comes from the factor $(-1)^n$ in
(\BO), taking into account that our current $n$ is by $k-1$ larger than the
$n$ in Case~3, and the factor $(xy)^{k-1}$ comes from the factor
$(xy)^{n-1}$ in the definitions of $c$ and~$d$, again taking into
account that the $n$ here differs from the one in Case~3.

Since the entries of $D$ were not affected by the row operations
from Case~3, they still equal the corresponding entries in
$(xy)^kN_1(n,k)$ (cf\. (\AS)). Consequently --- as we already stated
earlier --- $D$ is ``reflected upper triangular", with entries
$(xy)^n$ along the main antidiagonal. It follows that $\det D$
equals $(-1)^{\binom {k-1}2}(xy)^{n(k-1)}$. Hence, 
the determinant of the matrix in (\BN) differs from 
$\det\widetilde M(n,k)$ by a factor of
$$
(-1)^{k-1+\binom {k-1}2}(xy)^{(n+1)(k-1)}=
(-1)^{\binom k2}(xy)^{(n_1+1)(k^2-1)}.
$$
This is indeed exactly the factor by which the third expression on the
right-hand side of (\AB) differs from the second expression. 

\medskip
This completes the proof of the theorem.\quad \quad \qed
\enddemo

\subhead 6. Specialisations\endsubhead
In this section we list specialisations of Theorems~\TA\ and~\TB.
The special values of $x$ and $y$ that we choose 
are those that we discussed at the end of Section~2.
In all the results that we list in this section, the convention 
of Remark~(1) after the statements of Theorems~\TA\ and \TB\ applies
(in a slightly modified form):
for $k=0,1$, it is the first {\it applicable} case in the case distinctions on
the right-hand sides that produces the correct result.

\medskip
We begin by setting $x=-y=\sqrt{-1}$ in Theorem~\TA. 
Using (\SM), we obtain the following result.

\proclaim{Corollary \TJ}
For all positive integers $n$ and non-negative integers $k$, we have
$$
\det_{0\le i,j\le n-1}\(\sum_{l=0}^k\binom {i+j}{\fl{\frac {1} {2}(i+j+1-k)}+l}\)=
\cases 
(-1)^{n_1\binom {k +1}2},
&\text{if }n=(k +1)n_1,\\
(-1)^{n_1\binom {k +1}2},
&\text{if }n=(k +1)n_1+1,\\
0,&\text{if }n\not\equiv0,1~(\text{\rm mod }k +1).\vphantom{\hbox{$\Big($}}
\endcases
\tag\CA$$
\endproclaim

A noteworthy special case is the one for $k=0$,
$$
\det_{0\le i,j\le n-1}\(\binom {i+j}{\fl{\frac {1} {2}(i+j+1)}}\)
=\det_{0\le i,j\le n-1}\(\binom {i+j}{\fl{\frac {1} {2}(i+j)}}\)
=1.
\tag\CAa$$
In other words, this gives the ``Hankel transform" of the sequence
$\(\binom n{\fl{n/2}}\)_{n\ge0}$ of central and ``almost central" 
binomial coefficient. According to \cite{\OEIS, Sequence~{\tt A001405}},
this Hankel determinant evaluation had been observed by
Philippe Del\'eham in 2007.

On the other hand, for $k=1$ we get
$$
\det_{0\le i,j\le n-1}\(\binom {i+j+1}{\fl{\frac {1} {2}(i+j+2)}}\)=
\det_{0\le i,j\le n-1}\(\binom {i+j+1}{\fl{\frac {1} {2}(i+j+1)}}\)=
(-1)^{\fl{n/2}},
\tag\CAb$$
thus obtaining the ``Hankel transform" of the shifted sequence
$\(\binom n{\fl{n/2}}\)_{n\ge1}$ of central and ``almost central"  
binomial coefficients. We add that the choices of $k=2$ and $k=3$
provide Hankel determinant evaluations for the sequences
{\tt A026010} and {\tt A026023} in \cite{\OEIS}.

Next we set $x=-y=\sqrt{-1}$ in Theorem~\TB, upon using (\SM) again. 
This leads to the following determinant identity.

%$$
%\sum_{r=0}^{n_1}I^{(k+1)r}(-I)^{(k+1)(n_1-r)}
%=\cases 
%(-1)^{(k+1)n_1/2}(n_1+1),
%&\text{if $k$ is odd,}\\
%0,
%&\text{if $k$ is even and $n_1$ is odd.}\\
%(-1)^{(k+1)n_1/2},
%&\text{if $k$ and $n_1$ are even.}\\
%\endcases
%$$

\proclaim{Corollary \TK}
For all positive integers $n$ and non-negative integers $k$, we have
$$\multline
\det_{0\le i,j\le n-1}\(\sum_{l=0}^k\binom {i+j+1}{\fl{\frac {1} {2}(i+j+2-k)}+l}\)\\=
\cases 
2n_1+1,
&\text{if }n=(k +1)n_1\text{ and }k\equiv1~(\text{\rm mod }4),\\
1,
&\text{if }n=(k +1)n_1\text{ and }k\equiv3~(\text{\rm mod }4),\\
       (-1)^{n_1/2},
&\text{if }n=(k +1)n_1,\text{ and $k$ and $n_1$ are even},\\
(-1)^{(k+n_1-1)/2},
&\text{if }n=(k +1)n_1,\text{ $k$ is even, and $n_1$ is odd},\\
2n_1+2,
&\text{if }n=(k +1)n_1+1,\text{ and $k$ is odd},\\
2(-1)^{n_1/2},
&\text{if }n=(k +1)n_1+1,\text{ and $k$ and $n_1$ are even},\\
       (-1)^{(k-1)/2}(2n_1+2),
&\text{if }n=(k +1)n_1+k,\text{ and $k$ is odd},\\
2(-1)^{(k+n_1)/2},
&\text{if }n=(k +1)n_1+k,\text{ and $k$ and $n_1$ are even},\\
0,&\text{otherwise.}
\endcases\\
\endmultline
\tag\CB
$$
\endproclaim

For $k=0$, we obtain (\CAb) again, while for $k=1$ we get
$$
\det_{0\le i,j\le n-1}\(\binom {i+j+2}{\fl{\frac {1} {2}(i+j+3)}}\)=
\det_{0\le i,j\le n-1}\(\binom {i+j+2}{\fl{\frac {1} {2}(i+j+2)}}\)=
n+1,
\tag\CAc
$$
which is the ``Hankel transform" of the doubly shifted sequence
$\(\binom n{\fl{n/2}}\)_{n\ge2}$ of central and ``almost central"  
binomial coefficients. Clearly, for $k=2$ and $k=3$, Corollary~\TK\
provides Hankel determinant evaluations for the sequences
{\tt A026010} and {\tt A026023}, respectively, with the first
element of each sequence omitted.

\medskip
We continue setting $x=y^{-1}=\om$ in Theorem~\TA. 
We recall that this specialisation in $\P^+_n(k ,l )$ corresponds to
weighting each path by~$1$ --- which amounts to ordinary counting
of paths --- so that $\P^+_n(k,l)_{x=y^{-1}=\om}$ is simply equal
to the number of all three-step paths from $(0,k)$ to $(n,l)$ that
never run below the $x$-axis. Consequently,
$$
\sum_{l\ge0}\P^+_n(k,l)\Big\vert_{x=y^{-1}=\om}
\tag\CAd$$
equals the number of three-step paths starting at $(0,k)$, proceeding for 
$n$ steps, and never running below the $x$-axis. For $k=0$, these are
the Motzkin prefix numbers $\MP_n$ mentioned in the introduction.
For generic~$k$, these numbers can be considered as {\it generalised
Motzkin prefix numbers}, for which (\SN) provides an explicit formula.
We denote the number in (\CAd) by $\MP_n(k)$.

\proclaim{Corollary \TL}
For all positive integers $n$ and non-negative integers $k$, we have
$$\align
\det_{0\le i,j\le n-1}\big(\MP_{i+j}(k)\big)
&=
\det_{0\le i,j\le n-1}\(\sum _{\ell\ge0} ^{}
\sum_{l=-k}^{k+1}\binom {i+j} {\ell,\ell+l}\)\\
&=
\cases 
(-1)^{n_1\binom {k +1}2},
&\text{if }n=(k +1)n_1,\\
(-1)^{n_1\binom {k +1}2},
&\text{if }n=(k +1)n_1+1,\\
0,&\text{if }n\not\equiv0,1~(\text{\rm mod }k +1).
\endcases
\tag\CC
\endalign$$
\endproclaim

Clearly, the case $k=0$ provides the proof of (\CKc).
Further noteworthy special cases are the one for $k=1$,
$$
\det_{0\le i,j\le n-1}\big(\MP_{i+j}(1)\big)
=
(-1)^{\fl{n/2}},
\tag\CAe
$$
providing the ``Hankel transform" of Sequence {\tt A025566} in 
\cite{\OEIS}, and the one for $k=2$,
$$
\det_{0\le i,j\le n-1}\big(\MP_{i+j}(2)\big)
=
\cases 
(-1)^{\fl{n/3}},
&\text{if }n\equiv0,1~(\text{\rm mod }3),\\
0,&\text{if }n\equiv2~(\text{\rm mod }3),
\endcases
\tag\CAf
$$
providing the ``Hankel transform" of Sequence {\tt A005774}
in \cite{\OEIS}.

Specialisation of $x=y^{-1}=\om$ in Theorem~\TB\ yields
further Hankel determinant evaluations for (generalised)
Motzkin prefix numbers.

%$$\multline
%\sum_{r=0}^{n_1}\om^{(k+1)r}\om^{-(k+1)(n_1-r)}
%=
%\sum_{r=0}^{n_1}\om^{(k+1)r}\om^{-(k+1)(n_1-r)}
%\\
%=\cases 
%(-1)^{(k+1)n_1}(n_1+1),
%&\text{if }k\equiv2~\text{(mod 3),}\\
%(-1)^{(k+1)n_1},
%&\text{if }k\not\equiv2~\text{(mod 3) and }n_1\equiv0~\text{(mod 3),}\\
%(-1)^{(k+1)n_1+1},
%&\text{if }k\not\equiv2~\text{(mod 3) and }n_1\equiv1~\text{(mod 3),}\\
%0,&\text{otherwise.}
%\endcases
%\endmultline$$

\proclaim{Corollary \TM}
For all positive integers $n$ and non-negative integers $k$, we have
$$\multline
\det_{0\le i,j\le n-1}\big(\MP_{i+j+1}(k)\big)
=
\det_{0\le i,j\le n-1}\(\sum _{\ell\ge0} ^{}
\sum_{l=-k}^{k+1}\binom {i+j+1} {\ell,\ell+l}\)
\\
=
\cases 
(-1)^{n_1\binom {k +2}2},
&\text{if }n=(k +1)n_1\text{ and }k\equiv2~(\text{\rm mod }3),\\
(-1)^{n_1\binom {k +2}2},
&\text{if }n=(3k +3)n_1\text{ and }k\not\equiv2~(\text{\rm mod }3),\\
2(-1)^{(n_1+1)\binom {k +2}2+1},
&\text{if }n=(3k +3)n_1+k+1\text{ and }k\not\equiv2~(\text{\rm mod }3),\\
(-1)^{n_1\binom {k +2}2},
&\text{if }n=(3k +3)n_1+2k+2\text{ and }k\not\equiv2~(\text{\rm mod }3),\\
3(-1)^{n_1\binom {k +2}2}(n_1+1),
&\text{if }n=(k +1)n_1+1\text{ and }k\equiv2~(\text{\rm mod }3),\\
3(-1)^{n_1\binom {k +2}2},
&\text{if }n=(3k +3)n_1+1\text{ and }k\not\equiv2~(\text{\rm mod }3),\\
3(-1)^{(n_1+1)\binom {k +2}2+1},
&\text{if }n=(3k +3)n_1+k+2\text{ and }k\not\equiv2~(\text{\rm mod }3),\\
3(-1)^{(n_1+1)\binom {k +2}2+1}(n_1+1),
&\text{if }n=(k +1)n_1+k\text{ and }k\equiv2~(\text{\rm mod }3),\\
3(-1)^{(n_1+1)\binom {k +2}2+1},
&\text{if }n=(3k +3)n_1+k\text{ and }k\not\equiv2~(\text{\rm mod }3),\\
3(-1)^{(n_1+1)\binom {k +2}2},
&\text{if }n=(3k +3)n_1+2k+1\text{ and }k\not\equiv2~(\text{\rm mod }3),\\
0,&\text{otherwise}.
\endcases\\
\endmultline
\tag\CD
$$
\endproclaim

The special cases $k=0,1,2$ are explicitly
$$\align
\det_{0\le i,j\le n-1}\big(\MP_{i+j+1}\big)
&=
\cases 
(-1)^{\fl{n/3}},
&\text{if }n\equiv0,2~(\text{mod }3),\\
2(-1)^{\fl{n/3}},
&\text{if }n\equiv1~(\text{mod }3),
\endcases
\tag\CAg\\
\det_{0\le i,j\le n-1}\big(\MP_{i+j+1}(1)\big)
&=
\cases 
(-1)^{\fl{n/6}},
&\text{if }n\equiv0,4~(\text{mod }3),\\
3(-1)^{\fl{n/6}},
&\text{if }n\equiv1,3~(\text{mod }3),\\
2(-1)^{\fl{n/6}},
&\text{if }n\equiv2~(\text{mod }3),\\
0,
&\text{if }n\equiv5~(\text{mod }3),
\endcases
\tag\CAh\\
\det_{0\le i,j\le n-1}\big(\MP_{i+j+1}(2)\big)
&=
\cases 
1,
&\text{if }n\equiv0~(\text{mod }3),\\
3\cl{n/3},
&\text{if }n\equiv1~(\text{mod }3),\\
-3\cl{n/3},
&\text{if }n\equiv2~(\text{mod }3),
\endcases
\tag\CAi
\endalign$$
providing further Hankel determinant evaluations for the sequences
{\tt A005773}, {\tt A025566}, and {\tt A005774} in \cite{\OEIS}.

\medskip
Next we turn our attention to the specialisation $x=y=1$.
Use of (\SO) in Theorem~\TA\ yields the following result.

\proclaim{Corollary \TN}
For all positive integers $n$ and non-negative integers $k$, we have
$$\det_{0\le i,j\le n-1}\(\sum_{l=-k}^{k+1}\binom {2i+2j} {i+j+l}\)=
\cases 
(-1)^{n_1\binom {k +1}2},
&\text{if }n=(k +1)n_1,\\
(-1)^{n_1\binom {k +1}2},
&\text{if }n=(k +1)n_1+1,\\
0,&\text{if }n\not\equiv0,1~(\text{\rm mod }k +1).
\endcases
\tag\CE$$
\endproclaim

For $k=0$, Equation~(\CE) says that the ``Hankel transform" of the 
sequence\break $\(\binom {2n+1}{n+1}\)_{n\ge0}$ (which is \cite{\OEIS,
Sequence~{\tt A001700}}) is the all-1 sequence. 
This is a well-known result, and it is also covered by 
\cite{\CJKrAA, Theorem~21}. 

For $k=1$, Equation~(\CE) provides the ``Hankel transform" of the
sequence $\(\binom {2n+2}{n}\)_{n\ge0}$ (which is \cite{\OEIS,
Sequence~{\tt A001791}} up to a shift). 
Again, this is a known result, see e.g.~\cite{\CJKrAA, Cor.~20 with $k=1$}.

On the other hand, specialising $x=y=1$ in Theorem~\TB, we arrive
at the following Hankel determinant evaluation.

\proclaim{Corollary \TO}
For all positive integers $n$ and non-negative integers $k$, we have
$$\multline
\det_{0\le i,j\le n-1}\(\sum_{l=-k}^{k+1}\binom {2i+2j+2} {i+j+l+1}\)\\
=
\cases 
(-1)^{n_1\binom {k +1}2}(2n_1+1),
&\text{if }n=(k +1)n_1\text{ and $k$ is even,}\\
(-1)^{n_1\binom {k +1}2},
&\text{if }n=(k +1)n_1\text{ and $k$ is odd,}\\
      (-1)^{n_1\binom {k+1}2}(4n_1+4),
&\text{if }n=(k +1)n_1+1,\\
      (-1)^{(n_1+1)\binom {k+1}2+k}(4n_1+4)
&\text{if }n=(k +1)n_1+k,\\
0,&\text{if }n\not\equiv0,1,k~(\text{\rm mod }k +1).
\endcases\\
\endmultline
\tag\CF
$$
\endproclaim

Similarly to before, for $k=1$ this recovers 
\cite{\CJKrAA, Cor.~23 with $k=1$}, while for $k=0$ it proves
Conjecture~24 in \cite{\CJKrAA} for $k=0$ and $k=1$.

\medskip
Finally, we set $x=y=-1$ in Theorem~\TA. 
Using (\SR), we obtain a determinant evaluation which can be
considered to be in a row with \cite{\CJKrAA, Cors.~12, 15, 13, 18, 16}.

\proclaim{Theorem \TP}
For all positive integers $n$ and non-negative integers $k$, we have
$$\multline
\det_{0\le i,j\le n-1}\(\frac {2k+2} {i+j+k+1}
     \binom {2i+2j-1} {i+j+k}\)\\
=
\cases 
(-1)^{n_1\binom {k +1}2+1}(n_1-1),
&\text{if }n=(k +1)n_1\\
(-1)^{n_1\binom {k +1}2+k+1}n_1,
&\text{if }n=(k +1)n_1+1\\
0,&\text{if }n\not\equiv0,1~(\text{\rm mod }k +1).
\endcases
\endmultline
\tag\CG
$$
Here, the $(0,0)$-entry of the matrix on the left-hand side is zero by
definition.
\endproclaim

\demo{Proof}
Since
$$
\frac {2k+2} {i+j+k+1}
     \binom {2i+2j-1} {i+j+k}=\frac {k+1} {i+j}
     \binom {2i+2j} {i+j+k+1},
$$
the matrix on the left-hand side of (\CG), of which the determinant is taken,
is the $n\times n$ Hankel matrix corresponding to the sequence in (\SR),
up to the sign in (\SR), and up to the convention of how to interpret the
term for $n=0$ in (\SR). Applying the specialisation $x=y=-1$ in
Theorem~\TA\ and using (\SR), we are led to the determinant evaluation
$$\det_{0\le i,j\le n-1}\((-1)^{i+j+k}\frac {k+1} {i+j}
     \binom {2i+2j} {i+j+k+1}\)=
\cases 
(-1)^{n_1\binom {k +1}2},
&\text{if }n=(k +1)n_1,\\
(-1)^{n_1\binom {k +1}2},
&\text{if }n=(k +1)n_1+1,\\
0,&\text{if }n\not\equiv0,1~(\text{\rm mod }k +1),\vphantom{\hbox{$\Big($}}
\endcases
\tag\CH
$$
where the $(0,0)$-entry in the matrix on the left-hand side has to be
taken as~$1$. Let $M=(M_{i,j})_{0\le i,j\le n-1}$ be that matrix.
We write the $0$-th row of the matrix as
$$
(0,M_{0,1},M_{0,2},\dots,M_{0,n-1})+(1,0,\dots,0).
$$
Subsequently, we use the linearity of the determinant in this row to decompose the 
determinant into the sum
$$\multline
\det_{0\le i,j\le n-1}\((-1)^{i+j+k}\frac {k+1} {i+j}
     \binom {2i+2j} {i+j+k+1}\)\\+
\det_{1\le i,j\le n-1}\((-1)^{i+j+k}\frac {k+1} {i+j}
     \binom {2i+2j} {i+j+k+1}\),
\endmultline\tag\CI$$
where the $(0,0)$-entry in the first matrix is zero by definition.
The first determinant in (\CI) is thus the determinant in (\CG), up to a sign
of $(-1)^{nk}$. Because of (\CH), we know the total value of (\CI),
while the second determinant in (\CI) has been evaluated in 
\cite{\CJKrAA, Corollary~18}.\footnote{Unfortunately, the determinant
evaluation in Corollary~18 of \cite{\CJKrAA} is seriously mistyped:
on the right-hand side, every occurrence of~$k$ must be replaced 
by~$k-1$. Alternatively, we might replace every occurrence of~$k$ on the
left-hand side by~$k+1$ (and, in the sentence above Corollary~18 in
\cite{\CJKrAA} disregard ``and replacing $k$ by $k-1$".) The latter
leads in fact to the form of this evaluation that we need here.} 
Thus, this sets up an equation for the determinant in (\CG), 
which we just have to solve.\quad \quad \qed
\enddemo

Specialising $x=y=1$ in Theorem~\TB\ and using (\SO) again, we recover  
Corollary~15 in~\cite{\CJKrAA}.

\subhead Acknowledgements\endsubhead
The authors thank an anonymous referee for an extremely careful reading
of the original manuscript.

\enddocument